\documentclass[11pt]{article}

\title{Gluing pseudoholomorphic curves along branched covered
 cylinders I}
 \author{Michael Hutchings and Clifford Henry Taubes}
\date{}

\usepackage{amssymb}
\usepackage{latexsym}
\usepackage{amsmath}
\usepackage{amsthm}
\usepackage{amscd}
\usepackage[dvips]{graphics}

\numberwithin{equation}{section}

\newcommand{\mc}[1]{{\mathcal #1}}

\newtheorem{theorem}{Theorem}[section]
\newtheorem{proposition}[theorem]{Proposition}
\newtheorem{corollary}[theorem]{Corollary}
\newtheorem{lemma}[theorem]{Lemma}
\newtheorem{lemma-definition}[theorem]{Lemma-Definition}

\theoremstyle{definition}
\newtheorem{definition}[theorem]{Definition}
\newtheorem{remark}[theorem]{Remark}

\newtheorem{example}[theorem]{Example}

\newtheorem{notation}[theorem]{Notation}

\newtheorem{assumption}[theorem]{Assumption}

\newcommand{\floor}[1]{\left\lfloor #1 \right\rfloor}
\newcommand{\ceil}[1]{\left\lceil #1 \right\rceil}

\newcommand{\pin}{P^{\op{in}}}
\newcommand{\pout}{P^{\op{out}}}
\newcommand{\lin}{\Lambda^{\op{in}}}
\newcommand{\lout}{\Lambda^{\op{out}}}

\newcommand{\eqdef}{\;{:=}\;}
\newcommand{\fedqe}{\;{=:}\;}

\newcommand{\C}{{\mathbb C}}
\newcommand{\Q}{{\mathbb Q}}
\newcommand{\R}{{\mathbb R}}

\newcommand{\Z}{{\mathbb Z}}
\newcommand{\op}{\operatorname}
\newcommand{\dbar}{\overline{\partial}}
\newcommand{\zbar}{\overline{z}}

\newcommand{\Hom}{\op{Hom}}
\newcommand{\Ker}{\op{Ker}}
\newcommand{\Coker}{\op{Coker}}

\newcommand{\tensor}{\otimes}

\newcommand{\vu}{\nu}

\newcommand{\bpm}{\begin{pmatrix}}
\newcommand{\epm}{\end{pmatrix}}

\begin{document}

\maketitle

\begin{abstract}
This paper and its sequel prove a generalization of the usual gluing
theorem for two index $1$ pseudoholomorphic curves $u_+$ and $u_-$ in
the symplectization of a contact 3-manifold.  We assume that for each
embedded Reeb orbit $\gamma$, the total multiplicity of the negative
ends of $u_+$ at covers of $\gamma$ agrees with the total multiplicity
of the positive ends of $u_-$ at covers of $\gamma$.  However, unlike
in the usual gluing story, here the individual multiplicities are
allowed to differ.  In this situation, one can often glue $u_+$ and
$u_-$ to an index $2$ curve by inserting genus zero branched covers of
$\R$-invariant cylinders between them.  We establish a combinatorial
formula for the signed count of such gluings.  As an application, we
deduce that the differential $\partial$ in embedded contact homology
satisfies $\partial^2=0$.

This paper explains the more algebraic aspects of the story, and proves the
above formulas using some analytical results from part II.
\end{abstract}

\tableofcontents

\section{Statement of the gluing theorem}
\label{sec:SGT}

\subsection{Pseudoholomorphic curves in symplectizations}
\label{sec:GP}

Our gluing theorem concerns pseudoholomorphic curves in the
symplectization of a contact 3-manifold.  We now recall some mostly
standard definitions and introduce some notation regarding such curves.
The geometric setup here is essentially that of Hofer, Wysocki, and Zehnder
\cite{hofer}, and the four-dimensional case of the setup used to
define symplectic field theory \cite{egh}.

Let $Y$ be a closed oriented 3-manifold.  Let $\lambda$ be a contact
form on $Y$, i.e.\ a $1$-form $\lambda$ such that $\lambda\wedge
d\lambda >0$.  The associated contact structure is the 2-plane field
$\xi\eqdef\Ker(\lambda)$.  The contact form $\lambda$ determines a
vector field $R$ on $Y$, called the Reeb vector field, which is
characterized by $d\lambda(R,\cdot)=0$ and $\lambda(R)=1$.  A {\em
Reeb orbit\/} is a closed orbit of the flow $R$, i.e.\ a map
$\gamma:\R/T\Z\to Y$ for some $T\in\R$, modulo reparametrization, such
that $\partial_t\gamma(t)=R(\gamma(t))$.  We do not require $\gamma$
to be an embedding.  Define the homology class of the Reeb orbit by
$[\gamma]\eqdef \gamma_*[\R/T\Z]\in H_1(Y)$.

If $\gamma$ is a Reeb orbit passing through a point $y\in Y$, then the
linearized return map of the flow $R$ along $\gamma$ determines a
symplectic linear map $P_\gamma:\xi_y\to\xi_y$.  The eigenvalues of
$P_\gamma$ do not depend on $y$.  We say that the Reeb orbit $\gamma$
is {\em nondegenerate\/} if $P_\gamma$ does not have $1$ as an
eigenvalue.  We assume throughout that all Reeb orbits are
nondegenerate; this condition holds for generic contact forms
$\lambda$.

We now choose an almost complex structure $J$ on the 4-manifold
$\R\times Y$.  We always assume that $J$ is ``admissible'' in the following
sense.

\begin{definition}
An almost complex structure $J$ on $\R\times Y$ is {\em admissible\/}
if:
\begin{itemize}
\item $J(\partial_s)=R$, where $s$ denotes the $\R$ coordinate on
$\R\times Y$.
\item $J(\xi)=\xi$.
\item
$J$ rotates $\xi$ ``positively'' in the
sense that $d\lambda (v,Jv)\ge 0$ for all $v\in\xi$.
\item $J$ is invariant under the $\R$ action on $\R\times Y$ that
translates $s$.
\end{itemize}
\end{definition}

A {\em $J$-holomorphic curve\/} in $\R\times Y$ is a triple $(C,j,u)$
where $C$ is a smooth surface, $j$ is a complex structure on $C$, and
$u:C\to\R\times Y$ is a smooth map such that $J\circ du=du\circ j$.
The triple $(C,j,u)$ is equivalent to the triple $(C',j',u')$ iff
there is a biholomorphic map
$\phi:(C,j)\stackrel{\simeq}{\longrightarrow} (C',j')$ such that
$u'\circ \phi = u$.  We always assume that the domain $(C,j)$ is a
punctured compact Riemann surface, possibly disconnected.  We usually
denote a $J$-holomorphic curve simply by $u$.

A {\em positive end\/} of $u$ at a Reeb orbit $\gamma$ is an end of
$C$ that can be parametrized by $(s,t)\in[R,\infty)\times S^1$ for
some $R\in\R$, such that $u(s,t)=(s,y(s,t))$ and
$\lim_{s\to\infty}y(s,\cdot)$ is a reparametrization of $\gamma$.
A {\em negative end\/} of $u$ at $\gamma$ is defined analogously with
$s\in(-\infty,-R]$.

For each embedded Reeb orbit $\gamma$, fix a point $y\in Y$ in the
image of $\gamma$.  If $m$ is a positive integer, let $\gamma^m$
denote the Reeb orbit that is an $m$-fold cover of $\gamma$.  If $u$
has an end at $\gamma^m$, then the intersection of this end with
$\{s\}\times Y$ for $|s|$ large is an $m$-fold covering of $\gamma$
via a normal bundle projection.  An {\em asymptotic marking\/} of the
end is an inverse image of $y$ under this covering.  This notion does
not depend on the choice of normal bundle projection or on the choice
of $s$ with $|s|$ large.  Note that $\Z/m$ acts freely and
transitively on the set of asymptotic markings of an end of $u$ at
$\gamma^m$.

\begin{definition}
\label{def:MSFT}
Let $\alpha=(\alpha_1,\ldots,\alpha_k)$ and
$\beta=(\beta_1,\ldots,\beta_l)$ be ordered lists of Reeb orbits,
possibly repeated.  Define $\mc{M}^J(\alpha,\beta)$ to be the moduli
space of $J$-holomorphic curves $u:C\to\R\times Y$ as above, such that
$u$ has ordered and asymptotically marked positive ends at
$\alpha_1,\ldots,\alpha_k$, ordered and asymptotically marked negative
ends at $\beta_1,\ldots,\beta_l$, and no other ends.
\end{definition}

Note that $\mc{M}^J(\alpha,\beta)\neq\emptyset$ only if
$\sum_{i=1}^k[\alpha_i] = \sum_{j=1}^l[\beta_j] \in H_1(Y)$.  Also,
since the $\R$ action on $\R\times Y$ preserves $J$, it induces an
$\R$ action on $\mc{M}^J(\alpha,\beta)$.

\begin{definition}
\label{def:ind}
If $u\in \mc{M}^J(\alpha,\beta)$, define the (Fredholm) {\em index\/}
\begin{equation}
\label{eqn:ind}
\op{ind}(u) \eqdef -\chi(C) + 2 c_1(u^*\xi,\tau) + \sum_{i=1}^k
\op{CZ}_\tau(\alpha_i) - \sum_{j=1}^l \op{CZ}_\tau(\beta_j).
\end{equation}
\end{definition}
\noindent
The terms on the right hand side of equation \eqref{eqn:ind} are
defined as follows.  First, $C$ is the domain of $u$ as above.
Second, $\tau$ is a trivialization of $\xi$ over the Reeb orbits
$\alpha_i$ and $\beta_j$; it turns out that $\op{ind}(u)$ does not
depend on $\tau$, although individual terms in its definition
do. Next, $c_1(u^*\xi,\tau)$ denotes the relative first Chern class of
the complex line bundle $u^*\xi$ over $C$ with respect to the
trivializations $\tau$ at the ends.  This is defined by counting the
zeroes of a generic section which at the ends is nonvanishing and
constant with respect to the chosen trivializations.  Finally,
$\op{CZ}_\tau(\gamma)$ denotes the Conley-Zehnder index of $\gamma$
with respect to $\tau$.

In the present setting where $\dim(Y)=3$, this Conley-Zehnder index is
described explicitly as follows.  Let $\gamma$ be an embedded Reeb
orbit.  Let $\tau$ be a trivialization of $\xi$ over $\gamma$, and use
$\tau$ to trivialize $\xi$ over $\gamma^m$ for each positive integer
$m$.  Our assumption that all Reeb orbits are nondegenerate implies
that the linearized return map $P_{\gamma^m}=P_\gamma^m$ does not have
$1$ as an eigenvalue.  Let $\lambda,\lambda^{-1}$ denote the
eigenvalues of $P_\gamma$.  We say that $\gamma$ is {\em positive
hyperbolic\/} if $\lambda,\lambda^{-1}>0$, {\em negative hyperbolic\/}
if $\lambda,\lambda^{-1} <0$, and {\em elliptic\/} if
$\lambda,\lambda^{-1}$ are on the unit circle.  If $\gamma$ is
hyperbolic, then there is an integer $n$ such that the linearized Reeb
flow along $\gamma$ rotates the eigenspaces by angle $\pi n$ with
respect to $\tau$, and
\begin{equation}
\label{eqn:CZHyp}
\op{CZ}_\tau(\gamma^m) = mn.
\end{equation}
The integer $n$ is even when $\gamma$ is positive hyperbolic and odd
when $\gamma$ is negative hyperbolic.  If $\gamma$ is elliptic, then
there is an irrational number $\theta$, which we call the ``monodromy
angle'', such that
\begin{equation}
\label{eqn:CZEll}
\op{CZ}_\tau(\gamma^m) = 2\floor{m\theta} +1.
\end{equation}
Here $\tau$ is homotopic to a trivialization in which the linearized
Reeb flow along $\gamma$ rotates by angle $2\pi \theta$.

We say that $u\in\mc{M}^J(\alpha,\beta)$ is ``not multiply covered''
if $u$ does not multiply cover any component of its image.  We say
that $u$ is ``unobstructed'' if the linear deformation operator
associated to $u$ is surjective; then $\mc{M}^J(\alpha,\beta)$ is a
manifold near $u$.  The following proposition is the 3-dimensional
case of a result proved in \cite{dragnev}, using an index calculation
from \cite{schwarz}.

\begin{proposition}
\label{prop:dimension}
If $J$ is generic, and if $u\in\mc{M}^J(\alpha,\beta)$ is not multiply
covered, then $u$ is unobstructed, so that $\mc{M}^J(\alpha,\beta)$ is
a manifold near $u$.  Moreover, this manifold has dimension $\op{ind}(u)$.
\end{proposition}

\noindent
Assume henceforth that $J$ is generic in this sense.

Following \cite{bm}, one can ``coherently'' orient all the moduli
spaces of non-multiply covered $J$-holomorphic curves by making one
orientation choice for each Reeb orbit.  (We use slightly different
conventions from \cite{bm}, and we make a canonical choice for each
elliptic Reeb orbit; see \cite[\S9]{obg2} for details.)
Given
coherent orientations, if $M$ is a non-$\R$-invariant component of
such a moduli space, we orient $M$ using the $\R$ direction first.
That is, if $u\in M$, if $v_1\in T_uM$ denotes the derivative of the
$\R$ action on $M$, and if $(v_1,\ldots,v_n)$ is an oriented basis for
$T_uM$, then we declare that the projection of $(v_2,\dots,v_n)$ is an
oriented basis for $T_u(M/\R)$.  If $u$ has index $1$, then the above
convention defines a sign, which we denote by
$\epsilon(u)\in\{\pm1\}$.

\begin{remark}
\label{rem:hypor}
It follows from \cite{bm} that a system of coherent orientations
behaves as follows under the diffeomorphisms between moduli spaces
obtained by changing the orderings and asymptotic markings of the
ends.  If one switches the order of two ends, then this switches the
orientation if and only if both ends are at positive hyperbolic Reeb
orbits.  (Here an even cover of a negative hyperbolic orbit is
classified as positive hyperbolic.)  If $\gamma$ is an embedded Reeb
orbit, and if one acts on the asymptotic marking of an end at
$\gamma^m$ by a generator of $\Z/m$, then this switches the
orientation if and only if $m$ is even and $\gamma$ is negative
hyperbolic.
\end{remark}

\subsection{Branched covered cylinders}

To prepare for the statement of the gluing theorem, we now calculate
the index of branched covers of $\R$-invariant cylinders.

\begin{definition}
\label{def:indTheta}
If $a_1,\ldots,a_k$ and $b_1,\ldots,b_l$ are positive integers with
\[
\sum_{i=1}^ka_i = \sum_{j=1}^l b_j = M
\]
and if $\theta$ is an irrational
number, define
\[
\op{ind}_\theta(a_1,\ldots,a_k \mid b_1,\ldots,b_l) \eqdef
2\left(\sum_{i=1}^k \ceil{a_i\theta} - \sum_{j=1}^l \floor{b_j\theta}
- 1\right).
\]
Note that $\op{ind}_\theta\ge 0$, because $\sum_i\ceil{a_i\theta} \ge
\ceil{M\theta}$ and $\sum_j\floor{b_j\theta} \le \floor{M\theta}$.
\end{definition}

\begin{lemma}
\label{lem:BCCIndex}
Suppose $u\in\mc{M}^J(\alpha,\beta)$ is a branched cover of
$\R\times\gamma$, where $\gamma$ is an embedded Reeb orbit.  Then
$\op{ind}(u) \ge 0$, with equality only if:
\begin{description}
\item{(i)}
Each component of the domain $C$ of $u$ has genus $0$.
\item{(ii)}
If $\gamma$ is hyperbolic, then the covering
$u:C\to\R\times\gamma$ has no branch points.
\end{description}
\end{lemma}

\begin{proof}
Without loss of generality, $C$ is connected; let $g$ denote its
genus.  Write $\alpha=(\gamma^{a_1},\ldots,\gamma^{a_k})$ and
$\beta=(\gamma^{b_1},\ldots,\gamma^{b_l})$.  To calculate
$\op{ind}(u)$, choose a trivialization $\tau$ of $\gamma^*\xi$,
and use this to trivialize $\xi$ over all of the ends of $u$.  Since
$\tau$ extends to a trivialization of $\xi$ over $\R\times\gamma$, it
follows that $c_1(u^*\xi,\tau)=0$.  Thus
\begin{equation}
\label{eqn:BCCIndex}
\op{ind}(u) = -\chi(C) + \sum_{i=1}^k \op{CZ}_\tau(\gamma^{a_i}) -
\sum_{j=1}^l \op{CZ}_\tau(\gamma^{b_j}).
\end{equation}

If $\gamma$ is hyperbolic, then by equation \eqref{eqn:CZHyp}, the
Conley-Zehnder index terms in equation \eqref{eqn:BCCIndex} cancel, so
$\op{ind}(u)=-\chi(C)\ge 0$.  If equality holds, then $C$ is a
cylinder, and by Riemann-Hurwitz there are no branch points.

Now suppose that $\gamma$ is elliptic with monodromy angle $\theta$ with
respect to $\tau$.  Then by equations \eqref{eqn:CZEll} and
\eqref{eqn:BCCIndex},
\begin{equation}
\label{eqn:indug}
\begin{split}
\op{ind}(u) & = (2g-2+k+l) + \sum_{i=1}^k \left(2\ceil{a_i\theta}-1\right)  -
\sum_{j=1}^l \left(2\floor{b_j\theta} + 1\right)\\
& = 2g + \op{ind}_\theta(a_1,\ldots,a_k \mid b_1,\ldots,b_l).
\end{split}
\end{equation}
Since $\op{ind}_\theta\ge 0$, it follows that $\op{ind}(u)\ge 0$, with
equality only if $g=0$.
\end{proof}

Recall that a {\em partition\/} of a nonnegative integer $M$ is a
list of positive integers $(a_1,\ldots,a_k)$ modulo reordering,
possibly with repetitions, such that $\sum_{i=1}^ka_i=M$.  In
connection with the above index calculation, we now define a partial
order on the set of partitions of $M$.

\begin{definition}
\label{def:newpo}
Fix $\theta$ irrational.  We say that $(a_1,\ldots,a_k)\ge_\theta
(b_1,\ldots,b_l)$ if whenever $\gamma$ is an elliptic Reeb orbit with
monodromy angle $\theta$, there exists an index zero branched cover of
$\R\times\gamma$ in $\mc{M}^J((\gamma^{a_1},\ldots,\gamma^{a_k}),
(\gamma^{b_1},\ldots,\gamma^{b_l}))$.  It is an exercise (which we
will not need) to check that $\ge_\theta$ is a partial order.
\end{definition}

\subsection{Statement of the gluing problem}

The following definition specifies the kinds of pairs of curves that
we will be gluing.

\begin{definition}
\label{def:gluingPair}
A {\em gluing pair\/} is a pair of immersed $J$-holomorphic curves
$u_+\in\mc{M}^J(\alpha_+,\beta_+)$ and $u_-\in\mc{M}^J(\beta_-,
\alpha_-)$ such that:
\begin{description}
\item{(a)}
$\op{ind}(u_+)=\op{ind}(u_-)=1$.
\item{(b)}
$u_+$ and $u_-$ are not
multiply covered, except that they may contain unbranched covers of
$\R$-invariant cylinders.
\item{(c)}
For each embedded Reeb orbit $\gamma$, the total covering multiplicity
of Reeb orbits covering $\gamma$ in the list $\beta_+$ is the same as
the total for $\beta_-$.
(In contrast, for the usual form of gluing one would assume that
$\beta_+=\beta_-$.)
\item{(d)}
If $\gamma$ is an elliptic embedded Reeb orbit with monodromy angle
$\theta$, let $a_1',\ldots,a_{k'}'$ denote the covering multiplicities of
the $\R$-invariant cylinders over $\gamma$ in $u_+$, and let
$b_1',\ldots,b_{l'}'$ denote the corresponding multiplicities in $u_-$.
Then under the partial order $\ge_\theta$ in
Definition~\ref{def:newpo}, the partition $(a_1',\ldots,a_{k'}')$ is
minimal, and the partition $(b_1',\ldots,b_{l'}')$ is maximal.
\end{description}
\end{definition}
\noindent

Let $(u_+,u_-)$ be a gluing pair.  Our gluing theorem computes an
integer $\#G(u_+,u_-)$ which, roughly speaking, is a signed count of
ends of the index 2 part of the moduli space
$\mc{M}^J(\alpha_+,\alpha_-)/\R$ that break into $u_+$ and $u_-$ along
with some index zero branched covers of $\R$-invariant cylinders
between them.  The precise definition of $\#G(u_+,u_-)$ is a bit
technical and occupies the rest of this subsection.  There is some
subtlety here when $u_+$ or $u_-$ contain covers of $\R$-invariant
cylinders; in this case we will use condition (d) above in showing
that $\#G(u_+,u_-)$ is well-defined.

To prepare for the definition of the count $\#G(u_+,u_-)$, we first
define a set $\mc{G}_\delta(u_+,u_-)$ of index $2$ curves in
$\mc{M}^J(\alpha_+,\alpha_-)$ which, roughly speaking, are close to
breaking in the above manner.  For the following definition, choose an
arbitrary product metric on $\R\times Y$.

\begin{definition}
\label{def:Gdelta}
For $\delta>0$, define $\mc{C}_\delta(u_+,u_-)$ to be the set of
immersed (except possibly for finitely many singular points) surfaces
in $\R\times Y$ that can be decomposed as $C_-\cup C_0\cup C_+$, such
that the following hold:
\begin{itemize}
\item
There is a real number $R_-$, and a section $\psi_-$ of the normal
bundle to $u_-$ with $|\psi_-|<\delta$, such that $C_-$ is the
$s\mapsto s+R_-$ translate of the $s\le 1/\delta$ part of the
exponential map image of $\psi_-$.
\item
Likewise, there is a real number $R_+$, and a section $\psi_+$ of the
normal bundle to $u_+$ with $|\psi_+|<\delta$, such that $C_+$ is the
$s\mapsto s+R_+$ translate of the $s\ge -1/\delta$ part of the
exponential map image of $\psi_+$.
\item
$R_+-R_- > 2/\delta$.
\item
$C_0$ is contained in the union of the radius $\delta$ tubular
neighborhoods of the cylinders $\R\times\gamma$, where $\gamma$ ranges
over the embedded Reeb orbits covered by orbits in $\beta_\pm$.
\item
$\partial C_0 = \partial C_- \sqcup \partial C_+$, where the positive
boundary circles of $C_-$ agree with the negative boundary circles of
$C_0$, and the positive boundary circles of $C_0$ agree with the
negative boundary circles of $C_+$.
\end{itemize}
Let $\mc{G}_\delta(u_+,u_-)$ denote the set of index $2$ curves in
$\mc{M}^J(\alpha_+,\alpha_-) \cap \mc{C}_\delta(u_+,u_-)$.
\end{definition}

To see that this definition does what it is supposed to, we have:

\begin{lemma}
\label{lem:dds}
Given a gluing pair $(u_+,u_-)$, there exists $\delta_0>0$ with the
following property.  Let $\delta\in(0,\delta_0)$ and let
$\{[u_n]\}_{n=1,2,\ldots}$ be a sequence in
$\mc{G}_{\delta}(u_+,u_-)/\R$.  Then there is a subsequence which
converges in the sense of \cite{behwz} either to a curve in
$\mc{M}^J(\alpha_+,\alpha_-)/\R$, or to a broken curve in which the
top level is $u_+$, the bottom level is $u_-$, and all intermediate
levels are unions of index zero branched covers of $\R$-invariant
cylinders.
\end{lemma}

\begin{proof}
By the compactness theorem in \cite{behwz}, any sequence of index 2
curves in $\mc{M}^J(\alpha_+,\alpha_-)/\R$ has a subsequence which
converges to some broken curve.  Moreover, the indices of the levels
of the broken curve sum to $2$.

If the sequence is in $\mc{G}_\delta(u_+,u_-)/\R$ with $\delta>0$
sufficiently small, then by Lemma~\ref{lem:BCCIndex} and the
definition of $\mc{G}_\delta$, one of the following two scenarios occurs:
\begin{description}
\item{(i)} One level of the broken curve contains the index $1$
component of $u_+$, and some lower level contains the index $1$
component of $u_-$.
\item{(ii)} Some level contains two index $1$ components or one index
$2$ component.
\end{description}
Moreover, all other components of all levels are index zero branched
covers of $\R$-invariant cylinders.  By condition (d) in the
definition of gluing pair, any covers of $\R$-invariant cylinders in
the top and bottom levels of the broken curve must be unbranched.  It
follows that in case (i), the top level is $u_+$ and the bottom level
is $u_-$, while in case (ii), there are no other levels.
\end{proof}

\begin{definition}
\label{def:countG}
Fix coherent orientations and generic $J$ as in
Proposition~\ref{prop:dimension}, and let $(u_+,u_-)$ be a gluing
pair. If $\delta\in(0,\delta_0)$, then by Lemma~\ref{lem:dds} one can
choose an open set $U\subset\mc{M}^J(\alpha_+,\alpha_-)/\R$ such that:
\begin{itemize}
\item
$\mc{G}_{\delta'}(u_+,u_-)/\R \subset U \subset
\mc{G}_{\delta}(u_+,u_-)/\R$ for some $\delta'\in(0,\delta)$.
\item
The closure $\overline{U}$ has finitely many boundary points.
\end{itemize}
Define $\#G(u_+,u_-)\in\Z$ to be minus the signed count of boundary
points of $\overline{U}$.  By Lemma~\ref{lem:dds} , this does not
depend on the choice of $\delta$ or $U$.
\end{definition}

\subsection{Statement of the main theorem}

Let $(u_+,u_-)$ be a gluing pair.  The main result of this paper gives
a combinatorial formula for $\#G(u_+,u_-)$.  To state the formula,
note first that by Lemma~\ref{lem:BCCIndex}, if $\#G(u_+,u_-)\neq0$
then for each hyperbolic Reeb orbit $\gamma$, the multiplicities of
the negative ends of $u_+$ at covers of $\gamma$ agree, up to
reordering, with the multiplicities of the positive ends of $u_-$ at
covers of $\gamma$.  When this is the case, assume that the orderings
of the negative ends of $u_+$ and of the positive ends of $u_-$ are such
that for each positive hyperbolic orbit $\gamma$, the aforementioned
multiplicities appear in the same order for $u_+$ and for $u_-$.  With
this ordering convention, the statement of the main theorem is as
follows:

\begin{theorem}
\label{thm:main}
Fix coherent orientations.  If $J$ is generic and if $(u_+,u_-)$ is a
gluing pair, then
\begin{equation}
\label{eqn:CGG}
\#G(u_+,u_-) = \epsilon(u_+)\epsilon(u_-)\prod_\gamma c_\gamma(u_+,u_-).
\end{equation}
Here the product is over embedded Reeb orbits $\gamma$ such that $u_+$
has a negative end at a cover of $\gamma$.  The integer
$c_\gamma(u_+,u_-)$, defined below, depends only on $\gamma$ and on
the multiplicities of the $\R$-invariant and non-$\R$-invariant
negative ends of $u_+$ and positive ends of $u_-$ at covers of
$\gamma$.
\end{theorem}

To complete the statement of Theorem~\ref{thm:main}, we now define the
``gluing coefficients'' $c_\gamma(u_+,u_-)$ that appear in equation
\eqref{eqn:CGG}.  We will use the following notation.  Let
$a_1,\ldots,a_k$ denote the multiplicities of the non-$\R$-invariant
negative ends of $u_+$ at covers of $\gamma$ (in some arbitrary
order).  Let $a_1',\ldots,a_{k'}'$ denote the multiplicities of the
$\R$-invariant components of $u_+$ at covers of $\gamma$.  Likewise,
let $b_1,\ldots,b_l$ denote the multiplicities of the
non-$\R$-invariant positive ends of $u_-$ at covers of $\gamma$, and
let $b_1',\ldots,b_{l'}'$ denote the multiplicities of the
$\R$-invariant components of $u_-$ at covers of $\gamma$.  We will
define
\begin{equation}
\label{eqn:cgamma}
c_\gamma(u_+,u_-) \eqdef c_\gamma(a_1,\ldots,a_k;
a_1',\ldots,a_{k'}' \mid b_1,\ldots,b_l; b_1',\ldots,b_{l'}'),
\end{equation}
where the right hand side of \eqref{eqn:cgamma} is defined below.

\subsection{The gluing coefficients $c_\gamma$ for hyperbolic $\gamma$}

The gluing coefficient $c_\gamma$ is relatively straightforward when
$\gamma$ is hyperbolic.  In this case the gluing over $\gamma$ does
not involve any branch points, by Lemma~\ref{lem:BCCIndex}.  One just
needs to match up negative ends of $u_+$ at covers of $\gamma$ with
positive ends of $u_-$ at covers of $\gamma$ with the same
multiplicity.  Also, when gluing two ends at the $m$-fold cover
$\gamma^m$, there are $m$ possibilities for matching up the sheets.
The signs of these different matchings are related according to
Remark~\ref{rem:hypor}.  In many cases, the various possibilities all
cancel out because of the orientations; while in the remaining cases,
all possibilities have the same sign.

\begin{definition}
\label{def:HGC}
Suppose $\gamma$ is hyperbolic.  Then $c_\gamma=0$ unless:
\begin{description}
\item{(a)}
The list of multiplicities $(a_1,\ldots,a_k,a_1',\ldots,a_{k'}')$ is a
permutation of the list $(b_1,\ldots,b_l,b_1',\ldots,b_{l'}')$.
\item{(b)} If $\gamma$ is positive hyperbolic, then the numbers
$a_1,\ldots,a_{k'}'$ are distinct.
\item{(c)}
If $\gamma$ is negative hyperbolic, then the numbers
$a_1,\ldots,a_{k'}'$ are all odd.
\end{description}
If (a), (b), and (c) hold, then for each positive integer $m$, let
$r(m)$ denote the number of times that the number $m$ appears in the
list $a_1,\ldots,a_{k'}'$, and define
\[
c_\gamma \eqdef \prod_{m=1}^\infty m^{r(m)} \cdot r(m)!
\]
\end{definition}

\subsection{The gluing coefficients $c_\gamma$ for elliptic $\gamma$}
\label{sec:ctheta}

The interesting case of the gluing coefficient $c_\gamma$ is when
$\gamma$ is elliptic with monodromy angle $\theta$.  Here the only
relevant feature of $\gamma$ is the irrational number $\theta$, so we
denote $c_\gamma$ by $c_\theta$.  In this section we give a recursive
definition of $c_\theta$ which is easy to compute with.  An alternate
definition of $c_\theta$ as a sum over forests, which is useful for
proving certain symmetry properties of $c_\theta$, is given in
\S\ref{sec:trees}.

To simplify the notation, denote the
arguments of the function $c_\theta$ by
\begin{equation}
\label{eqn:S}
S \eqdef (a_1,\ldots,a_k;
a_1',\ldots,a_{k'}' \mid b_1,\ldots,b_l; b_1',\ldots,b_{l'}').
\end{equation}
When $k'=0$ or $l'=0$ we drop the corresponding semicolon from the
notation \eqref{eqn:S}.  It is always assumed that
\begin{equation}
\label{eqn:sumCondition}
\sum_{i=1}^ka_i + \sum_{i=1}^{k'} a_i' = \sum_{j=1}^l b_j +
\sum_{j=1}^{l'} b_j'.
\end{equation}

\begin{definition}
If $S$ as in \eqref{eqn:S} satisfies \eqref{eqn:sumCondition}, define
a positive integer
\begin{equation}
\label{eqn:kappa}
\kappa_\theta(S) \eqdef
\sum_{i=1}^k\ceil{a_i\theta} +
\sum_{i=1}^{k'}\ceil{a_i'\theta} - \sum_{j=1}^l\floor{b_j\theta} -
\sum_{j=1}^{l'}\floor{b_j'\theta}.
\end{equation}
\end{definition}

The significance of $\kappa_\theta(S)$ is that by the calculation
\eqref{eqn:indug}, any index zero branched cover of $\R\times\gamma$
with positive ends of multiplicities $a_1,\ldots,a_{k'}'$ and negative
ends of multiplicities $b_1,\ldots,b_{l'}'$ must consist of
$\kappa_\theta(S)$ genus zero components.

To define $c_\theta(S)$, we first reduce to the case where
$\kappa_\theta(S)=1$.  We need to consider the different ways that the
ends of a branched cover can be divided among $\kappa_\theta(S)$
different components.

\begin{definition}
\label{def:TD}
A {\em $\theta$-decomposition\/} of $S$ is
decomposition
\begin{equation}
\label{eqn:thetadecomposition}
\begin{split}
\{1,\ldots,k\} &= I_1\sqcup \cdots \sqcup
I_{\kappa_\theta(S)},\\
\{1,\ldots,k'\} &= I_1'\sqcup \cdots \sqcup
I_{\kappa_\theta(S)}',\\
\{1,\ldots,l\} &= J_1\sqcup \cdots \sqcup
J_{\kappa_\theta(S)},\\
\{1,\ldots,l'\} &= J_1'\sqcup \cdots \sqcup
J_{\kappa_\theta(S)}',
\end{split}
\end{equation}
such that for each $\nu=1,\ldots,\kappa_\theta(S)$, the sets $I_\nu$,
$I'_\nu$, $J_\nu$, and $J'_\nu$ are not all empty, and
\begin{equation}
\label{eqn:Snu}
S_\nu \eqdef ((a_i\mid i\in I_\nu) ; (a_i' \mid i\in I_\nu') \mid
(b_j \mid j\in J_\nu) ; (b_j' \mid j\in J_\nu'))
\end{equation}
satisfies the sum condition \eqref{eqn:sumCondition}.  Note that since
$\kappa_\theta$ is always positive, we must have
$\kappa_\theta(S_\nu)=1$ for each $\nu$. 
  Declare two $\theta$-decompositions to
be equivalent iff they differ by applying a permutation of the set
$\{1,\ldots,\kappa_\theta(S)\}$ to the indexing on the right hand side
of \eqref{eqn:thetadecomposition}.  We sometimes abuse notation and
denote a $\theta$-decomposition by $\{S_\nu\}$.
\end{definition}

\begin{lemma}
\label{lem:potd}
With the notation of \eqref{eqn:S}, a $\theta$-decomposition of $S$
exists if and only if $(a_1,\ldots,a_{k'}')
\ge_\theta (b_1,\ldots,b'_{l'})$.
\end{lemma}

\begin{proof}
This follows directly from the above discussion.
\end{proof}

\begin{definition}
\label{def:GC4}
For any $S$ as in \eqref{eqn:S} satisfying the sum condition
\eqref{eqn:sumCondition}, define
\[
c_\theta(S) \eqdef \sum_{\substack{\mbox{\scriptsize equivalence
classes of}\\\mbox{\scriptsize $\theta$-decompositions of $S$}}}
\prod_{\nu=1}^{\kappa_\theta(S)} c_\theta(S_\nu).
\]
\end{definition}

To complete this definition, the rest of this subsection defines
$c_\theta(S)$ when $\kappa_\theta(S)=1$.  We first need to define some
auxiliary functions.

\begin{notation}
If $a$ and $b$ are positive integers, define a positive integer
\[
\delta_\theta(a,b) \eqdef b\ceil{a\theta} - a \floor{b\theta}.
\]
\end{notation}

\begin{definition}
\label{def:ftheta}
Given {\em ordered\/} lists of positive integers $a_1,\ldots,a_k$ and
$b_1,\ldots,b_l$ with the same sum
$\sum_{i=1}^ka_i=\sum_{j=1}^lb_j$, define a positive integer
$f_\theta(a_1,\ldots,a_k\mid b_1,\ldots,b_l)$ recursively as follows.
To start the recursion, if $k=l=0$, then $f_\theta(\mid)\eqdef 1$.

For $k\ge 1$, the recursion involves a sum over subsets
\[
I = \{i_1<\cdots < i_q\} \subset \{1,\ldots,l\}
\]
such that
\begin{equation}
\label{eqn:ICondition}
\sum_{j=1}^{q-1} b_{i_j} < a_1 \le \sum_{j=1}^q b_{i_j}.
\end{equation}
We also require that equality holds in \eqref{eqn:ICondition} only
when $k=1$.  (This requirement is automatically satisfied in the case
of interest where $\kappa_\theta=1$.)  The formula is now
\begin{equation}
\label{eqn:FTS}
f_\theta(a_1,\ldots,a_k \mid b_1,\ldots,b_l) \eqdef \sum_I
f_\theta(a_2,\ldots,a_k\mid b_I) \prod_{n=1}^q \delta_\theta
\bigg(a_1- \sum_{j=1}^{n-1} b_{i_j} \, , \, b_{i_n}\bigg).
\end{equation}
Here $b_I$ denotes the arguments $b_i$ for $i\notin I$, arranged in
order, together with (when $k>1$) one additional argument equal to
$\sum_{i=2}^ka_i - \sum_{i\notin I}b_i$, inserted in the position that
$b_{i_q}$ would occupy in the order.
\end{definition}

\begin{remark}
If $\kappa_\theta=1$, then $f_\theta$ is always a positive integer,
because the sum \eqref{eqn:FTS} always has at least one term.  This
follows by induction, since one can find a subset $I$ satisfying
\eqref{eqn:ICondition} by just taking $I=\{1,\ldots,q\}$, where $q$ is the
smallest integer such that $\sum_{j=1}^q b_j \ge a_1$.
\end{remark}

The definition of $c_\theta(S)$ when $\kappa_\theta(S)=1$ is now divided into
several cases depending on the value of $k'+l'$.

\begin{definition}
\label{def:ctheta}
If $\kappa_\theta(S)=1$ and $k'=l'=0$, then
$c_\theta(S)$ is defined as follows.  Choose a reordering of the
$a_i$'s and $b_j$'s so that
\begin{equation}
\label{eqn:reordering}
\frac{\ceil{a_i\theta}}{a_i} \le \frac{\ceil{a_{i+1}\theta}}{a_{i+1}},
\quad \quad \frac{\floor{b_j\theta}}{b_j} \ge
\frac{\floor{b_{j+1}\theta}}{b_{j+1}}.
\end{equation}
Then
\begin{equation}
\label{eqn:CTS}
c_\theta(S)
\eqdef
f_\theta(a_1,\ldots,a_k \mid b_1,\ldots,b_l).
\end{equation}
\end{definition}

\begin{remark}
\label{rem:notobvious}
It is not obvious from Definition~\ref{def:ctheta} that $c_\theta(S)$
is independent of the choice of reordering satisfying
\eqref{eqn:reordering}.  This fact follows from the analysis used to
prove Theorem~\ref{thm:main}, in particular Corollary~\ref{cor:REC}
and Proposition~\ref{prop:count}.  It can also be proved
combinatorially, as described in Remark~\ref{rem:IHX}.

Another nonobvious property which follows from the analysis
is the symmetry
\begin{equation}
\label{eqn:cthetasym}
c_\theta(a_1,\ldots,a_k \mid b_1,\ldots,b_l) =
c_{-\theta}(b_1,\ldots,b_l \mid a_1,\ldots,a_k).
\end{equation}
We will give a combinatorial proof of this in \S\ref{sec:trees}.
\end{remark}

\begin{remark}
The geometric significance of the condition \eqref{eqn:reordering} is that
if a $J$-holomorphic curve $u$ has ordered negative ends of
multiplicities $a_1,\ldots,a_k$ at a given elliptic Reeb orbit of
monodromy angle $\theta$, then generically, the $i^{th}$ negative end
of $u$ decays no faster than the $(i+1)^{st}$
negative end.  Likewise, if $u$ has ordered positive ends of multiplicities
$b_1,\ldots,b_l$ at this Reeb orbit, then generically, the $j^{th}$
positive end decays no faster than the $(j+1)^{st}$ positive end.
\end{remark}

\begin{definition}
\label{def:change}
If $\kappa_\theta(S)=1$ and $k'+l'=1$, re-order the $a_i$'s and
$b_j$'s in accordance with \eqref{eqn:reordering}.  If $k'=1$, define
\[
c_\theta(S) \eqdef f_\theta(a_1,\ldots,a_k,a_1'\mid b_1,\ldots,b_l).
\]
Likewise, if $l'=1$ define
\[
c_\theta(S) \eqdef f_\theta(a_1,\ldots,a_k \mid b_1,\ldots,b_l,b_1').
\]
\end{definition}

\begin{definition}
\label{def:OTE}
If $\kappa_\theta(S)=1$ and $k'+l'\ge 2$, define $c_\theta(S)\eqdef 0$ unless
$k=l=0$ and $k'=l'=1$, in which case define
\[
c_\theta(;a'|;a')\eqdef a'.
\]
\end{definition}

\subsection{Examples and applications}

The following are some important examples of elliptic gluing coefficients
$c_\theta(S)$ where $k'=l'=0$.

\begin{example}
If $k=l=1$, then $c_\theta(a\mid a)=a$.  In the gluing theorem this
corresponds to a situation with no branch points, and reflects the fact
that there are $a$ ways to match up the sheets of a negative end of
$u_+$ and a positive end of $u_-$ along an $a$-fold cover of Reeb
orbit.
\end{example}

\begin{example}
Suppose $k=1$ and $\theta\in (0,1/a_1)$.  (Elliptic orbits with
$\theta$ close to zero arise when $\lambda$ is a perturbation of a
Morse-Bott contact form.)  Then $\kappa_\theta(a_1\mid
b_1,\ldots,b_l)=1$, the sum in \eqref{eqn:FTS} has only one term, and
we find that
\[
c_\theta(a_1 \mid b_1,\ldots,b_l) = \prod_{j=1}^l b_j.
\]
\end{example}

\begin{example}
In \S\ref{sec:ech} we will use Theorem~\ref{thm:main} to prove that
the differential $\partial$ in embedded contact homology (ECH)
satisfies $\partial^2=0$.  For this purpose, one has to calculate
$c_\theta$ in a certain special (but nontrivial) case as follows.
Given an irrational number $\theta$, for each nonnegative integer $M$
there are two distinguished partitions of $M$, called the ``incoming
partition'' and the ``outgoing partition'', and denoted here by
$\pin_\theta(M)$ and $\pout_\theta(M)$ respectively, see
\S\ref{sec:pinpout}.  If $u$ is a $J$-holomorphic curve contributing
to the ECH differential, then the multiplicities of the negative
(resp.\ positive) ends of $u$ at each elliptic Reeb orbit are
determined by the corresponding incoming (resp.\ outgoing) partitions.
Thus a key part of the proof that $\partial^2=0$ is to show that if
\begin{equation}
\label{eqn:delicate1}
\pin_\theta(M)=(a_1,\ldots,a_k), \quad\quad
\pout_\theta(M)=(b_1,\ldots,b_l),
\end{equation}
then
\begin{equation}
\label{eqn:delicate2}
c_\theta(a_1,\ldots,a_k \mid b_1,\ldots,b_l) = 1.
\end{equation}
We will do so in Proposition~\ref{prop:pinpout}.  In fact, similar
arguments show that the converse is also true, i.e.\ \eqref{eqn:delicate2}
implies \eqref{eqn:delicate1}.  Thus the proof here that
$\partial^2=0$ is quite delicate.
\end{example}

\begin{remark}
Symplectic field theory \cite{egh} defines a differential $D$ on a
supercommutative algebra over $\Q$ generated by all ``good'' (not necessarily
embedded) Reeb orbits.  The differential $D$ counts points in some
abstract perturbations of the compactified moduli spaces of index $1$
$J$-holomorphic curves in $\R\times Y$.  Even if one knows all
$J$-holomorphic curves, it is a nontrivial problem to read off the
differential $D$.  Theorem~\ref{thm:main} gives some constraints on
the answer in our three-dimensional case (SFT is defined for contact
manifolds of any dimension).  To give the simplest example, let
$\gamma_1$ be an embedded elliptic Reeb orbit with monodromy angle
$\theta\in(0,1/2)$, and let $\gamma_2$ denote the double cover of
$\gamma_1$.  Suppose that the only index $1$ $J$-holomorphic curves
modulo translation with ends at $\gamma_1$ or $\gamma_2$ are a curve
$u_+\in\mc{M}^J((\gamma_+),(\gamma_2))$ and a curve
$u_-\in\mc{M}^J((\gamma_1,\gamma_1)),(\gamma_-))$, say with
$\epsilon(u_\pm)=1$.  There is some contribution $D_+\in\Q$ to the
differential coefficient $\langle D\gamma_+,\gamma_1^2\rangle$ arising
from broken curves consisting of $u_+$ together with a branched double
cover of $\R\times\gamma_1$.  Likewise there is some contribution
$D_-\in\Q$ to the differential coefficient $\langle
D\gamma_2,\gamma_-\rangle$ arising from broken curves consisting of a
branched double cover of $\R\times\gamma_1$ together with $u_-$.
Presumably $D_+$ and $D_-$ may depend on the choice of abstract
perturbations.  However, since $D^2=0$, Theorem~\ref{thm:main} requires
that
\[
D_+ + D_- = c_\theta(2\mid 1,1)=1.
\]
\end{remark}

\subsection{Overview of the proof of the main theorem}
\label{sec:overview}

We now describe the proof of Theorem~\ref{thm:main}.  To simplify
notation, we will restrict attention to the special case where
conditions (i) and (ii) below hold:
\begin{description}
\item{(i)}
There is an embedded elliptic Reeb orbit $\alpha$ such that all
negative ends of $u_+$ and all positive ends of $u_-$ are at covers of
$\alpha$.
\end{description}
Let $\theta$ denote the monodromy angle of $\alpha$.  Suppose that
$u_+$ has negative ends at $\alpha^{a_1},\ldots,\alpha^{a_{N_+}}$, and
$u_-$ has positive ends at $\alpha^{a_{-1}},\ldots,\alpha^{a_{-N_-}}$.
The second condition is then:
\begin{description}
\item{(ii)}
$\kappa_\theta(a_1,\ldots,a_{N_+}\mid a_{-1},\ldots,a_{-N_-}) = 1$.
\end{description}

The strategy for gluing $u_+$ and $u_-$ is as follows.  Fix large
constants $R>>r>>0$.  Let $\Sigma$ be a connected genus zero branched
cover of $\R\times \alpha$ which has positive ends of covering
multiplicities $a_1,\ldots,a_{N_+}$ and negative ends of covering
multiplicities $a_{-1},\ldots,a_{-N_+}$, such that all ramification
points have $|s|\le R$.  Form a ``preglued'' curve by using appropriate
cutoff functions to patch the negative ends of the $s\mapsto s+R+r$
translate of $u_+$ to the positive ends of $\Sigma$, and the positive
ends of the $s\mapsto s-R-r$ translate of $u_-$ to the negative ends
of $\Sigma$.  Now try to perturb the preglued curve to a
$J$-holomorphic curve, where near the ramification points of the
branched cover we only perturb in the directions normal to $\R\times\alpha$.

It turns out that there is an ``obstruction bundle'' $\mc{O}$ over the
moduli space $\mc{M}_R$ of branched covers $\Sigma$ as above, and a
section $\mathfrak{s}:\mc{M}_R\to\mc{O}$ of this bundle, such that the
preglued curve determined by $\Sigma\in\mc{M}_R$ can be perturbed as
above to a $J$-holomorphic curve if and only if $\mathfrak{s}(\Sigma)=0$.
In this way we will identify the count $\#G(u_+,u_-)$ in
Theorem~\ref{thm:main} with $\epsilon(u_+)\epsilon(u_-)$ times an
appropriate count of the zeroes of $\mathfrak{s}$.  The section $\mathfrak{s}$
is defined rather indirectly from the analysis, but there is an
approximation $\mathfrak{s}_0$ to $\mathfrak{s}$ which is given by an explicit
formula.  We will see that if $J$ is generic, then the sections
$\mathfrak{s}_0$ and $\mathfrak{s}$ have the same count of zeroes, because one
can deform $\mathfrak{s}$ to $\mathfrak{s}_0$ without any zeroes crossing the
boundary of the moduli space $\mc{M}_R$.  We will then use a detailed
analysis of the obstruction bundle to count the zeroes of $\mathfrak{s}_0$
and recover the combinatorial gluing coefficient $c_\theta(u_+,u_-)$.

Without conditions (i) and (ii) above, one also needs to keep track of
the different Reeb orbits where gluing takes place, and also to
consider disconnected branched covers of $\R$ cross an elliptic Reeb
orbit.  Since this does not involve any additional analysis, in an
attempt to keep the notation manageable we will continue to assume (i)
and (ii) below and in \cite{obg2}.

The harder analytic parts of the above proof are carried out in the sequel
\cite{obg2}.  The present paper explains the more algebraic aspects and
is organized as follows.  In \S\ref{sec:obstructionBundle} we define
the obstruction bundle over the moduli space of branched covers and
discuss its basic properties.  In \S\ref{sec:section} we define the
section $\mathfrak{s}_0$ of the obstruction bundle and quote results from
\cite{obg2} relating $\#G(u_+,u_-)$ to an
appropriate count of zeroes of $\mathfrak{s}_0$.  \S\ref{sec:trees} is
almost completely independent of the previous two sections, and
discusses the combinatorics of the gluing coefficients $c_\theta$ in
detail.  \S\ref{sec:count} brings the analysis and the combinatorics
together to count the zeroes of $\mathfrak{s}_0$, thereby completing the
proof of Theorem~\ref{thm:main} modulo the aforementioned results from
\cite{obg2}.
\S\ref{sec:approx} then ties up loose ends by proving some estimates
on the obstruction bundle which were used in
\S\ref{sec:obstructionBundle} and
\S\ref{sec:count}.  Finally, \S\ref{sec:ech} explains the
application to embedded contact homology;  this section is independent
of \S\ref{sec:obstructionBundle}--\S\ref{sec:approx}.

For some other obstruction bundle calculations in symplectic field
theory, concerning index $1$ branched covers of $\R$-invariant
cylinders, see \cite{fabert}.

\section{The obstruction bundle}
\label{sec:obstructionBundle}

As described in \S\ref{sec:overview}, the number of gluings in
Theorem~\ref{thm:main} is determined by counting zeroes of a certain
section of an ``obstruction bundle'' $\mc{O}$ over a moduli space
$\mc{M}$ of genus zero branched covers of $\R\times S^1$.  We now
define this bundle and discuss its basic properties.  In
\S\ref{sec:BCT}  we define the moduli space
$\mc{M}$, and in \S\ref{sec:asop} we review the asympotic operator
associated to a Reeb orbit.  We then define the bundle $\mc{O}$ in
\S\ref{sec:DSigma}; the fiber of $\mc{O}$ over a branched cover
$\Sigma$ is the dual of the cokernel of a certain operator $D_\Sigma$.
In \S\ref{sec:SCE} we introduce some special elements of
$\Coker(D_\Sigma)$; later we will study the section of $\mc{O}$ by
evaluating it on these.  In \S\ref{sec:ENCE} we give some estimates on
the behavior of a special cokernel element in terms of the
combinatorics of the branched cover on which it is defined.  Finally,
\S\ref{sec:OBO} defines an orientation of $\mc{O}$, and
\S\ref{sec:cmr} defines a useful compactification of $\mc{M}/\R$.

In this section fix positive integers $a_1,\ldots,a_{N_+}$ and
$a_{-1},\ldots,a_{-N_-}$ with
\begin{equation}
\label{eqn:M}
\sum_{i=1}^{N_+} a_i = \sum_{j=-1}^{-N_-}a_{j} = M.
\end{equation}
Write $N\eqdef N_++N_-$, and to avoid trivialities assume that $N>2$.
Also, fix an irrational number $\theta$, and assume that
\begin{equation}
\label{eqn:kappa1}
\kappa_\theta(a_1,\ldots,a_{N_+} \mid a_{-1},\ldots,a_{-N_-}) = 1.
\end{equation}
Finally, fix an admissible almost complex structure $J$ on $\R\times Y$
and an embedded Reeb orbit $\alpha$.  In
\S\ref{sec:DSigma}--\S\ref{sec:OBO} we assume that $\alpha$ is
elliptic with monodromy angle $\theta$.

\subsection{Branched covers and trees}
\label{sec:BCT}

The following basic definitions will be used throughout the paper.

\begin{definition}
\label{def:M}
Let $\mc{M}=\mc{M}(a_1,\ldots,a_{N_+} \mid a_{-1},\ldots,a_{-N_-})$
denote the moduli space of degree $M$ branched covers $\pi:\Sigma\to
\R\times S^1$ such that:
\begin{itemize}
\item
$\Sigma$ is connected and has genus zero.
\item
The positive ends of $\Sigma$ are labeled by $1,\ldots,N_+$, and the
negative ends of $\Sigma$ are labeled by $-1,\ldots,-N_-$.
\item
The end of $\Sigma$ labeled by $i$ has covering multiplicity $a_i$.
\item
The ends are asymptotically marked.  That is, an
identification is chosen between the $i^{th}$ positive end of $\Sigma$ and
$[R,\infty)\times (\R/2\pi a_i\Z)$, respecting the projection to
$\R\times (\R/2\pi\Z)$.  Likewise for the negative ends.
\end{itemize}
We declare $\pi:\Sigma\to\R\times S^1$ to be equivalent to $\pi':
\Sigma' \to \R\times S^1$ if there is a diffeomorphism
$\phi:\Sigma\stackrel{\simeq}{\to} \Sigma'$ such that
$\pi'\circ\phi=\pi$, and $\phi$ respects the labelings and asymptotic
markings of the ends.  We often abuse notation and denote an element
of $\mc{M}$ by $\Sigma$.
\end{definition}

Note that $\mc{M}$ is a finite-sheeted covering space of the space of
meromorphic functions on ${\mathbb C}{\mathbb P}^1$ with poles of
order $a_1,\ldots,a_{N_+}$ and zeroes of order
$a_{-1},\ldots,a_{-N_-}$, modulo automorphisms of ${\mathbb C}{\mathbb
P}^1$.  In particular, $\mc{M}$ is a complex manifold of
dimension
\begin{equation}
\label{eqn:dimCM}
\dim_\C(\mc{M}) = N-2.
\end{equation}

We now explain how to associate, to each branched cover $\Sigma\in\mc{M}$,
a tree with certain additional structure. In this paper, a {\em
tree\/} is a finite, connected, simply connected graph $T$, such that
each vertex has degree either one (a {\em leaf\/}) or at least three
(an {\em internal vertex\/}).  We denote the set of internal vertices
by $\dot{V}(T)$.  The tree $T$ is {\em trivalent\/} if every internal
vertex has degree three.

For any two vertices $v$ and $w$ in a tree, let $P_{v,w}$ denote the
unique (nonbacktracking combinatorial) path from $v$ to $w$.  Given
three distinct leaves $i$, $j$, and $k$, the triple intersection of
the paths $P_{i,j}$, $P_{i,k}$, and $P_{j,k}$ consists of a single
vertex, which we call the {\em central vertex\/} for $i$, $j$, and
$k$.

\begin{definition}
\label{def:OWT}
An {\em oriented weighted tree\/} is a tree $T$ such that:
\begin{itemize}
\item
Each edge $e$ has an orientation $\mathfrak{o}(e)$ and a positive integer
weight $m(e)$, which we call the ``multiplicity'' of $e$.
\item
For each internal vertex, the sum of the multiplicities of the
outgoing edges equals the sum of the multiplicities of the incoming
edges.
\end{itemize}
\end{definition}

In an oriented tree, we call a leaf {\em positive\/} if the incident
edge points towards the leaf, and {\em negative\/} otherwise.  An
``upward'' path will mean a positively oriented path, and a
``downward'' path will mean a negatively oriented path.  A vertex $v$
is a {\em splitting vertex\/} if it has at least two outgoing edges,
and a {\em joining vertex\/} if it has at least two incoming edges.

\begin{definition}
\label{def:T(S)}
Let $T(a_1,\ldots,a_{N_+}\mid a_{-1},\ldots,a_{-N_-})$ denote the set
of oriented weighted trees such that:
\begin{itemize}
\item
The positive leaves are labeled by
$1,\ldots,N_+$, and the negative leaves are labeled by
$-1,\ldots,-N_-$.
\item
The (edge incident to the) $i^{th}$ leaf
has multiplicity $a_i$.
\end{itemize}
\end{definition}

\begin{definition}
Let $\mc{T}=\mc{T}(a_1,\ldots,a_{N_+}\mid a_{-1},\ldots,a_{-N_-})$
denote the set of oriented weighted trees $T\in
T(a_1,\ldots,a_{N_+}\mid a_{-1},\ldots,a_{-N_-})$ such that:
\begin{itemize}
\item
Each internal vertex $v$ is labeled by a real number
$\rho(v)$.
\item
If $v$ and $w$ are internal vertices, and if there is an oriented edge
from $v$ to $w$, then $\rho(v)<\rho(w)$.
\end{itemize}
\end{definition}

\begin{figure}
\begin{center}
\scalebox{0.55}{\includegraphics{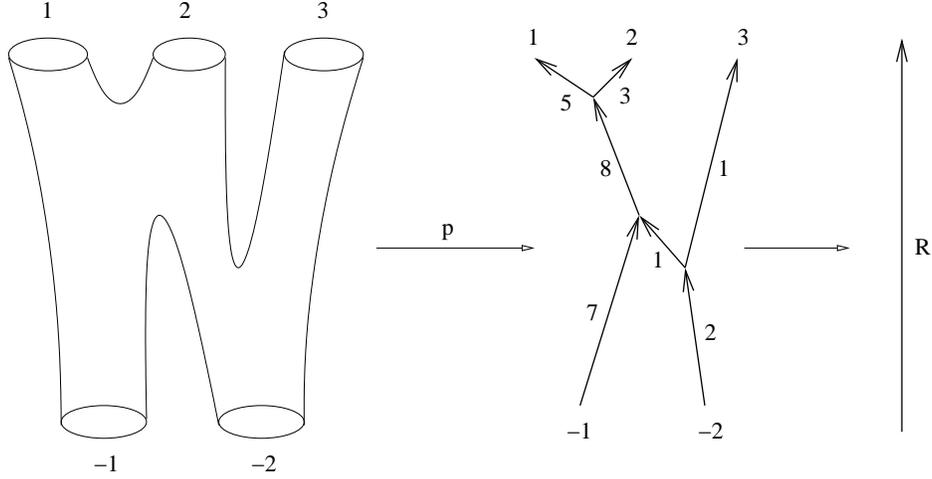}}
\end{center}
\caption{A branched cover $\Sigma\in\mc{M}(5,3,1\mid 7,2)$, its
associated oriented weighted tree $\tau(\Sigma)$, and their projections
to $\R$.  The end labels and edge weights are shown.}
\label{fig:tree}
\end{figure}

\begin{definition}
\label{def:RC}
Define a map
\[
\tau:\mc{M} \longrightarrow \mc{T}
\]
as follows.  Given a branched cover $\pi:\Sigma\to\R\times S^1$ in
$\mc{M}$, let $\rho$ denote the composition $\Sigma\to\R\times
S^1\to\R$.  Define two points in $\Sigma$ to be equivalent if they are
connected by a path on which $\rho$ is constant.  The quotient
space of $\Sigma$ by this equivalence relation is a one-dimensional
CW complex $\tau(\Sigma)$, which is a tree with a continuous map
\begin{equation}
\label{eqn:rhotree}
\rho:\tau(\Sigma)\longrightarrow\R.
\end{equation}
In the tree $\tau(\Sigma)$, a vertex $v$ of degree $d\ge 3$ corresponds to
an equivalence class $R(v)$ in $\Sigma$ containing ramification points
with total ramification index $d-2$.  The complement $\Sigma\setminus
\cup_v R(v)$ is a collection of cylinders, which correspond to the
edges of $\tau(\Sigma)$.  We orient the edges via the direction in
which $\rho$ increases, and define the multiplicity of an edge to be
the covering multiplicity of the corresponding cylinder in $\Sigma$.
An example is shown in Figure~\ref{fig:tree}.
\end{definition}

We next consider the extent to which the tree $T\eqdef \tau(\Sigma)$
determines the branched cover $\Sigma$.  In the ``generic'' case when
$T$ is trivalent, each internal vertex corresponds to a unique
ramification point in $\Sigma$, so there is a well-defined map
\begin{equation}
\label{eqn:phiT}
\phi_T: \tau^{-1}(T) \longrightarrow (S^1)^{\dot{V}(T)}
\end{equation}
which sends a branched cover $\Sigma$ to the $S^1$-coordinates of
$\pi$ of the ramification points.  Let $E(T)$ denote the set of edges
of $T$.

\begin{lemma}
\label{lem:covering}
If $T$ is trivalent, then the map $\phi_T$ in \eqref{eqn:phiT} is a
covering of degree
\begin{equation}
\label{eqn:degpiT}
\deg(\phi_T) = \prod_{e\in E(T)}m(e).
\end{equation}
\end{lemma}

\begin{proof}
Given a trivalent tree $T$ and an element of $S^1$ for each vertex, a
corresponding branched cover $\Sigma$ is obtained by taking a pair of
pants for each vertex and gluing them together as dictated by the
internal edges of $T$.  For each internal edge $e$ there are $m(e)$
possible gluings, and for each external edge $e$ there are $m(e)$
possible asymptotic markings of the corresponding end of $\Sigma$.
\end{proof}

The following subset of $\mc{M}$ will play a crucial role.

\begin{definition}
\label{def:M_R}
Given $R>0$, define
\[
\mc{M}_R=\mc{M}_R(a_1,\ldots,a_{N_+}\mid
a_{-1},\ldots,a_{-N_-})
\]
to be the set of $\pi:\Sigma\to\R\times S^1$ in $\mc{M}$ such that if
$x\in\Sigma$ is a ramification point and $\pi(x)=(s,t)$, then $|s|\le
R$.  Let $\partial\mc{M}_R$ denote the set of $\Sigma\in\mc{M}_R$
having a ramification point with $|s|=R$.
\end{definition}

\begin{lemma}
\label{lem:mrc}
$\mc{M}_R$ is compact\footnote{Note that we will use the assumption
\eqref{eqn:kappa1} here.  Thanks to A. Cotton-Clay for pointing out a
mistake in this regard in an earlier draft of this paper.}.
\end{lemma}

\begin{proof}
Let $X$ denote the symmetric product $\op{Sym}^{N-2}([-R,R]\times S^1)$.
Consider the map $\phi:\mc{M}_R\to X$ which sends a branched cover
$\pi:\Sigma\to\R\times S^1$ in $\mc{M}_R$ to the set of $\pi$-images
of the ramification points in $\Sigma$, repeated according to their
ramification indices.  Note that the symmetric product $X$ is compact,
the map $\phi$ is continuous, and each point in $X$ has only finitely
many inverse images under $\phi$.  Furthermore, $\phi$ defines a
covering space over each stratum in the symmetric product.  Thus to prove
that $\mc{M}_R$ is compact, it is enough to show that if
$\eta:[0,1]\to X$ is a path that maps all of $[0,1)$ to the same
stratum, then $\eta$ has a lift to $\mc{M}_R$ starting at any given
$\pi\in\phi^{-1}(\eta(0))$.  The issue is to check that whenever two
branch points in $[-R,R]\times S^1$ collide, the corresponding
ramification points in $\Sigma$ either do not interact or can be
merged.  More precisely, it is enough to show that if $x_0,
x_1\in\Sigma$ are two distinct ramification points, and if $\gamma$ is
an embedded path in $[-R,R]\times S^1$ from $\pi(x_0)$ to $\pi(x_1)$,
then $\gamma$ has at most one lift to a path in $\Sigma$.

Suppose to the contrary that $\gamma$ has two distinct lifts
$\widetilde{\gamma}_1$ and $\widetilde{\gamma}_2$.  We can then make a
new branched cover $\pi':\Sigma'\to\R\times S^1$ by cutting $\Sigma$
along the paths $\widetilde{\gamma}_1$ and $\widetilde{\gamma}_2$, and
gluing each side of $\widetilde{\gamma}_1$ to the opposite side of
$\widetilde{\gamma}_2$.  This operation reduces the total ramification
index by $2$.  Since the original branched cover $\Sigma$ had genus
zero, it follows by Riemann-Hurwitz that the new branched cover
$\Sigma'$ is disconnected.  On the other hand, $\Sigma'$ still has
positive ends of multiplicities $a_1,\ldots,a_{N_+}$, and negative
ends of multiplicities $a_{-1},\ldots,a_{-N_-}$.  Hence there are
decompositions  $\{1,\ldots,N_+\}=I_1\sqcup I_2$ and
$\{-1,\ldots,-N_-\}=J_1\sqcup J_2$ into proper subsets such that
$\sum_{i\in I_1}a_i =
\sum_{j\in J_1}a_j$. It follows that $\kappa_\theta(a_1,\ldots,a_{N_+}
\mid a_{-1},\ldots,a_{-N_-}) \ge 2$, contradicting the assumption
\eqref{eqn:kappa1}. 
\end{proof}

\subsection{The asymptotic operator}
\label{sec:asop}

We now review the asymptotic operator associated to a Reeb orbit.
This operator plays a fundamental role in the analysis.

Recall that we are fixing an embedded Reeb orbit $\alpha$.  By
rescaling the $t$ coordinate, we may assume that $\alpha$ is
parametrized by $S^1\eqdef\R/2\pi \Z$.  The linearized Reeb flow on
the contact planes along $\alpha$ defines a symplectic connection
$\nabla^R$ on the $2$-plane bundle $\alpha^*\xi$ over $S^1$.

\begin{definition}
Define the {\em asymptotic operator\/}
\[
L \eqdef L_\alpha \eqdef J\nabla^R_t
\]
acting on sections of $\alpha^*\xi$ over $S^1$.  More generally, if
$m$ is a positive integer, let $L_m\eqdef L_{\alpha^m}$ denote the pullback of
$L$ to $\R/2\pi m \Z$.
\end{definition}

To describe the operator $L_m$ more explicitly, fix
a complex linear, symplectic trivialization of $\alpha^*\xi$.  For
$t\in\R$, let $\Psi(t):\R^2\to\R^2$ denote the linearized Reeb flow
with respect to this trivialization, as the $S^1$ coordinate increases
from $0$ to $t$.  Let $J_0$ denote the standard complex structure on
$\R^2$.  For $t\in S^1$, define a matrix $S(t)$ by writing the
derivative of the linearized Reeb flow as
\begin{equation}
\label{eqn:psis}
\frac{d\Psi(t)}{dt}\Psi(t)^{-1} \fedqe J_0S(t).
\end{equation}
Then in the above trivialization,
\[
L_m = J_0\frac{d}{dt} + S(t)
\]
acting on complex functions on $\R/2\pi m\Z$.  Since the connection
$\nabla^R$ on $\alpha^*\xi$ is symplectic, it follows from
\eqref{eqn:psis} that the matrix $J_0S(t)$ is in the Lie algebra of
the symplectic group, which means that the matrix $S(t)$ is symmetric.
In particular, the operator $L_m$ is self-adjoint.

Our standing assumption that all Reeb orbits are nondegenerate implies
that $0\notin\op{Spec}(L_m)$.  The reason is that if $L_m\gamma=0$,
then it follows from equation \eqref{eqn:psis} that
$\gamma(t)=\Psi(t)\gamma(0)$.  Thus $\gamma(0)$ is an eigenvector of
$\Psi(2\pi m)$ with eigenvalue $1$, and if $\gamma\neq 0$ this
contradicts the nondegeneracy of $\alpha^m$.

To describe more spectral properties of the operator $L_m$, note that
if $\gamma$ is an eigenfunction with eigenvalue $\lambda$, then
$\gamma$ solves the ODE
\begin{equation}
\label{eqn:eigenvalue}
\frac{d\gamma(t)}{dt} = J_0(S(t)-\lambda)\gamma(t).
\end{equation}
It follows from \eqref{eqn:eigenvalue} and the uniqueness of solutions
to ODE's that if $\gamma$ is nonzero, then it is nonvanishing.  Then
the loop $\gamma:\R/2\pi m\Z\to\C$ has a well-defined winding number
around $0$.  We denote this winding number by $\eta(\gamma)\in\Z$.

\begin{example}
An important special case is where $S(t)=\theta$ for all $t$.  In this
case $\Phi(t)=e^{i \theta t}$, the operator $L_m$ is complex linear,
and the eigenfunctions are complex multiples of the functions
$\gamma(t)=e^{i\eta t/m}$ for $\eta\in\Z$.  Such an eigenfunction has
eigenvalue $\theta-\eta/m$ and winding number $\eta$.
\end{example}

It follows by analytic perturbation theory that in the general case,
eigenvalues are related to winding numbers as follows;  for details see
\cite[\S3]{hwz2}.

\begin{lemma}
\label{lem:hwz2}
\begin{description}
\item{(a)}
For each integer $\eta$, the sum of the eigenspaces whose nonzero
eigenfunctions have winding number $\eta$ is $2$-dimensional.
\item{(b)}
If $\gamma$ and $\gamma'$ are eigenfunctions of $L_m$ with eigenvalues
$\lambda\le\lambda'$, then $\eta(\gamma)\ge \eta(\gamma')$.
\item{(c)}
Suppose the Reeb orbit $\alpha$ is elliptic with monodromy angle
$\theta$.  If $\gamma$ is an eigenfunction of $L_m$
with eigenvalue $\lambda$, then
\[
\lambda > 0 \Longleftrightarrow \eta(\lambda) < m\theta.
\]
\end{description}
\end{lemma}

If $\eta$ and $m$ are integers with $m>0$, let $E_{\eta/m}^- \le
E_{\eta/m}^+$ denote the two eigenvalues of $L_m$, given by
Lemma~\ref{lem:hwz2}(a), whose associated eigenfunctions have winding
number $\eta$.

\begin{remark}
\label{rem:eigenvalues}
These eigenvalues depend only on the rational number
$\eta/m$, because if $\gamma$ is an eigenfunction of $L_m$ with
winding number $\eta$, and if $d$ is a positive integer, then the
covering $\R/2\pi d m \Z \to \R/2\pi m \Z$ pulls back $\gamma$ to an
eigenfunction of $L_{dm}$ with the same eigenvalue and with winding
number $d\eta$.
Also
\[
\frac{\eta_1}{m_1} < \frac{\eta_2}{m_2} \Longrightarrow E_{\eta_1/m_1}^- >
E_{\eta_2/m_2}^+;
\]
this follows by pulling back eigenfunctions to $\R/2\pi m_1m_2\Z$ and
applying Lemma~\ref{lem:hwz2}(b).
By Lemma~\ref{lem:hwz2}(c), if $\alpha$ is elliptic with monodromy
angle $\theta$, then the largest negative eigenvalue of $L_m$ is
$E_{\ceil{m\theta}/m}^+$, while the smallest positive eigenvalue of
$L_m$ is $E_{\floor{m\theta}/m}^-$.
\end{remark}

\subsection{The operator $D_\Sigma$ and the obstruction bundle}
\label{sec:DSigma}

We now introduce an operator $D_\Sigma$, which arises in connection
with deformations of $J$-holomorphic curves given by branched covers
of $\R\times\alpha$ in directions normal to $\R\times\alpha$.  

For the analysis to follow, fix a Hermitian metric on each
$\Sigma\in\mc{M}$ which varies smoothly over $\mc{M}$, which agrees
with the pullback of the standard metric on $\R\times S^1$ at points
in $\Sigma$ with distance $\ge 1$ from the ramification points, and
which within distance $1$ of a ramification point depends only on the
local structure of the branched cover within distance $2N$.  (Here the
``distance'' between two points $x,y\in\Sigma$ is defined to be the
infimum, over all paths $P$ in $\Sigma$ from $x$ to $y$, of the length
of the projection of $P$ to $\R\times S^1$.)

Now let $\pi:\Sigma\to\R\times S^1$ be a branched cover in $\mc{M}$.
Let $(s,t)$ denote the usual coordinates on $\R\times S^1$, and write
$z\eqdef s+it$.  Recall the notation $L_m$ and $S(t)$ from
\S\ref{sec:asop}.

\begin{definition}
Define a real linear operator 
\begin{gather}
\label{eqn:L2norms}
D_\Sigma: L^2_1(\Sigma,\C) \longrightarrow L^2(T^{0,1}\Sigma),\\
\nonumber
D_\Sigma \eqdef \dbar + \frac{1}{2}\pi^*(S(t)d\zbar).
\end{gather}
\end{definition}

Note that over an end of $\Sigma$ of multiplicity $m$, identified with
$[R,\infty)\times\R/2\pi m\Z$ or $(-\infty,-R]\times\R/2\pi m\Z$, we
have
\[
D_\Sigma f = (\partial_s + L_m)f\tensor\frac{d\zbar}{2}.
\]
Since $0\notin\op{Spec}(L_m)$, it follows by standard arguments that
the operator $D_\Sigma$ is Fredholm.

{\bf Assume henceforth that the Reeb orbit $\alpha$ is elliptic with
monodromy angle $\theta$.\/} Recall that $\theta$ is assumed to
satisfy \eqref{eqn:kappa1}.  The index of $D_\Sigma$ is then given as
follows:

\begin{lemma}
\label{lem:indDSigma}
$
\op{ind}(D_\Sigma) =
-
\dim_\R(\mc{M}).
$
\end{lemma}

\begin{proof}
By a standard index formula (cf.\ \cite{schwarz}),
\[
\op{ind}(D_\Sigma) = \chi(\Sigma) + \sum_{i=1}^{N_+} \mu(a_i) -
\sum_{j=-1}^{-N_-} \mu(a_j).
\]
Here $\mu(m)$, for a positive integer $m$, denotes the Maslov index of
the path of symplectic matrices $\{\Psi(t) \mid t\in[0,2\pi m]\}$,
which is given explicitly by
\[
\mu(m)=2\floor{m\theta}+1.
\]
Of course,
\[
\chi(\Sigma) = 2 - N.
\]
The lemma follows directly from the above three equations, together
with \eqref{eqn:kappa1} and \eqref{eqn:dimCM}.
\end{proof}

We now consider the cokernel of $D_\Sigma$.  We can identify
$\Coker(D_\Sigma)$ with the space of smooth $(0,1)$-forms $\sigma$ on
$\Sigma$ that are in $L^2$ and annihilated by the formal adjoint
$D_\Sigma^*$ of $D_\Sigma$.

A nonzero cokernel element $\sigma$ has the following asymptotic
behavior.  Over the $i^{th}$ positive end of $\Sigma$, which the
asymptotic marking identifies with $[R,\infty)\times\R/2\pi a_i\Z$,
write $\sigma=\sigma_i(s,t)d\zbar$.  The function $\sigma_i$ satisfies
the equation
\[
(\partial_s - L_{a_i})\sigma_i = 0.
\]
Since $\sigma$ is in $L^2$, we can expand
\begin{equation}
\label{eqn:PE}
\sigma_i(s,t) = \sum_{\lambda < 0} e^{\lambda
  s}\gamma_{i,\lambda}(t)
\end{equation}
where the sum is over negative eigenvalues $\lambda$ of $L_{a_i}$, and
$\gamma_{i,\lambda}$ is a (possibly zero) eigenfunction with
eigenvalue $\lambda$.  Let $\lambda_i$ denote the largest negative
eigenvalue for which $\gamma_{i,\lambda_i}$ is nonzero, and write
$\gamma_i\eqdef \gamma_{i,\lambda_i}$.  If $\kappa>0$ is the
difference between $\lambda_i$ and the second largest negative
eigenvalue, then by \eqref{eqn:PE}, there is an $s$-independent
constant $A$ such that
\[
\left|\sigma_i(s,t) - e^{\lambda_i s}\gamma_i(t)\right| \le A
e^{(\lambda_i-\kappa)s}.
\]
It follows that when $s$ is large, $\sigma_i$ has no zeroes, and has
winding number $\eta(\gamma_i)$ around the $i^{th}$ positive end of
$\Sigma$.  We denote this winding number by $\eta_i^+(\sigma)$.
Likewise, for $j\in\{-1,\ldots,-N_-\}$, one can expand $\sigma$ on the
$j^{th}$ negative end of $\Sigma$ as
\begin{equation}
\label{eqn:NE}
\sigma_j(s,t) = \sum_{\lambda > 0} e^{\lambda
  s}\gamma_{j,\lambda}(t),
\end{equation}
and $\sigma$ has a well-defined winding number $\eta_j^-(\sigma)$
around the $j^{th}$ negative end.  Lemma~\ref{lem:hwz2}(c) then gives
the winding bounds
\begin{equation}
\label{eqn:windbound}
\eta_i^+(\sigma) \ge \ceil{a_i\theta}, \quad\quad \eta_j^{-}(\sigma) \le
\floor{a_j\theta}.
\end{equation}

Given a nonzero cokernel element $\sigma$, let $Z(\sigma)$ denote the
number of ends of $\Sigma$ for which the inequalities
\eqref{eqn:windbound} are strict.

\begin{lemma}
\label{lem:windingBounds}
\begin{description}
\item{(a)} If $0\neq \sigma\in\Coker(D_\Sigma)$, then the zeroes of
$\sigma$ are isolated and have negative multiplicity, and the
algebraic count of zeroes is bounded by
\[
\#\sigma^{-1}(0) \ge \frac{\op{ind}(D_\Sigma)}{2} + 1 + Z(\sigma).
\]
\item{(b)} 
$\dim(\Coker(D_\Sigma))=-\op{ind}(D_\Sigma)$, or equivalently
$\Ker(D_\Sigma)=\{0\}$.
\end{description}
\end{lemma}

\begin{proof}
(a) Since $D_\Sigma^*$ is $\dbar^*$ plus a zeroth order term, the zeroes
of $\sigma$ are isolated and have negative multiplicity.  For any
$(0,1)$-form $\sigma$ with finitely many zeroes, the algebraic count
of zeroes is given by
\begin{equation}
\label{eqn:ACZ}
\#\sigma^{-1}(0) = \chi(\Sigma) + \sum_{i=1}^{N_+} \eta_i^+(\sigma) -
\sum_{j=-1}^{-N_-}
\eta_j^-(\sigma).
\end{equation}
If $\sigma$ is a nonzero cokernel element, then putting the winding
bounds \eqref{eqn:windbound} into \eqref{eqn:ACZ} and using
Lemma~\ref{lem:indDSigma} proves part (a).

(b) If $\dim(\Coker(D_\Sigma))>-\op{ind}(D_\Sigma)$, then one can find
a nonzero cokernel element $\sigma$ with zeroes at
$-\op{ind}(D_\Sigma)/2$ given points in $\Sigma$.  Since all zeroes of
$\sigma$ have negative multiplicity, this contradicts part (a).
\end{proof}

Lemma~\ref{lem:windingBounds}(b) implies that the cokernels of the
operators $D_\Sigma$ for $\Sigma\in\mc{M}$ comprise a smooth real
vector bundle over $\mc{M}$, which we denote by $\mc{O}^*$.

\begin{definition}
Define the {\em obstruction bundle\/} $\mc{O}\to\mc{M}$ to be the dual
of the bundle of cokernels $\mc{O}^*\to\mc{M}$.  Thus the fiber of
$\mc{O}$ over $\Sigma$ is
\[
\mc{O}_\Sigma = \Hom(\Coker(D_\Sigma),\R).
\]
\end{definition}

By Lemmas~\ref{lem:indDSigma} and \ref{lem:windingBounds}(b), the rank
of $\mc{O}$ equals the dimension of $\mc{M}$.

\subsection{Special cokernel elements}
\label{sec:SCE}

We now introduce some special elements of $\Coker(D_\Sigma)$ which
will play a key role in the obstruction bundle calculations in
\S\ref{sec:count}.

To define the special cokernel elements, for
$\sigma\in\Coker(D_\Sigma)$ we need to consider the ``leading terms''
of the asymptotic expansions \eqref{eqn:PE} and \eqref{eqn:NE}.
Namely, to the $i^{th}$ end of $\Sigma$ we associate a real vector
space $\mc{A}_i$ as follows.  For $i\in\{1,\ldots,N_+\}$, define
$\mc{A}_i$ to be the direct sum of those eigenspaces of $L_{a_i}$
whose nonzero elements have winding number $\ceil{a_i\theta}$.  By
Lemma~\ref{lem:hwz2}, $\mc{A}_i$ is two-dimensional, and consists of
the eigenspace for the largest negative eigenvalue of $L_{a_i}$,
together with the eigenspace for the second largest negative
eigenvalue when the largest negative eigenvalue has multiplicity one.
Likewise, for $j\in\{-1,\ldots,-N_-\}$, define $\mc{A}_j$ to be the
direct sum of those eigenspaces of $L_{a_j}$ whose nonzero elements
have winding number $\floor{a_j\theta}$.

Given $\sigma\in\Coker(D_\Sigma)$, for $i\in\{1,\ldots,N_+\}$, define
$\Phi_i(\sigma)\in\mc{A}_i$ to be the sum of those eigenfunctions
$\gamma_{i,\lambda}$ in the expansion \eqref{eqn:PE} that have winding
number $\ceil{a_i\theta}$.  Likewise, for $j\in\{-1,\ldots,-N_-\}$,
define $\Phi_j(\sigma)\in\mc{A}_j$ to be the sum of those
eigenfunctions in the expansion \eqref{eqn:NE} that have winding
number $\floor{a_j\theta}$.

We can now define the special cokernel elements:

\begin{definition}
If $i$, $j$, and $k$ label distinct ends of $\Sigma$, define a subspace
$V_{i,j,k}$ of $\Coker(D_\Sigma)$ by
\[
V_{i,j,k} \eqdef \left\{\sigma\in\Coker(D_\Sigma) \mid \Phi_l(\sigma) = 0
\;\;\;\forall l \notin \{i,j,k\} \right\}.
\]
\end{definition}

\begin{lemma}
\label{lem:SCE}
\begin{description}
\item{(a)}
$\Phi_i$ restricts to an isomorphism
$V_{i,j,k}\stackrel{\simeq}{\longrightarrow}\mc{A}_i$.
\item{(b)}
Every nonzero element of $V_{i,j,k}$ is nonvanishing.
\end{description}
\end{lemma}

\begin{proof}
By Lemmas~\ref{lem:indDSigma} and \ref{lem:windingBounds}(b), we know that
\begin{equation}
\label{eqn:dimCoker}
\dim(\Coker(D_\Sigma))=2(N -2).
\end{equation}
It follows that $\dim(V_{i,j,k})\ge 2$.  On the other hand,
Lemma~\ref{lem:windingBounds}(a) implies that if $0\neq\sigma\in
V_{i,j,k}$, then $\sigma$ is nonvanishing and $\Phi_i(\sigma)\neq 0$.
Assertions (a) and (b) follow.
\end{proof}

Consider now the oriented weighted tree $\tau(\Sigma)$ associated to the
branched cover $\pi:\Sigma\to\R\times S^1$.  A nonzero special
cokernel element $\sigma\in V_{i,j,k}$, being nonvanishing, has a
well-defined winding number around the cylinder in $\Sigma$ corresponding
to each edge $e$ of the tree $\tau(\Sigma)$, which we denote by
$\eta(\sigma,e)$.  We now derive a useful formula for these winding
numbers, which will play an essential role in the calculations in
\S\ref{sec:count}.  First a preliminary lemma:

\begin{lemma}
\label{lem:LFI}
Let $T\in T(a_1,\ldots,a_{N_+} \mid a_{-1},\ldots,a_{-N_-})$;  and let $v$
be an internal vertex of $T$ with outgoing edges of multiplicities
$m_1^+,\ldots,m_p^+$ and incoming edges of multiplicities
$m_1^-,\ldots,m_q^-$.  Then
\begin{equation}
\label{eqn:LFI}
\sum_{l=1}^p\ceil{m_l^+\theta} - \sum_{l=1}^q\floor{m_l^-\theta} - 1 = 0.
\end{equation}
\end{lemma}

\begin{proof}
Let $\op{ind}_\theta(v)$ denote twice the left hand side of
\eqref{eqn:LFI}.  By the sum condition on the weights,
$\sum_{l=1}^pm_l^+$ and $\sum_{l=1}^qm_l^-$ are equal, say to $r$, so
\[
\op{ind}_\theta(v) \ge 2(\ceil{r\theta} - \floor{r\theta} - 1) = 0.
\]
On the other hand, a straightforward calculation shows that the
quantity $\op{ind}_\theta$ from Definition~\ref{def:indTheta} satisfies
\[
\op{ind}_\theta(a_1,\ldots,a_{N_+} \mid a_{-1},\ldots,a_{-N_-}) =
\sum_{v\in\dot{V}(T)} \op{ind}_\theta(v).
\]
Thanks to our assumption \eqref{eqn:kappa1}, the left hand side of the
above equation is zero, and this completes the proof.
\end{proof}

\begin{lemma}
\label{lem:specialWinding}
Let $\sigma\in V_{i,j,k}$ be a nonzero special cokernel element, and
let $v$ denote the central vertex for $i$, $j$, and $k$.
\begin{itemize}
\item
If $e$ is an edge on one of the paths $P_{v,i}$, $P_{v,j}$, or $P_{v,k}$, then
\begin{equation}
\label{eqn:winding1}
\eta(\sigma,e) = \left\{\begin{array}{cl} \lceil m(e)\theta\rceil,
& \mbox{$e$ points away from $v$,} \\
\lfloor m(e)\theta\rfloor, & \mbox{$e$ points towards $v$.}
\end{array}
\right.
\end{equation}
\item
If $e$ is not on one of the paths
$P_{v,i}$, $P_{v,j}$, or $P_{v,k}$, then
\begin{equation}
\label{eqn:winding2}
\eta(\sigma,e) = \left\{\begin{array}{cl} \lceil m(e)\theta\rceil +
1,
& \mbox{$e$ points away from $v$,} \\
\lfloor m(e)\theta\rfloor - 1, & \mbox{$e$ points towards $v$.}
\end{array}
\right.
\end{equation}
\end{itemize}
\end{lemma}

\begin{proof}
Suppose first that $e$ is an external edge, incident to the $l^{th}$
leaf.  If $l\in\{i,j,k\}$, then by Lemma~\ref{lem:SCE}(a),
we have $\eta(\sigma,e)=\ceil{a_l\theta}$ when $l>0$ and
$\eta(\sigma,e)=\floor{a_l\theta}$ when $l>0$.  Thus equation
\eqref{eqn:winding1} holds in this case.  If $l\notin\{i,j,k\}$, then
the definition of $V_{i,j,k}$ implies that $\eta(\sigma,e)\ge
\ceil{a_l\theta}+1$ when $l>0$, and $\eta(\sigma,e)\le
\floor{a_l\theta}-1$ when $l<0$.  These inequalities must be
equalities, or else equation \eqref{eqn:ACZ} would give
$\#\sigma^{-1}(0)>0$, a contradiction.  Thus equation
\eqref{eqn:winding2} also holds in this case.

To prove that \eqref{eqn:winding1} and \eqref{eqn:winding2} hold for
internal edges $e$, we will use downward induction on the distance
(i.e.\ number of edges on the path) from $e$ to the central vertex
$v$.

To carry out the inductive step, let $w\neq v$ be an internal vertex
with outgoing edges $e_1^+,\ldots,e_p^+$ and incoming edges
$e_1^-,\ldots,e_q^-$.  Since $w\neq v$, there is a unique edge
incident to $w$ which is closest to $v$.  By symmetry, we may assume
that this edge is incoming, say $e_q^-$.  We may inductively assume
that the winding numbers of $\sigma$ around all other edges incident
to $w$ are given by
\eqref{eqn:winding1} and \eqref{eqn:winding2}. 

Counting zeroes of $\sigma$ as in
\eqref{eqn:ACZ} over a neighborhood in $\Sigma$ of the circle $R(w)$
(see Definition~\ref{def:RC}), and using the assumption that $\sigma$
is nonvanishing, shows that
\begin{equation}
\label{eqn:formZeroes}
\sum_{l=1}^p \eta(\sigma,e_l^+) - \sum_{l=1}^q
\eta(\sigma,e_l^-) = p+q-2.
\end{equation}
Now $e_q^-$ points
away from $v$, all other incoming edges of $w$ point towards $v$, and
all outgoing edges of $w$ point away from $v$.  If $e_q^-$ is on one
of the paths $P_{v,i}$, $P_{v,j}$, or $P_{v,k}$, then so is exactly
one other edge incident to the vertex $w$, whence by inductive
hypothesis,
\begin{equation}
\label{eqn:windingInduction}
\sum_{l=1}^p \eta(\sigma,e_l^+) - \sum_{l=1}^{q-1}
\eta(\sigma,e_l^-) = \sum_{l=1}^p \lceil
m(e_l^+)\theta\rceil - \sum_{l=1}^{q-1} \lfloor
m(e_l^-)\theta\rfloor + p + q - 2.
\end{equation}
Combining this with equation \eqref{eqn:formZeroes} and
Lemma~\ref{lem:LFI} gives
\[
\eta(\sigma,e_q^-) = \lceil m(e_q^-)\theta\rceil,
\]
as desired.  If $e_q^-$ is not on one of the paths $P_{v,i}$,
$P_{v,j}$, or $P_{v,k}$, then neither is any other edge adjacent to
$w$, so a modification of equation \eqref{eqn:windingInduction} holds
where we add $1$ to the right hand side, giving
\[
\eta(\sigma,e_q^-) = \lceil m(e_q^-)\theta\rceil + 1.
\qedhere
\]
\end{proof}

\subsection{Estimates on nonvanishing cokernel elements}
\label{sec:ENCE}

Given a nonvanishing cokernel element $\sigma\in\Coker(D_\Sigma)$,
we now state an estimate on the relative sizes of the restrictions of
$\sigma$ to different parts of $\Sigma$, and some other related
estimates, which will be used in \S\ref{sec:count} and \cite{obg2}.

To state the first estimate, recall that $\tau(\Sigma)$ denotes the oriented
weighted tree associated to $\Sigma$, which is a one-dimensional CW
complex with continuous maps $p:\Sigma\to\tau(\Sigma)$ and
$\rho:\tau(\Sigma)\to\R$, such that $\rho\circ p$ equals the composition
$\Sigma\stackrel{\pi}{\longrightarrow}\R\times S^1\to\R$.  We give
$\tau(\Sigma)$ the metric for which $\rho$ restricts to an isometry on
each edge.

Given $x,y\in \tau(\Sigma)$, let $P_{x,y}$ denote the path in
$\tau(\Sigma)$ from $x$ to $y$.  By an ``edge of $P_{x,y}$'', we mean
an edge $e$ of $\tau(\Sigma)$ such that a positive length subset of
$e$ is on the path $P_{x,y}$.  Let $P^+_{x,y}$ denote the set of edges
of $P_{x,y}$ that are oriented in the direction pointing from $x$ to
$y$, and let $P^-_{x,y}$ denote the set of edges of $P_{x,y}$ whose
orientation points from $y$ to $x$. If $e$ is an edge of $P_{x,y}$, let
$\ell(e)>0$ denote the length of the portion of $e$ that is on the
path $P_{x,y}$.

If $e$ is an edge of $\tau(\Sigma)$, then in the notation of
\S\ref{sec:asop}, define
\[
E_\pm(\sigma,e) \eqdef E^\pm_{\eta(\sigma,e)/m(e)}.
\]

\begin{proposition}
\label{prop:RS}
There exists $r'>0$ with the following property.  Let
 $\Sigma\in\mc{M}$ and suppose $\sigma\in\Coker(D_\Sigma)$ is
 nonvanishing.  Let $z,w\in \Sigma$ and define $x\eqdef p(z)$ and
 $y\eqdef p(w)$. Then
\begin{equation}
\label{eqn:RS}
\log|\sigma(w)| - \log|\sigma(z)| \le \sum_{e\in
P_{x,y}^+} \ell(e) E_+(\sigma,e) - \sum_{e\in P_{x,y}^-} \ell(e)
E_-(\sigma,e) + r'.
\end{equation}
\end{proposition}

The proof of Proposition~\ref{prop:RS} is given in
\S\ref{sec:RSProof}.  The following is an important special case:

\begin{corollary}
\label{cor:rsc}
If $S(t)=\theta$, then in Proposition~\ref{prop:RS} we can replace
\eqref{eqn:RS} by
\[
\left|\log|\sigma(w)| - \log|\sigma(z)|
-
\bigg(\sum_{e\in P_{x,y}^+} -
\sum_{e\in P_{x,y}^-} \bigg)
 \ell(e) \left(\theta-\frac{\eta(\sigma,e)}{m(e)}\right)
\right|
\le r'.
\]
\end{corollary}

\begin{proof}
In general, by switching the role of $z$ and $w$ in \eqref{eqn:RS} we
 get
\begin{equation}
\label{eqn:RS'}
\log|\sigma(w)| - \log|\sigma(z)| \ge \sum_{e\in
P_{x,y}^+} \ell(e) E_-(\sigma,e) - \sum_{e\in
P_{x,y}^-} \ell(e) E_+(\sigma,e) - r'.
\end{equation}
The assumption $S(t)=\theta$ implies that
\[
E_{\eta/m}^+ = E_{\eta/m}^- = \theta-\frac{\eta}{m}.
\]
Combining this with \eqref{eqn:RS} and \eqref{eqn:RS'} proves the
corollary.
\end{proof}

In the general case, another useful corollary of
Proposition~\ref{prop:RS} is that a special cokernel element decays
away from the central vertex, in the following sense:

\begin{corollary}
\label{cor:dacv}
There exist $c,\kappa>0$ with the following property.  Let
$\Sigma\in\mc{M}$, let $\sigma\in V_{i,j,k}$ be a special cokernel
element, let $z,w\in\Sigma$, and suppose that $p(z)$ is the central
vertex for $i,j,k$.  Then
\begin{equation}
\label{eqn:dacv}
|\sigma(w)| \le c \exp(-\kappa\cdot\op{dist}(p(w),p(z)))|\sigma(z)|.
\end{equation}
\end{corollary}

\begin{proof}
Let $x\eqdef p(z)$ and $y\eqdef p(w)$.  The winding number
calculations in Lemma~\ref{lem:specialWinding}, and the relation
between winding numbers and signs of eigenvalues in
Lemma~\ref{lem:hwz2}(c), show that $E_+(\sigma,e)<0$ for each $e\in
P_{x,y}^+$, and $E_-(\sigma,e)>0$ for each $e\in P_{x,y}^-$.  Now for
our given $\alpha$ and $a_i$'s,
there are only finitely many possible values of $E_\pm(\sigma,e)$.
We can then find $\kappa>0$ such that $E_+(\sigma,e)<-\kappa$ for each
$e\in P_{x,y}^+$ and $E_-(\sigma,e)>\kappa$ for each $e\in
P_{x,y}^-$.  Putting these inequalities into \eqref{eqn:RS}
proves \eqref{eqn:dacv}.
\end{proof}

We now state one more estimate which we will need.

\begin{definition}
\label{def:PiW}
Let $\Sigma\in\mc{M}$ and let $\sigma\in\Coker(D_\Sigma)$ be
nonvanishing.  Define a $(0,1)$-form $\Pi_W\sigma$ on those points
$z\in\Sigma$ for which $p(z)$ is not a vertex of $\tau(\Sigma)$, as
follows.  Let $x$ be a point in the CW-complex $\tau(\Sigma)$ which is
in the interior of an edge $e$.  Then $p^{-1}(x)$ is a circle in
$\Sigma$; choose an identification $p^{-1}(x)\simeq\R/2\pi m(e)\Z$
commuting with the projections to $\R\times S^1$.  Let $W(e)$ denote
the sum of the eigenspaces of $L_{m(e)}$ whose eigenfunctions have winding
number $\eta(\sigma,e)$.  On $p^{-1}(x)$, use $d\overline{z}$ to
identify $(0,1)$-forms with complex-valued functions, and define
$\Pi_W\sigma$ to be the $L^2$-orthogonal projection of $\sigma$ onto
$W(e)$.
\end{definition}

The following proposition, which is proved in \S\ref{sec:PAR},
asserts roughly that away from the ramification points, $\sigma$ is
well approximated by $\Pi_W\sigma$.

\begin{proposition}
\label{prop:approx1}
Given $\varepsilon_0>0$, there exists $R>1$ with the following property.
Let $\Sigma\in\mc{M}$, let $\sigma\in\Coker(D_\Sigma)$ be
nonvanishing, let $z\in\Sigma$, and suppose that $p(z)$ has distance
at least $R$ from all vertices in $\tau(\Sigma)$.  Then
\begin{equation}
\label{eqn:approx1}
|\sigma(z) - \Pi_W\sigma(z)| < \varepsilon_0|\Pi_W\sigma(z)|.
\end{equation}
\end{proposition}


\subsection{Orientation of the obstruction bundle}
\label{sec:OBO}

We now specify an orientation of the obstruction bundle
$\mc{O}\to\mc{M}$ associated to an elliptic Reeb orbit.  This will be
needed later to define various signs.

In the special case $S(t)=\theta$, the operators $L_m$ and $D_\Sigma$
are complex linear, and so the real vector bundle $\mc{O}^*$ has a
canonical orientation, which determines an orientation of the dual
real vector bundle $\mc{O}$.

To orient the obstruction bundle in the general case, note that for
any elliptic Reeb orbit, we can deform $J$ to a different admissible
almost complex structure $J'$ for which $S(t)=\phi(t)$ for some
function $\phi:S^1\to\R$.  One can then deform the trivialization of
$\alpha^*\xi$ to arrange that $\phi(t)=\theta$.  This process gives
rise to a continuous family of Fredholm operators $D_\Sigma$, cf.\
\cite{fh}.  Since the space of admissible almost complex structures
and trivializations of $\alpha^*\xi$ is contractible, we obtain a
canonical bijection between orientations of the obstruction bundle for
a general $S(t)$ and orientations of the obstruction bundle in the
special case $S(t)=\theta$.

\subsection{A compactification of $\mc{M}/\R$}
\label{sec:cmr}

As in \S\ref{sec:GP}, $\R$ acts on
$\mc{M}=\mc{M}(a_1,\ldots,a_{N_+}\mid a_{-1},\ldots,a_{-N_-})$ by
translating the $\R$ coordinate on $\R\times S^1$.  Given
$\Sigma\in\mc{M}$, let $[\Sigma]$ denote the equivalence class of
$\Sigma$ in $\mc{M}/\R$.  We now define a compactification of
$\mc{M}/\R$, which is slightly different from the symplectic field
theory compactification in \cite{behwz,egh}.  This compactification
will be used in the analysis in \S\ref{sec:PAR} and in \cite{obg2}.

\begin{definition}
\label{def:compactification}
An element of $\overline{\mc{M}/\R}$ is a tuple
$(T;[\Sigma_{*1}],\ldots,[\Sigma_{*p}])$ where:
\begin{itemize}
\item
$T$ is an oriented weighted tree in $T(a_1,\ldots,a_{N_+}\mid
a_{-1},\ldots,a_{-N_-})$ with $p$ internal vertices, and with
orderings of the edges and internal vertices such that the edge
ordering restricts to the given orderings of the positive and negative
leaves.
\item
Let $m_{j,1},m_{j,2},\ldots$ denote the multiplicities of the outgoing
edges of the $j^{th}$ internal vertex, and let
$n_{j,1},n_{j,2},\ldots$ denote the multiplicities of the incoming
edges of the $j^{th}$ internal vertex, in their given order.  Then
\[
\Sigma_{*j}\in\mc{M}^{(j)}\eqdef \mc{M}(m_{j,1},m_{j,2},\ldots \mid
n_{j,1},n_{j,2},\ldots).
\]
\end{itemize}
Two such tuples are equivalent if they differ by the following
operations:
\begin{itemize}
\item
Reordering the edges and internal vertices.
\item
For an internal edge $e$ from vertex $j'$ to vertex $j$, acting on the
asymptotic markings of the corresponding positive end of $\Sigma_{*j'}$
and negative end of $\Sigma_{*j}$ by the same element of $\Z/m(e)$.
\end{itemize}
There is an inclusion $\mc{M}/\R\to\overline{\mc{M}/\R}$ sending
$[\Sigma]\mapsto (T,[\Sigma])$, where $T$ has only one internal
vertex.
\end{definition}

\begin{definition}
\label{def:convergence}
A sequence $\{[\pi_k:\Sigma_k\to\R\times S^1]\}_{k=1,2,\ldots}$ in
$\mc{M}/\R$ converges to
$(T;[\Sigma_{*1}],\ldots,[\Sigma_{*p}])\in\overline{\mc{M}/\R}$ if for
all $k$ sufficiently large, there are disjoint closed subsets
$\Sigma_{k1},\ldots,\Sigma_{kp}\subset\Sigma_k$ such that:
\begin{description}
\item{(a)}
Each ramification point in $\Sigma_k$ is contained in some
$\Sigma_{kj}$.
\item{(b)}
Each $\Sigma_{kj}$ is a component of
the $\pi_k$-inverse image of a cylinder in $\R\times S^1$, and the length of
$\pi_k(\Sigma_{kj})$ goes to infinity as $k\to\infty$.
\item{(c)}
Let $s_{kj}$ denote the $s$ coordinate of the central circle of
$\pi_k(\Sigma_{kj})$.  Then the function $\pi_k^*s - s_{kj}$ on the set
of ramification points in $\Sigma_{kj}$ has a $k$-independent upper
bound.
\item{(d)}
$T$ is obtained from $\tau(\Sigma_k)$ by for each $j\in\{1,\ldots,p\}$
collapsing the vertices corresponding to the ramification points in
$\Sigma_{kj}$ to a single vertex.
\item{(e)}
For each internal edge $e$ of $T$, there is an identification
$\Phi_{k,e}$ of the corresponding cylinder in $\Sigma_k$ with an
interval cross $\R/2\pi m(e)\Z$, commuting with the projections to
$\R\times S^1$, such that the following holds.  For each internal
vertex $j$ of $T$, let $\widehat{\Sigma}_{kj}\in \mc{M}^{(j)}$ be
obtained by attaching half-infinite cylinders to the boundary circles
of $\Sigma_{kj}$, and asymptotically marking the ends corresponding to
internal edges using the identifications $\Phi_{k,e}$.  Let
$T_{-s_{kj}}$ denote the translation $s\mapsto s - s_{kj}$.  Choose
the representatives $\Sigma_{*j}$ so that their ramification points
are centered at $s=0$. Then
$\lim_{k\to\infty}T_{-s_{kj}}(\widehat{\Sigma}_{kj})=\Sigma_{*j}$ in
$\mc{M}^{(j)}$.
\end{description}
\end{definition}

\begin{lemma}
\label{lem:convergence}
Any sequence in $\mc{M}/\R$ has a subsequence which converges in
$\overline{\mc{M}/\R}$.
\end{lemma}

\begin{proof}
Let $\{[\pi_k:\Sigma_k\to\R\times S^1]\}_{k=1,2,\ldots}$ be a sequence
in $\mc{M}/\R$.  The number of ramification points in $\Sigma_k$
counted with multiplicity is $N-2$, which is independent of
$k$.  Hence we can pass to a subsequence so that there is an integer
$p\in\{1,\ldots,N-2\}$, and for each $k$ a partition of the
ramification points in $\Sigma_k$ into subsets
$\{\Lambda_{kj}\}_{j=1,\ldots,p}$, such that:
\begin{itemize}
\item
The diameter in $\Sigma_k$ of $\Lambda_{kj}$ has a $k$-independent
upper bound.
\item
If $j\neq j'$, then the distance in $\Sigma_k$ between $\Lambda_{kj}$
and $\Lambda_{kj'}$ is greater than $k$.
\end{itemize}
We can then refine the sequence so that there are disjoint closed subsets
$\Sigma_{k1},\ldots,\Sigma_{kp}\subset\Sigma_k$ satisfying conditions
(a)--(c) in Definition~\ref{def:convergence}.

To keep track of the combinatorics of the $\Sigma_{kj}$'s, for each
$k$ define an oriented weighted tree $T_k\in T(a_1,\ldots,a_{N_+}\mid
a_{-1},\ldots,a_{-N_+})$ as follows.  The internal vertices of $T_k$
are labeled by $1,\ldots,p$ and identified with
$\Sigma_{k1},\ldots,\Sigma_{kp}$.  Edges correspond to components of
$\Sigma_k\setminus \bigcup_k\Sigma_{kj}$, and the multiplicity of an
edge is the degree of $\pi_k$ on the corresponding cylinder.  The
internal (resp.\ external) vertices at the endpoints of an edge are
determined by the boundary circles (resp.\ ends) of the associated
cylinder.  Edges are oriented in the increasing $s$ direction.

Since there are only finitely many oriented weighted trees in the set
$T(a_1,\ldots,a_{N_+}\mid a_{-1},\ldots,a_{-N_+})$, we can pass to a
subsequence such that all the trees $T_k$ are isomorphic, preserving
the leaf labels, to a single oriented weighted tree $T$.  Note that
these isomorphisms $T_k\simeq T$ are canonical since a tree has no
nontrivial automorphisms fixing the leaves.  Choose orderings of the
edges and internal vertices of $T$ as in
Definition~\ref{def:compactification}, and use the vertex ordering to
order $\Sigma_{k1},\ldots,\Sigma_{kp}$.  We have now achieved
requirement (d) in Definition~\ref{def:convergence}.

Next, for each internal edge $e$ of $T$, fix an identification
$\Phi_{k,e}$ of the corresponding cylinder in $\Sigma_k$ with an
interval cross $\R/2\pi m(e)\Z$, commuting with the projections to
$\R\times S^1$.  Then $\widehat{\Sigma}_{kj}\in\mc{M}^{(j)}$ is
defined.  Moreover, since condition (c) in
Definition~\ref{def:convergence} holds, there is a $k$-independent
constant $R$ such that each $T_{-s_{kj}}(\widehat{\Sigma}_{kj})$ is in
$\mc{M}_R(m_{j,1},m_{j,2},\ldots\mid n_{j,1},n_{j,2},\ldots)$.  
By Lemma~\ref{lem:LFI}, we have
\[
\kappa_\theta(m_{j,1},m_{j,2},\ldots\mid n_{j,1},n_{j,2},\ldots) = 1.
\]
Hence we can apply Lemma~\ref{lem:mrc} to find a further subsequence
satisfying condition (e) in Definition~\ref{def:convergence}.
\end{proof}

\section{The linearized section of the obstruction bundle}
\label{sec:section}

In this section, fix an admissible almost complex structure $J$ on
$\R\times Y$, an elliptic Reeb orbit $\alpha$ with monodromy angle
$\theta$, and
\begin{equation}
\label{eqn:fixS}
S = (a_1,\ldots,a_{\overline{N}_+};a_{\overline{N}_++1},\ldots,a_{N_+} \mid
a_{-1},\ldots,a_{\overline{N}_-};a_{-(\overline{N}_-+1)},\ldots,a_{-N_-})
\end{equation}
satisfying the sum condition \eqref{eqn:sumCondition}.
Assume as usual that $N\eqdef N_++N_->2$.  Assume also that
$\kappa_\theta(S)=1$.  Finally, we make:
\begin{assumption}[Partition minimality]
\label{asm:pm}
Under the partial order $\ge_\theta$ in Definition~\ref{def:newpo},
the partition $(a_{\overline{N}_++1}, \ldots,a_{N_+})$ is minimal, and
the partition $(a_{-(\overline{N}_-+1)},\ldots,a_{-N_-})$ is maximal.
\end{assumption}
By \S\ref{sec:obstructionBundle},
there is an obstruction bundle $\mc{O}$ over the moduli space of
branched covers $\mc{M}=\mc{M}(a_1,\ldots,a_{N_+}\mid
a_{-1},\ldots,a_{-N_-})$ with $\text{rank}(\mc{O}) =
\dim(\mc{M})$.

We will now define a section $\mathfrak{s}_0$ of $\mc{O}$ over $\mc{M}_R$
(see Definition~\ref{def:M_R}).  We will then quote results from
\cite{obg2} which relate an appropriate count of zeroes of
$\mathfrak{s}_0$ to the count of gluings in Theorem~\ref{thm:main}.

\subsection{Definition of the linearized section $\mathfrak{s}_0$}
\label{sec:s0}

For $i=1,\ldots,\overline{N}_+$, let $\lambda_i$ denote the largest
negative eigenvalue of the operator $L_{a_i}$.  Likewise, for
$i=-1,\ldots,-\overline{N}_-$, let $\lambda_i$ denote the smallest
positive eigenvalue of the operator $L_{a_i}$.  Let $\mc{B}_i$ denote
the $\lambda_i$-eigenspace of $L_{a_i}$, and let $\mc{B}$ denote the
direct sum of the $\mc{B}_i$'s for $i=1,\ldots,\overline{N}_+$ and
$i=-1,\ldots,-\overline{N}_-$.  Given $\Sigma\in\mc{M}$ and
$\sigma\in\Coker(D_\Sigma)$, over the $i^{th}$ end of $\Sigma$
(identified with $[R,\infty)\times \R/2\pi a_i\Z$ or
$(-\infty,-R]\times \R/2\pi a_i\Z$ via the asymptotic marking), write
$\sigma=\sigma_i(s,t)d\zbar$.

\begin{definition}
\label{def:s0}
Fix $\gamma=\{\gamma_i\}\in\mc{B}$ and real numbers $R,r>0$.  Define
the {\em linearized section\/} $\mathfrak{s}_0$ of the obstruction bundle
$\mc{O}\to\mc{M}_R$ as follows.  If $\Sigma\in\mc{M}_R$ and
$\sigma\in\Coker(D_\Sigma)$, then
\begin{equation}
\label{eqn:s0}
\mathfrak{s}_0(\Sigma)(\sigma) \eqdef \sum_{i=1}^{\overline{N}_+}
\langle \gamma_i,\sigma_i(R+r,\cdot)\rangle
-
\sum_{i=-1}^{-\overline{N}_-}
\langle\gamma_i, \sigma_i(-R-r,\cdot)\rangle
 \in \R.
\end{equation}
Here the brackets denote the real inner product on $L^2(\R/2\pi
a_i\Z,\R^2)$.
\end{definition}

We will study the linearized section under the assumption that the
eigenvectors $\gamma$ are ``admissible'' in the following sense.  If
$i,j\in\{1,\ldots,\overline{N}_+\}$ and
$\ceil{a_i\theta}/a_i=\ceil{a_j\theta}/a_j$, then we can identify
$\mc{B}_i=\mc{B}_j$ using coverings as follows.  Let
$\phi\eqdef\ceil{a_i\theta}/a_i$, so that
$\lambda_i=\lambda_j=E_\phi^+\fedqe\lambda$.  Write $\phi=\eta/m$
where $\eta,m\in\Z$ and $m>0$ is as small as possible; then the
integers $a_i$ and $a_j$ are both divisible by $m$.  For any positive
integer $d$, we can pull back an eigenfunction of $L_m$ with
eigenvalue $\lambda$ to an eigenfunction of $L_{md}$ with the same
eigenvalue.  In this way we identify the $\lambda$-eigenspaces of
$L_{md}$ for different $d$ with each other.  Likewise, if
$i,j\in\{-1,\ldots,-\overline{N}_-\}$, then we can identify
$\mc{B}_i=\mc{B}_j$ when
$\floor{a_i\theta}/a_i=\floor{a_j\theta}/a_j$.

Note also that the cyclic group $\Z/a_i$ acts linearly on the
eigenspace $\mc{B}_i$, via pullback from its action on $\R/2\pi a_i\Z$
by deck transformations of the covering map to $S^1=\R/2\pi\Z$.

\begin{definition}
\label{def:aa}
An element $\gamma=\{\gamma_i\}\in\mc{B}$ is {\em admissible\/} if:
\begin{description}
\item{(a)} $\gamma_i\neq 0$ for all $i=1,\ldots,\overline{N}_+$ and
$i=-1,\ldots,-\overline{N}_-$.
\item{(b)}
If $i,j\in\{1,\ldots,\overline{N}_+\}$ and
$\ceil{a_i\theta}/a_i=\ceil{a_j\theta}/a_j$, or
$i,j\in\{-1,\ldots,-\overline{N}_-\}$ and
$\floor{a_i\theta}/a_i=\floor{a_j\theta}/a_j$, then for all
$g_i\in\Z/a_i$ and $g_j\in\Z/a_j$ we have $g_i\cdot\gamma_i\neq
g_j\cdot \gamma_j$ in $\mc{B}_i=\mc{B}_j$.
\end{description}
\end{definition}

\subsection{Counting zeroes of $\mathfrak{s}_0$}
\label{sec:REC}

We now want to count zeroes of $\mathfrak{s}_0$, for which purpose we will use
the following formalism.

\begin{definition}
Let $\psi$ be a section of $\mc{O}$ over $\mc{M}_R$.  Suppose
that $\psi$ is nonvanishing on
$\partial\mc{M}_R$.  Then define the {\em
relative Euler class\/}
\[
e(\mc{O}\to\mc{M}_R,\psi)\in\Z
\]
as follows.
Let $\psi'$ be a section of $\mc{O}$ over $\mc{M}_R$ such that
$\psi=\psi'$ on $\partial\mc{M}_R$, and all zeroes of
$\psi'$ are nondegenerate.  Define $e(\mc{O}\to\mc{M}_R,\psi)$
to be the signed count of zeroes of $\psi'$, where the signs are determined
using the orientation of $\mc{M}_R$ as a complex
manifold and the orientation of $\mc{O}$ defined in
\S\ref{sec:OBO}.  By Lemma~\ref{lem:mrc}, this count is well-defined
and depends only on $\psi|_{\partial\mc{M}_R}$.  We usually denote
this relative Euler class by $\#\psi^{-1}(0)$, even though the zeroes
of $\psi$ itself may be degenerate.
\end{definition}

Note that if $\{\psi_t\mid t\in[0,1]\}$ is a homotopy of sections
of $\mc{O}\to\mc{M}_R$ such that $\psi_t|_{\partial\mc{M}_R}$ is
nonvanishing for all $t\in[0,1]$, then $\#\psi_0^{-1}(0) =
\#\psi_1^{-1}(0)$.

The following lemma gives sufficient conditions for $\mathfrak{s}_0$
to be nonvanishing on $\partial\mc{M}_R$.  To state it, let
$\lambda\eqdef\min\{|\lambda_i|\}$ and
$\Lambda\eqdef\max\{|\lambda_i|\}$.

\begin{lemma}
\label{lem:zab}
Given admissible $\gamma=\{\gamma_i\}\in\mc{B}$, if $r$ is
sufficiently large with respect to $\gamma$, and if $R>3\Lambda
r/\lambda$, then $\mathfrak{s}_0$ has no zeroes on
$\mc{M}_R\setminus\mc{M}_{R-r}$.
\end{lemma}

\begin{proof}
This follows from \cite[Prop.\ 8.2]{obg2}, as in the proof of
\cite[Cor.\ 8.6]{obg2}.
\end{proof}

It follows that if $\gamma$, $r$, and $R$ satisfy the hypotheses of
Lemma~\ref{lem:zab}, then the relative Euler class
$\#\mathfrak{s}_0^{-1}(0)$ is defined.  Moreover:

\begin{corollary}
\label{cor:REC}
Under the assumptions of Lemma~\ref{lem:zab}, the relative Euler class
$\#\mathfrak{s}_0^{-1}(0)$ depends only on $S$ and $\theta$.
\end{corollary}

\begin{proof}
As in \S\ref{sec:OBO}, we can deform the almost complex structure $J$
on $\R\times Y$ and the trivialization of $\alpha^*\xi$ so as to
deform the operator $D_\Sigma$ to a complex linear operator depending
only on $S$ and $\theta$, in which the operators $L_{a_i}$ are also
complex linear.  We can simultaneously deform the collection of
eigenvectors $\gamma$ while preserving the admissibility conditions.
Once the operators $L_{a_i}$ are complex linear, the set of admissible
$\gamma$ is connected, because now the conditions $\gamma_i=0$ and
$g_i\cdot\gamma_i=g_j\cdot\gamma_j$ have real codimension $2$.  Thus
all of the different versions of the section $\mathfrak{s}_0$ are
homotopic.  For any given homotopy of the almost complex structure,
the trivialization, and the admissible $\gamma$, if we fix $r$ and $R$
sufficiently large, then Lemma~\ref{lem:zab} will apply throughout the
homotopy so that the count of zeroes does not change.  Increasing $r$
and $R$ does not change the count of zeroes for the same reason.
\end{proof}

\subsection{Zeroes of $\mathfrak{s}_0$ and gluings}

The significance of $\#\mathfrak{s}_0^{-1}(0)$ for the gluing story is as
follows.  Let $(u_+,u_-)$ be a gluing pair satisfying conditions (i)
and (ii) in \S\ref{sec:overview}.  Order the negative ends of $u_+$
and the positive ends of $u_-$ such that the $\R$-invariant negative
ends of $u_+$ are those labeled by $\overline{N}_++1,\ldots,N_+$, and
the $\R$-invariant positive ends of $u_-$ are those labeled by
$-(\overline{N}_-+1),\ldots,-N_-$.  The main result of \cite{obg2} can
be stated as follows:

\begin{theorem}
\label{thm:sequel}
\cite[Thm.\ 1.1]{obg2}
Fix coherent orientations and generic $J$.  If $(u_+,u_-)$ is a gluing
pair as above, then
\[
\#G(u_+,u_-)=\epsilon(u_+)\epsilon(u_-)\#\mathfrak{s}_0^{-1}(0),
\]
where $\#\mathfrak{s}_0^{-1}(0)$ is defined as in Corollary~\ref{cor:REC}.
\end{theorem}

A few words are in order concerning the proof of this theorem.  As
sketched in
\S\ref{sec:overview}, we have
\[
\#G(u_+,u_-)=\epsilon(u_+)\epsilon(u_-)\#\mathfrak{s}^{-1}(0),
\]
where $\mathfrak{s}$ is a section of $\mc{O}$ over $\mc{M}_R$ arising from
the gluing analysis.  The idea is to relate the section $\mathfrak{s}$ to
the linearized section $\mathfrak{s}_0$, where the latter is defined using
certain eigenfunctions $\gamma\in\mc{B}$ that are determined by the
asymptotics of the negative ends of $u_+$ and the positive ends of
$u_-$.  If $J$ is generic, then the asymptotic eigenfunctions $\gamma$
are admissible.  In this case a generalization of Lemma~\ref{lem:zab}
shows that the homotopy of sections $\mathfrak{s}_t\eqdef
t\mathfrak{s}+(1-t)\mathfrak{s}_0$ has no zeroes on $\partial\mc{M}_R$ for
$t\in[0,1]$, so that $\#\mathfrak{s}^{-1}(0)=\#\mathfrak{s}_0^{-1}(0)$.

\subsection{Consequences of partition minimality}

By Theorem~\ref{thm:sequel}, to prove Theorem~\ref{thm:main} we just need
to compute $\#\mathfrak{s}_0^{-1}(0)$.  Before doing so, we note two basic
facts about the numbers $\overline{N}_+$ and $\overline{N}_-$ in
\eqref{eqn:fixS}.

\begin{lemma}
\label{lem:cpm}
The assumption $\kappa_\theta(S)=1$ and Assumption~\ref{asm:pm} imply that
$\overline{N}_+\ge N_+-1$ and $\overline{N}_-\ge N_--1$.
\end{lemma}

\begin{proof}
By symmetry it is enough to show that $\overline{N}_+\ge N_+-1$.
Suppose to the contrary that $N_+-\overline{N}_+\ge 2$.  Then we can
construct a tree $T\in T(S)$ which has a splitting vertex $v$ with one
incoming edge and with outgoing edges incident to the positive leaves
labeled by $N_+-1$ and $N_+$.  Then by Lemma~\ref{lem:LFI},
\[
\ceil{a_{N_+-1}\theta} + \ceil{a_{N_+}\theta} =
\ceil{(a_{N_+-1}+a_{N_+})\theta}.
\]
It follows from this that
\[
(a_{\overline{N}_++1},\ldots,a_{N_+}) >_\theta
(a_{\overline{N}_++1},\ldots,a_{N_+-2},a_{N_+-1}+a_{N_+}).
\]
This contradicts the assertion in Assumption~\ref{asm:pm} that
$(a_{\overline{N}_++1},\ldots,a_{N_+})$ is minimal.
\end{proof}

\begin{lemma}
\label{lem:OTE}
If $\gamma$ is admissible and $\mathfrak{s}_0(\Sigma)=0$, then
$\overline{N}_+=N_+$ or $\overline{N}_-=N_-$.
\end{lemma}

\begin{proof}
If not, then by Lemma~\ref{lem:cpm}, $\overline{N}_+=N_+-1$ and
$\overline{N}_-=N_--1$.  Since we are assuming that $N>2$, there
is an end labeled by $i\notin\{N_+,-N_-\}$.  Since $\gamma$ is
admissible, $\gamma_i\neq 0$.  By Lemma~\ref{lem:SCE}(a), the
projection $V_{i,N_+,-N_-}\to\mc{B}_i$ is surjective.  Hence we can
find a special cokernel element $\sigma\in V_{i,N_+,-N_-}$ such that the
$\sigma_i$ term in \eqref{eqn:s0} is nonzero.  Then
$\mathfrak{s}_0(\Sigma)(\sigma)\neq 0$, because all other terms on the
r.h.s. of \eqref{eqn:s0} are zero.
\end{proof}

\begin{corollary}
\label{cor:trivialCase}
If $\overline{N}_+\neq N_+$ and $\overline{N}_-\neq N_-$, then
\begin{equation}
\label{eqn:trivialCase}
\#\mathfrak{s}_0^{-1}(0) = c_\theta(S).
\end{equation}
\end{corollary}

\begin{proof}
The left side of \eqref{eqn:trivialCase} is zero by
Lemma~\ref{lem:OTE}, while the right hand side of
\eqref{eqn:trivialCase} is zero by Definition~\ref{def:OTE} and our
assumption that $N>2$.
\end{proof}

\section{Combinatorics of the elliptic gluing coefficients}
\label{sec:trees}

This section gives a diagrammatic reinterpretation of the
combinatorial gluing coefficient $c_\theta(S)$ from
\S\ref{sec:ctheta}.  Recall that when $\kappa_\theta(S)=1$, the
integer $c_\theta(S)$ is determined by an auxiliary function
$f_\theta$, evaluated on a certain reordering of the list $S$.  The
main result of this section, Lemma~\ref{lem:IFT} below, expresses
$f_\theta(S)$ as a sum, over ``admissible'' trivalent trees with
``edge pairings'', of certain positive integer weights.  (A
straightforward extension of this interprets $c_\theta(S)$ when
$\kappa_\theta(S)>1$ as a sum over forests.)  This alternate
definition is lengthier but more symmetric, and leads to a proof that
$c_\theta(S)$ satisfies the symmetry property \eqref{eqn:cthetasym}.
The combinatorics introduced here will be used in the computation of
$\#\mathfrak{s}_0^{-1}(0)$ in
\S\ref{sec:count}.

For the rest of this section fix $S=(a_1,\ldots,a_{N_+} \mid
a_{-1},\ldots,a_{-N_-})$, where $a_1,\ldots,a_{N_+}$ and
$a_{-1},\ldots,a_{-N_-}$ are positive integers satisfying
\eqref{eqn:M}.

\subsection{Edge pairings, weights, and admissible trees}
\label{sec:ATEP}

We begin with some combinatorial definitions.  

\begin{definition}
\label{def:EP}
Let $T$ be a trivalent tree.  An {\em edge pairing\/} $P$ on $T$ is an
assignment, to each internal vertex $v$ of $T$, of two distinct edges
$e_v^+$ and $e_v^-$ incident to $v$, such that the sets
$\{e_v^+,e_v^-\}$ and $\{e_w^+,e_w^-\}$ are disjoint whenever $v$ and
$w$ are adjacent (and hence whenever $v$ and $w$ are distinct)
internal vertices.
\end{definition}

Note that in a trivalent tree, the number of edges equals twice the
number of internal vertices plus one. Hence for any edge pairing $P$
on $T$, there is a distinguished edge $e_0$ which is not one of the
edges $e_v^+$ or $e_v^-$ for any internal vertex $v$.

\begin{definition}
\label{def:weight}
If $T$ is an oriented weighted trivalent tree, if $P$ is an edge
pairing on $T$, and if $\theta$ is an irrational number, define the
{\em weight\/} $W_\theta(T,P)\in\Z^{>0}$ as follows.  For an internal
vertex $v$, define $m_v^+$ to be ``the outward flow along the edge
$e_v^+$'', namely $m(e_v^+)$ if $e_v^+$ points outward from $v$, and
$-m(e_v^+)$ if $e_v^+$ points inward toward $v$.  Similarly define
$m_v^-$ to be ``the inward flow along the edge $e_v^-$'', namely
$m(e_v^-)$ if $e_v^-$ points inward toward $v$, and $-m(e_v^-)$ if $e_v^-$
points outward from $v$.  Then
\begin{equation}
\label{eqn:weight}
W_\theta(T,P) \eqdef m(e_0) \prod_{v\in\dot{V}(T)} \Big(m_v^-
\ceil{m_v^+\theta}
- m_v^+\floor{m_v^-\theta}\Big).
\end{equation}
\end{definition}

The above combinatorial notions arise naturally in the analysis, as
explained in Remark~\ref{rem:EPN}.  

Next, recall that $T(S)$ denotes the set of oriented weighted trees
whose positive leaves are labeled by $1,\ldots,N_+$ and whose negative
leaves are labeled by $-1,\ldots,-N_-$, such that the $i^{th}$ leaf
has multiplicity $a_i$.

\begin{figure}
\begin{center}
\scalebox{0.6}{\includegraphics{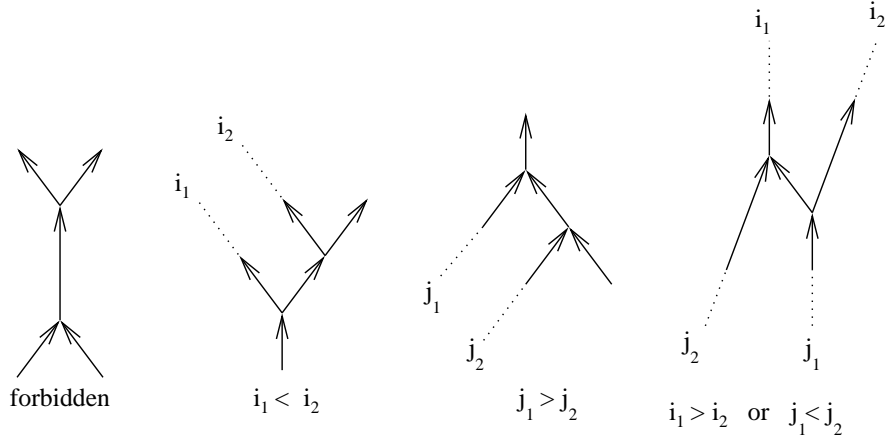}}
\end{center}
\caption{The rules for an admissible tree.}
\label{fig:admissible}
\end{figure}

\begin{definition}
A trivalent tree $T\in T(S)$ is {\em admissible\/} if the following
conditions hold (see Figure~\ref{fig:admissible}):
\begin{description}
\item{(a)}
No oriented edge starts at a joining vertex and ends at a splitting
vertex.
\item{(b)}
Let $v$ be a splitting vertex with outgoing edges $e_1$ and $e_2$.
Suppose there is an upward path starting along $e_1$ and ending at the
positive leaf ${i_1}$.  Suppose $e_2$ is incident to another splitting
vertex $w$, from which there is an upward path leading to the positive
leaf ${i_2}$.  Then $i_1<i_2$.
\item{(c)}
Symmetrically to (b), let $v$ be a joining vertex with incoming edges
$e_1$ and $e_2$.  Suppose there is a downward path starting along
$e_1$ and ending at the negative leaf ${j_1}$.  Suppose $e_2$ is
incident to another joining vertex $w$, from which there is a downward
path leading to the negative leaf ${j_2}$.  Then $j_1 > j_2$.
\item{(d)}
Let $e$ be an edge from a splitting vertex $w$ to a joining vertex
$v$.  Suppose there is an upward path from $v$ to the positive leaf
${i_1}$ and a downward path from $w$ to the negative leaf ${j_1}$.
Suppose there is an upward path from $w$, not containing $e$, to the
positive leaf ${i_2}$, and a downward path from $v$, not containing
$e$, to the negative leaf ${j_2}$.  Then $i_1>i_2$ or $j_1 < j_2$ (or
both).
\end{description}
\end{definition}

\begin{figure}
\begin{center}
\scalebox{0.6}{\includegraphics{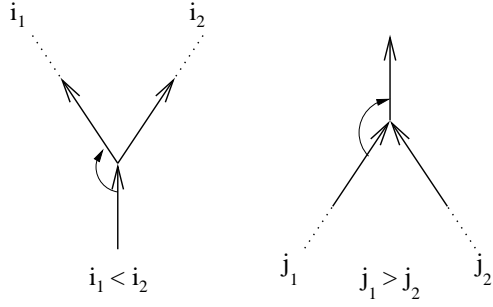}}
\end{center}
\caption{The canonical edge pairing for an admissible tree.  The
symbolism is that at each internal vertex $v$, the small arc starts at
$e_v^-$ and ends at $e_v^+$.}
\label{fig:edgepairing}
\end{figure}

\begin{lemma}
If $T$ is admissible, then there is a unique edge pairing $P_T$ on $T$
such that (see Figure~\ref{fig:edgepairing}):
\begin{description}
\item{(i)}
If $v$ is a splitting vertex, then $e_v^-$ is the incoming edge;
there is a unique upward path starting along $e_v^+$, say to the
positive leaf $i_1$; and if $i_2$ is a positive leaf reached by an
upward path starting along the other outgoing edge, then $i_1<i_2$.
\item{(ii)}
Symmetrically to (i), if $v$ is a joining vertex, then $e_v^+$ is the
outgoing edge; there is a unique downward path starting along $e_v^-$,
say to the negative leaf $j_1$; and if $j_2$ is a negative leaf
reached by a downward path starting along the other incoming edge,
then $j_1 > j_2$.
\end{description}
\end{lemma}

\begin{proof}
Let $v$ be a splitting vertex.  By admissibility condition (b), at
least one outgoing edge of $v$ is incident to a joining vertex or to a
positive leaf.  By admissibility condition (a), there is a unique
upward path starting along such an edge.  By condition (b) again, there is a
unique such edge $e_v^+$ satisfying condition (i) above.

Symmetrically, if $v$ is a joining vertex, then there is a unique
incoming edge $e_v^-$ satisfying condition (ii) above.

To see that $P_T$ satisfies the disjointness condition in the
definition of edge pairing, we must check that an edge $e$ from an
internal vertex $v$ to an internal vertex $w$ cannot be in both of the
sets $\{e_v^+,e_v^-\}$ and $\{e_w^+,e_w^-\}$.  There are four cases to
consider, depending on whether the vertices $v$ and $w$ are joining or
splitting.  These four cases are precisely covered by admissibility
conditions (a)--(d).
\end{proof}

Let $\mc{A}(S)$ denote the set of admissible trees in $T(S)$.  We can
now interpret $f_\theta(S)$ as a sum over admissible trees of the
weights associated to their canonical edge pairings:

\begin{lemma}
\label{lem:IFT}
\[
f_\theta(S) = \sum_{T\in
\mc{A}(S)} W_\theta(T,P_T).
\]
\end{lemma}

The proof of this lemma is postponed to \S\ref{sec:EAT}.  We can now
prove the symmetry property \eqref{eqn:cthetasym}.  By the definition
of $c_\theta$, it is enough to show:

\begin{corollary}
\label{cor:fsym1}
The function $f_\theta$ satisfies the symmetry property
\begin{equation}
\label{eqn:FS}
f_\theta(a_1,\ldots,a_{N_+} \mid a_{-1},\ldots,a_{-N_-}) =
f_{-\theta}(a_{-1},\ldots,a_{-N_-} \mid a_1,\ldots,a_{N_+}).
\end{equation}
\end{corollary}

\begin{proof}
The definition of admissible tree is
symmetric, in that reversing edge orientations defines a bijection
\[
\phi:\mc{A}(a_1,\ldots,a_{N_+} \mid a_{-1},\ldots,a_{-N_-}) \longrightarrow
\mc{A}(a_{-1},\ldots,a_{-N_-} \mid a_1,\ldots,a_{N_+}).
\]
The canonical edge pairing is also symmetric, in that $P_{\phi(T)}$
is obtained from $P_T$ by switching $e_v^+$ and $e_v^-$ for all
$v$.  It now follows from equation
\eqref{eqn:weight}, using the identity $-\ceil{x}=\floor{-x}$, that
\[
W_{-\theta}(\phi(T),P_{\phi(T)}) = W_{\theta}(T,P_T).
\]
Equation \eqref{eqn:FS} then follows from Lemma~\ref{lem:IFT}.
\end{proof}


\begin{remark}
\label{rem:IHX}
Assume that $N_+>1$ and $\kappa_\theta(S)=1$.  Let $S'$ be obtained
from switching $a_1$ and $a_2$:
\[
S' \eqdef (a_2,a_1,a_3,\ldots,a_{N_+} \mid a_{-1},\ldots,a_{-N_-}).
\]
Let $\widehat{S}$ be obtained from $S$ by adding the first two
entries:
\[
\widehat{S} \eqdef (a_1+a_2,a_3,\ldots,a_{N_+} \mid a_{-1},\ldots,a_{-N_-}).
\]
One can then show that
\begin{equation}
\label{eqn:IHX}
f_\theta(S) - f_\theta(S') = (a_2\ceil{a_1\theta} -
a_1\ceil{a_2\theta})\cdot f_\theta(\widehat{S}).
\end{equation}
Equation \eqref{eqn:IHX} implies that $c_\theta(S)$ does not depend on
the choice of reordering satisfying \eqref{eqn:reordering}, as
discussed in Remark~\ref{rem:notobvious}.  To prove \eqref{eqn:IHX},
one can first show that when $\kappa_\theta(S)=1$, the function
$W_\theta$ on trees with edge pairings satisfies a version of the IHX
relation.  (For the usual IHX relation see e.g.\ \cite{dbn}.)  One can
then expand both sides of equation
\eqref{eqn:IHX} using Lemma~\ref{lem:IFT} and show that the difference
is a linear combination of IHX relations.  We omit the details,
because we will give a different proof that $c_\theta(S)$ is
well-defined in \S\ref{sec:count}.
\end{remark}

\subsection{Enumerating admissible trees}
\label{sec:EAT}

To prove Lemma~\ref{lem:IFT}, we now introduce a way to enumerate
admissible trees, which is less symmetric, but more closely related to
the definition of $f_\theta$ and to the obstruction bundle
calculations in \S\ref{sec:count}.

\begin{definition}
\label{def:EAT}
Let $\mc{E}(S)$ denote the
set of $N_+$-tuples $(E_1,\ldots,E_{N_+})$ of nonempty subsets of
$\{-1,\ldots,-N_-\}$, such that:
\begin{description}
\item{(i)} For $i=1,\ldots,N_+-1$, let
$E_i'\eqdef E_i\setminus\{\min(E_i)\}$.  Then
\begin{equation}
\label{eqn:EC1}
\{-1,\ldots,-N_-\} = E_{N_+}\sqcup\bigsqcup_{i=1}^{N_+-1} E_i'.
\end{equation}
\item{(ii)}
If $N_+>1$, then
\begin{equation}
\label{eqn:EC2}
\sum_{j\in E_{1}'} a_j < 
a_{1} <  \sum_{j\in E_{1}} a_j.
\end{equation}
Also, let
\[
\overline{S} \eqdef (a_2,\ldots,a_{N_+} \mid a_{-E_1})
\]
where $a_{-E_1}$ denotes the arguments $a_j$ for $0>j\notin E_1$ arranged
in decreasing order, with one additional argument equal to $\sum_{j\in
E_1}a_j-a_1$,
inserted into the position that $\min(E_1)$ would occupy in the order.
Let
\[
\xi:\{-1,\ldots,-N_-\}\setminus E_1'\longrightarrow\{-1,\ldots,-(N_--|E_1'|)\}
\]
denote the order-preserving bijection.  Then
\begin{equation}
\label{eqn:EC3}
\overline{E} \eqdef (\xi(E_2),\ldots,\xi(E_{N_+})) \in \mc{E}(\overline{S}).
\end{equation}
\end{description}
\end{definition}

\begin{figure}
\begin{center}
\scalebox{0.6}{\includegraphics{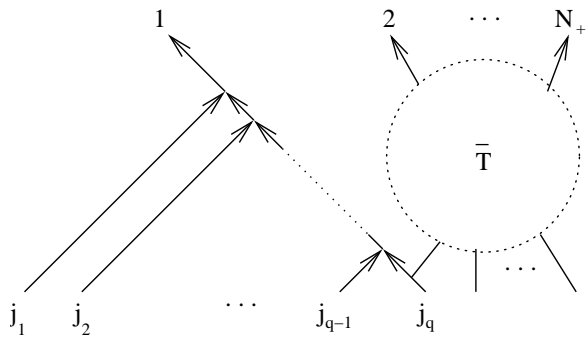}}
\end{center}
\caption{The inductive construction of an admissible tree when
$N_+>1$.  Here $E_1=\{j_1>\cdots>j_q\}$.}
\label{fig:induction}
\end{figure}

\begin{definition}
\label{def:EAT2}
We define a function
\[
\phi: \mc{E}(S) \longrightarrow
\mc{A}(S)
\]
as follows.  If $N_+=1$, then by \eqref{eqn:EC1} there is a unique
element $E\in\mc{E}(S)$, given by $E_1=\{-1,\ldots,-N_-\}$.  In the tree
$\phi(E)$, the path from the positive leaf to the negative
leaf $-N_-$ goes through $N_--1$ trivalent joining vertices.  These joining
vertices are adjacent to the negative leaves
$-1,\ldots,-(N_--1)$ in that order.  See Figure~\ref{fig:unique}.

\begin{figure}
\begin{center}
\scalebox{0.6}{\includegraphics{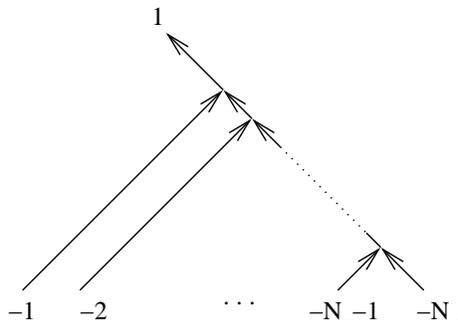}}
\end{center}
\caption{When $N_+=1$, there is a unique admissible tree.}
\label{fig:unique}
\end{figure}

If $N_+>1$, then given $(E_1,\ldots,E_{N_+})\in\mc{E}$, construct the
admissible tree $T=\phi(E_1,\ldots,E_{N_+})$ inductively as follows.  We
first define an oriented weighted tree $T_1$ with two positive leaves
and with negative leaves indexed by $E_1$.  To construct $T_1$, draw a
downward path $\gamma$ from the first positive leaf to the negative
leaf indexed by $\min(E_1)$.  This path will go through $|E_1'|$
trivalent joining vertices.  For each of these joining vertices, the
downward edge not on $\gamma$ is incident to a negative leaf.  These
negative leaves are those indexed by $E_1'$, in decreasing order.  Now add a
trivalent splitting vertex below the lowest joining vertex; the new
outgoing edge of this vertex is incident to the second positive leaf
of $T_1$.  The negative leaf of $T_1$ indexed by $j\in E_1$ has weight
$a_j$.  The first positive leaf has weight $a_1$ and the second
positive leaf has weight $\sum_{j\in E_1}a_j-a_1$.  The leaf weights
of $T_1$ extend to unique weights on the internal edges satisfying the
conservation condition at the vertices of $T_1$.  Condition
\eqref{eqn:EC2} insures that the edges of $T_1$ all have
positive weight.

Next, by \eqref{eqn:EC3} and induction, $\overline{E}$ determines an
oriented weighted tree $\overline{T} =
\phi(\overline{E})\in\mc{A}(\overline{S})$.
To construct the tree $T$, glue the second positive leaf of $T_1$ to
the negative leaf of $\overline{T}$ indexed by $\xi(\min(E_1))$.  (See
Figure~\ref{fig:induction}.)  It is straightforward to verify that the
tree $T$ is admissible.
\end{definition}

The above notions are connected to the definition of $f_\theta$ as
follows.

\begin{lemma}
\label{lem:EATSum}
\begin{equation}
\label{eqn:EF}
f_\theta(S) = \sum_{E\in\mc{E}(S)} W_\theta(\phi(E),P_{\phi(E)}).
\end{equation}
\end{lemma}

\begin{proof}
Summary: Unraveling the definitions shows that if one expands the
recursive formula \eqref{eqn:FTS} for $f_\theta(S)$, then one obtains
a sum indexed by elements $E=(E_1,\ldots,E_{N_+})\in\mc{E}(S)$, and the
summand corresponding to $E$ agrees with the summand on the right hand
side of \eqref{eqn:EF}.

Details:  If $N_+=1$ then \eqref{eqn:EF} is immediate
from the definitions since there is just one summand.  Suppose now
that $N_+>1$.  Write $E_1=\{j_1>\cdots>j_q\}$.  By equation
\eqref{eqn:FTS} and induction on $N_+$ we can write
\begin{equation}
\label{eqn:EATSum}
f_\theta(S) = \sum_{E_1} \sum_{\overline{E} \in \mc{E}(\overline{S})}
W_\theta\left(\overline{T},P_{\overline{T}}\right)
\prod_{n=1}^q \delta_\theta
\left(a_1- \sum_{k=1}^{n-1} a_{j_k} \, , \, a_{j_n}\right).
\end{equation}
Here the sum is over $E_1$ satisfying condition
\eqref{eqn:EC2}.  To clarify this, write $\overline{E}\fedqe
(\xi(E_2),\ldots,\xi(E_{N_+}))$ and $E\eqdef (E_1,E_2,\ldots,E_{N_+})$.  By
\eqref{eqn:EC3},  $\overline{E}\in\mc{E}(\overline{S})$ is equivalent to
$E\in\mc{E}(S)$.  Thus \eqref{eqn:EATSum} can be regarded as a sum
over $E\in\mc{E}(S)$.

Now the tree $T_1$ in the definition of $\phi(E)$ is admissible, and
with its canonical edge pairing has weight
\begin{equation}
\label{eqn:EATSummand1}
W_\theta(T_1,P_{T_1}) = m(e_0)
\prod_{n=1}^q \delta_\theta
\left(a_1- \sum_{k=1}^{n-1} a_{j_k} \, , \, a_{j_n}\right).
\end{equation}
Here $e_0$ denotes the unpaired edge in the canonical edge pairing
$P_{T_1}$.  Since $e$ is also the edge along which $T_1$ and
$\overline{T}$ are glued together, we  have
\begin{equation}
\label{eqn:EATSummand2}
W_\theta(\phi(E),P_{\phi(E)})
=
\frac{W_\theta(T_1,P_{T_1})}{m(e_0)}
W_\theta\left(\overline{T},P_{\overline{T}}\right).
\end{equation}
Regarding \eqref{eqn:EATSum} as a sum over $E\in\mc{E}(S)$ and
plugging in \eqref{eqn:EATSummand1} and \eqref{eqn:EATSummand2} gives
\eqref{eqn:EF}.
\end{proof}

Lemma~\ref{lem:IFT} now follows from Lemma~\ref{lem:EATSum} and:

\begin{lemma}
\label{lem:EATBijection}
The function $\phi:\mc{E}(S)\to\mc{A}(S)$ is a bijection.
\end{lemma}

\begin{proof}
Observe that if $E=(E_1,\ldots,E_{N_+})\in\mc{E}$, then $E_i$ is the
set of negative leaves of $\phi(E)$ that are accessible by downward
paths starting at the $i^{th}$ positive leaf.  Thus it is enough to
show the following: Given an admissible tree $T\in\mc{A}(S)$, let
$E_i$ denote the set of negative leaves of $T$ that are accessible by
downward paths starting at the $i^{th}$ positive leaf.  Then
$(E_1,\ldots,E_{N_+})\in\mc{E}$ and $T=\phi(E_1,\ldots,E_{N_+})$.

To begin the proof of this, let $\gamma$ denote the downward path from
the first positive leaf to the negative leaf indexed by $\min(E_1)$.
By the definition of $E_1$, for each $j\in E_1'$ there must be a
joining vertex $w_j$ on $\gamma$ such that there is a downward path
$\gamma_j$ from $w_j$ to the negative leaf indexed by $j$, where
$\gamma_j$ does not intersect $\gamma$ except at $w_j$. By
admissibility condition (c), the first edge on $\gamma_j$ cannot be
incident to a joining vertex (other than $w_j$).  Then by
admissibility condition (a), $w_j$ is the only joining vertex on the
path $\gamma_j$.  Hence the vertices $w_j$ are distinct.  By condition
(c) again, if $j>j'$ then $w_j$ is above $w_{j'}$ on $\gamma$.  If
$N_+=1$ then we are done, so assume henceforth that $N_+>1$.

By the definition of $E_1$, the path $\gamma$ cannot meet any joining
vertex except for the $w_j$'s.  We have seen above that $\gamma_j$
cannot meet any joining vertex except for $w_j$.  By admissibility
condition (d), the path $\gamma_j$ cannot meet any splitting vertex.
Hence the path $\gamma_j$ consists of a single edge incident to $w_j$
and the negative leaf indexed by $j$.  Since the tree $T$ is connected, the
path $\gamma$ must meet at least one splitting vertex.  By
admissibility condition (a), any such splitting vertex is below all of
the joining vertices on $\gamma$.  By admissibility condition (b),
there is only one splitting vertex on $\gamma$, call it $v_1$.

Let $e$ denote the outgoing edge of $v_1$ that is not on $\gamma$.
Cut $T$ along $e$ to obtain two oriented weighted trees $T_1$ and
$\overline{T}$, where $T_1$ contains $\gamma$ and $\overline{T}$ does
not.  Order the positive and negative leaves of $\overline{T}$ with
the orderings induced from those of $T$ (where the negative leaf of
$\overline{T}$ takes the place of $\min(E_1)$ in the ordering).  Then
$\overline{T}$ is admissible, so by induction,
\[
(\xi(E_2),\ldots,\xi(E_{N_+})) \in \mc{E}(a_2,\ldots,a_{N_+} \mid a_{-E_1})
\]
and $\overline{T}=\phi(\xi(E_2),\ldots,\xi(E_{N_+}))$.  Clearly
$E_1'$ is disjoint from $E_2,\ldots,E_{N_+}$, so
$(E_1,\ldots,E_{N_+})\in\mc{E}$.  (Condition \eqref{eqn:EC2} holds because
the two outgoing edges of $v_1$ have positive multiplicity.)  Now $T_1$
and $\overline{T}$ here are the same as in the definition of $\phi$,
so $T=\phi(E_1,\ldots,E_{N_+})$.
\end{proof}

\section{Counting zeroes of the linearized section}
\label{sec:count}

Continue with the assumptions from the first paragraph of
\S\ref{sec:section}.  We will now show that the relative Euler class
$\#\mathfrak{s}_0^{-1}(0)$ defined in \S\ref{sec:REC} agrees with the
combinatorial quantity $c_\theta(S)$ defined in \S\ref{sec:ctheta}, and
also that the latter does not depend on the choice of ordering in
Definition~\ref{def:ctheta} or \ref{def:change}.

To start, by Corollary~\ref{cor:trivialCase} we may assume that
$\overline{N}_+=N_+$ or $\overline{N}_-=N_-$.  Next, observe that
$\#\mathfrak{s}_0^{-1}(0)$ does not depend on the ordering of the
$a_i$'s for $i\in\{1,\ldots,\overline{N}_+\}$ or on the ordering of
the $a_i$'s for $i\in\{-1,\ldots,-\overline{N}_-\}$.  Let us then
choose these two orderings such that, as in
\eqref{eqn:reordering}, we have:
\begin{itemize}
\item
If $0 < i < j \le \overline{N}_+$, then $\ceil{a_i\theta}/a_i \le
\ceil{a_j\theta}/a_j$.
\item
If $0 > i > j \ge -\overline{N}_-$, then $\floor{a_i\theta}/a_i \ge
\floor{a_j\theta}/a_j$.
\end{itemize}
Then by Definitions~\ref{def:ctheta} and \ref{def:change}, to prove that
$\#\mathfrak{s}_0^{-1}(0)= c_\theta(S)$ and that the latter is
well-defined, it suffices to prove:
\begin{proposition}
\label{prop:count}
Under the above assumptions,
\begin{equation}
\label{eqn:count}
\# \mathfrak{s}_0^{-1}(0) = f_\theta(S).
\end{equation}
\end{proposition}

The rest of this section proves Proposition~\ref{prop:count}.  As
explained in \S\ref{sec:overview}, this will complete the proof of
Theorem~\ref{thm:main}, assuming the results from \cite{obg2} that are
quoted in \S\ref{sec:section}.

\subsection{Setup for the proof of Proposition~\ref{prop:count}}
\label{sec:countSetup}

To start, note that the statement of Proposition~\ref{prop:count} is
symmetric under switching positive ends with negative ends and
replacing $\theta$ by $-\theta$.  (The symmetry for $f_\theta(S)$
holds by Corollary~\ref{cor:fsym1}, while the symmetry for
$\#\mathfrak{s}_0^{-1}(0)$ is a straightforward consequence of the
definitions.)  By this symmetry, we can assume that
\[
\overline{N}_-=N_-.
\]

Now fix admissible $\gamma=\{\gamma_i\}$, $r$, and $R$ for use in the
definition of $\mathfrak{s}_0$.  (We will be more particular about these
choices later.)  By Corollary~\ref{cor:REC}, for the purposes of
computing $\#\mathfrak{s}_0^{-1}(0)$ we may assume that
\begin{equation}
\label{eqn:CPA2}
S(t)=\theta.
\end{equation}
Then the operator $D_\Sigma$ is complex linear, and $\Coker(D_\Sigma)$
is a complex vector space.  Let $\mc{O}_\C\to\mc{M}_R$ denote the
complex vector bundle whose fiber over $\Sigma$ is
$\Hom_{\C}(\Coker(D_\Sigma),\C)$.  There is a natural identification
of real vector bundles
\begin{equation}
\label{eqn:OOC}
\mc{O}=\mc{O}_\C.
\end{equation}
Under this identification, the section $\mathfrak{s}_0$ of $\mc{O}$
corresponds to a section $\mathfrak{s}_{\C}$ of $\mc{O}_{\C}$ defined by
\[
\mathfrak{s}_{\C}(\Sigma)(\sigma) \eqdef \mathfrak{s}_0(\Sigma)(\sigma) +
i\mathfrak{s}_0(\Sigma)(-i\sigma).
\]
Equivalently, $\mathfrak{s}_{\C}$ is defined as in
Definition~\ref{def:s0}, but with real inner products replaced by
complex inner products.
Under the identification \eqref{eqn:OOC}, the orientation of
$\mc{O}_\C$ as a complex vector bundle differs from the orientation of
$\mc{O}$ defined in \S\ref{sec:OBO} by
$(-1)^{\text{rank}_\C(\mc{O}_\C)}$.  Thus
\begin{equation}
\label{eqn:sLsC}
\#\mathfrak{s}_0^{-1}(0) = (-1)^{N}\#\mathfrak{s}_{\C}^{-1}(0).
\end{equation}

We will now compute $\#\mathfrak{s}_\C^{-1}(0)$.  To describe
$\mathfrak{s}_\C$ more concretely, choose an isomorphism of each
$\mc{B}_i$ with $\C$ as a complex vector space, and use these
isomorphisms to regard the $\gamma_i$'s as complex numbers.  If
$\Sigma\in\mc{M}_R$ and $\sigma\in\Coker(D_\Sigma)$, then for each $i$
in $\{1,\ldots,{N}_+\}$ or $\{-1,\ldots,-{N}_-\}$ labeling an end of
$\Sigma$, the projection of $\sigma_i(\pm(R+r),\cdot)$ to $\mc{B}_i$
(where $\pm$ denotes the sign of $i$) also corresponds to a complex
number, which we denote simply by $\sigma_i$.  In this notation,
\begin{equation}
\label{eqn:sC}
\mathfrak{s}_{\C}(\Sigma)(\sigma) =
\sum_{i=1}^{{N}_+}\overline{\gamma_i}\sigma_i -
\sum_{i=-1}^{{N}_-} \overline{\gamma_i}\sigma_i \in\C.
\end{equation}
Here we interpret $\gamma_{N_+}=0$ when $\overline{N}_+\neq N_+$.

\subsection{Outline of the argument}
\label{sec:CO}

The relative Euler class $\#\mathfrak{s}_\C^{-1}(0)$ is determined by the
restriction of $\mathfrak{s}_{\C}$ to the boundary of $\mc{M}_R$.  The
prototypical fact is that given a generic smooth function $f:D^2\to\C$
which does not vanish on $\partial D^2$, the algebraic count of points
$x\in D^2$ with $f(x)=0$ is equal to a count of points $x\in\partial
D^2$ with $f(x)>0$.  Roughly speaking, our strategy for computing
$\#\mathfrak{s}_\C^{-1}(0)$ is to understand the relevant boundary
behavior by induction on the dimension of boundary strata.

The precise procedure is as follows.  We will choose a large constant
$r_1>0$, and assume that the constant $r$ in \S\ref{sec:REC} is chosen
sufficiently large that $r > Nr_1$.  If $k$ is a positive integer,
define
\[
R_k \eqdef R - k r_1.
\]

\begin{definition}
For $k=0,\ldots,N-2$, define $\mc{M}^k$ to be the set
of $\Sigma\in\mc{M}_R$ such that:
\begin{itemize}
\item
The tree $\tau(\Sigma)$ has
trivalent vertices $v_1,\ldots,v_{k}$ with $\rho(v_i)\in\{\pm R_i\}$.
\item
All other vertices $w$ of $\tau(\Sigma)$ have
$\rho(w)\in[-R_{k+1},R_{k+1}]$.
\end{itemize}
Let $\partial\mc{M}^k$ denote the set of $\Sigma\in\mc{M}^k$ such that
there is at least one vertex $w$ with $\rho(w)\in\{\pm R_{k+1}\}$.
\end{definition}

Note that the interior of $\mc{M}^k$ is a smooth manifold of dimension
$2N-4-k$.  Also $\mc{M}^k$ is contained in the interior of
$\partial\mc{M}^{k-1}$, and this inductively determines an orientation
of the interior of $\mc{M}^k$.  We will later define a smaller space
$\mc{N}^k$ obtained by discarding certain components from $\mc{M}^k$.
This will satisfy $\mc{N}^0=\mc{M}^0$ and $\mc{N}^k
\subset \partial\mc{N}^{k-1}\eqdef\mc{N}^{k-1}\cap\partial\mc{M}^{k-1}$.
 Next, for each $k=1,\ldots,N-2$ and each component of $\mc{N}^{k-1}$,
we will pick a suitable $\epsilon_k\in
\{1,\ldots,N_+\}\cup\{-1,\ldots,-N_-\}$ labeling an end of the
$\Sigma$'s.  Then to each component of $\mc{N}^{k-1}$ we will have
associated $k$ ends $\epsilon_1,\ldots,\epsilon_k$ (since
$\mc{N}^{k-1}\subset\mc{N}^i$ for $i<k$), and these will be chosen to
be distinct.

\begin{definition}
Let $n_k$ denote the algebraic count of points $\Sigma\in\mc{N}^k$ such
that there exist constants $\Lambda_1,\ldots,\Lambda_k>0$ satisfying
the {\em $k$-boundary equation\/}
\begin{equation}
\label{eqn:boundaryk}
\forall\sigma\in\Coker(D_\Sigma):\quad \mathfrak{s}_{\C}(\Sigma)(\sigma) =
\Lambda_1\sigma_{\epsilon_1} +
\cdots + \Lambda_k\sigma_{\epsilon_k}.
\end{equation}
\end{definition}

To be more precise, $n_k$ is defined by making a small perturbation of
$\mathfrak{s}_\C$ to a section $\mathfrak{s}_\C'$ on $\mc{N}^k$ so that all
solutions to the $\mathfrak{s}_\C'$ analogue of
\eqref{eqn:boundaryk} on $\mc{N}^k$ with
$\Lambda_1,\ldots,\Lambda_k\ge 0$ have $\Lambda_1,\ldots,\Lambda_k>0$
and are cut out transversely, and then counting these solutions with
signs.  We will now specify the sign convention and then explain why
the count does not depend on the perturbation.

If $\epsilon_{k+1},\ldots,\epsilon_N$ are the remaining ends in any
order, then a complex basis for $\Coker(D_\Sigma)$ is given by
$(\sigma^{(1)},\ldots,\sigma^{(N-2)})$, where $\sigma^{(i)}$ denotes the
special cokernel element satisfying
\begin{equation}
\label{eqn:psi(i)}
\sigma^{(i)} \in V_{\epsilon_i,\epsilon_{N-1},\epsilon_N}, \quad \quad
\sigma^{(i)}_{\epsilon_i}=1.
\end{equation}
The $k$-boundary equation \eqref{eqn:boundaryk} for the perturbed
section $\mathfrak{s}_\C'$ is then equivalent to the open condition
\[
\mathfrak{s}_\C'(\Sigma)(\sigma^{(i)})\neq 0, \quad\quad i=1,\ldots,k
\]
together with the equations
\begin{align}
\label{eqn:simplified'}
\op{arg}\left(\mathfrak{s}_{\C}'(\Sigma)\left(\sigma^{(i)}\right)\right) = 0, &
 \quad
\quad i=1,\ldots,k,\\
\nonumber
\mathfrak{s}_{\C}'(\Sigma)(\sigma^{(i)}) = 0, & \quad
\quad i=k+1,\ldots,N-2.
\end{align}
Writing the equations in this order determines the sign convention for
$n_k$.

Now $n_0$ is well-defined and equal to the integer that we want to
compute, namely
\begin{equation}
\label{eqn:nsC}
n_0=\#\mathfrak{s}_{\C}^{-1}(0),
\end{equation}
because Lemma~\ref{lem:zab} guarantees that all zeroes of
$\mathfrak{s}_{\C}$ over $\mc{M}_R$ are in the interior of $\mc{N}^0$, and
the sign conventions for counting agree.  The following lemma provides
an inductive strategy for computing $n_0$.

\begin{lemma}
\label{lem:strategy}
Suppose that the following hold  for all $k=1,\ldots,N-1$:
\begin{description}
\item{(Ind1)} If $\Sigma\in\partial\mc{N}^{k-1}$
solves the $k$-boundary equation with
  $\Lambda_1,\ldots,\Lambda_{k-1}>0$ and $\Lambda_k\ge 0$, then
  $\Sigma\in\mc{N}^k$.
\item{(Ind2)} If $\Sigma\in \mc{N}^{k-1}$ solves the $k$-boundary
  equation\footnote{When $k=N-1$, statement (Ind1) is vacuously true
since $\partial\mc{N}^{N-2}=\emptyset$, while statement (Ind2) is to
be interpreted as saying that if $\Sigma\in\mc{N}^{N-2}$ satisfies the
$(N-2)$-boundary equation with $\Lambda_1,\ldots,\Lambda_{N-2}\ge 0$,
then $\Lambda_1,\ldots\Lambda_{N-2}>0$.} with
$\Lambda_1,\ldots,\Lambda_k\ge 0$, then
$\Lambda_1,\ldots,\Lambda_{k-1}>0$.
\end{description}
Then for all $k=1,\ldots,N-2$:
\begin{description}
\item{(a)}
$n_k$ is well-defined, independent of the small perturbation
$\mathfrak{s}_\C'$ of $\mathfrak{s}_\C$.
\item{(b)}
$n_k=(-1)^{k-1}n_{k-1}$.
\end{description}
\end{lemma}

\begin{proof}
First note that for all $k=1,\ldots,N-2$, by combining statement
(Ind1) for $k$ with statement (Ind2) for $k+1$, we have:
\begin{description}
\item{(Ind1$'$)}
If $\Sigma\in\partial\mc{N}^{k-1}$
solves the $k$-boundary equation with
  $\Lambda_1,\ldots,\Lambda_{k-1}>0$ and $\Lambda_k\ge 0$, then
  $\Sigma\in\mc{N}^k$ and $\Lambda_k>0$.
\end{description}

(a) To see that $n_k$ is well-defined, we need to show that if
$\Sigma\in\mc{N}^k$ solves the $k$-boundary equation with
$\Lambda_1,\ldots,\Lambda_k\ge 0$, then (i)
$\Lambda_1,\ldots,\Lambda_k>0$ and (ii) $\Sigma\notin\partial\mc{N}^k$.
Assertion (i) follows from statement (Ind2) for $k+1$.  Assertion (ii)
then follows from statement (Ind1$'$) for $k+1$.

(b) Consider the set
\[
Z \eqdef \{\Sigma\in\mc{N}^{k-1}\mid \mbox{$\Sigma$ solves the
$k$-boundary equation with $\Lambda_1,\ldots,\Lambda_k\ge 0$}\}.
\]
Conditions (Ind2) and (Ind1$'$) assert that:
\begin{itemize}
\item
Every $\Sigma\in Z$ has $\Lambda_1,\ldots,\Lambda_{k-1}>0$.
\item
If $\Sigma\in Z\cap\partial\mc{N}^{k-1}$, then $\Sigma\in\mc{N}^k$
(whence $\Sigma\in\op{int}(\partial\mc{N}^{k-1})$) and $\Lambda_k>0$.
\end{itemize}
It follows that we can choose the small perturbation $\mathfrak{s}_\C'$ of
the section $\mathfrak{s}_{\C}$ over $\mc{N}^{k-1}$ to arrange not only
that the points counted by $n_{k-1}$ and $n_k$ are cut out
transversely, but also that the $\mathfrak{s}_\C'$ version of $Z$, call it
$Z'$, is a one-manifold with boundary
\[
\partial Z' = (Z' \cap \mc{N}^k) \bigsqcup \left\{ \Sigma
\in Z' \mid \Lambda_k = 0 \right\}.
\]
Also $Z'$ is compact by Lemma~\ref{lem:mrc}.  Thus $Z'$ is a cobordism
between the set of solutions to the perturbed $k$-boundary equation on
$\mc{N}^k$ with $\Lambda_1,\ldots,\Lambda_k>0$, and the set of
solutions to the perturbed $(k-1)$-boundary equation on $\mc{N}^{k-1}$
with $\Lambda_1,\ldots,\Lambda_{k-1}>0$.  After an orientation check
it follows that $n_k=(-1)^{k-1}n_{k-1}$.
\end{proof}
  
We will see that if the $\epsilon_k$'s and $\mc{N}^k$'s and the
various constants are chosen carefully, then points (Ind1) and (Ind2)
hold for each $k$.  We will then be reduced to the problem of
computing $n_{N-2}$, i.e.\ counting solutions to the $(N-2)$-boundary
equation on $\mc{N}^{N-2}$.  Since every $\Sigma\in\mc{N}^{N-2}$ has a
trivalent tree $\tau(\Sigma)$, it follows that $\mc{N}^{N-2}$ is a
union of $(N-2)$-dimensional tori (cf.\ Lemma~\ref{lem:covering}), and
counting the solutions to the equations \eqref{eqn:simplified'} will
reduce to a determinant calculation.

\subsection{Decay estimates}
\label{sec:sCDecay}

To work with the $k$-boundary equation \eqref{eqn:boundaryk}, we need
a preliminary discussion of the relative sizes of the different
contributions to $\mathfrak{s}_{\C}(C)(\sigma)$ in equation \eqref{eqn:sC}.
By Corollary~\ref{cor:REC}, we can choose the $\gamma_i$'s
such that
\begin{equation}
\label{eqn:CPA4}
\begin{split}
|\gamma_i| > r_2|\gamma_{i+1}|, & \quad \quad i>0,\\
|\gamma_{j}| >
r_2|\gamma_{j-1}|, & \quad \quad j<0,
\end{split}
\end{equation}
where $r_2>1$ is a large constant.  We can also assume that
\begin{equation}
\label{eqn:CPA3}
i\in\{1,\ldots,\overline{N}_+\} \cup \{-1,\ldots,-N_-\}
\Longrightarrow \gamma_i \notin \R.
\end{equation}

Let $i$, $j$, and $k$ be distinct ends.  Recall from \eqref{eqn:sC}
that for a special cokernel element $\sigma\in V_{i,j,k}$, we have
$\sigma_l=0$ for $l\notin\{i,j,k\}$, so only three terms contribute to
$\mathfrak{s}_{\C}(C)(\sigma)$. Often one term dominates the other two, in
the following sense.  By
\eqref{eqn:sLsC} and
\eqref{eqn:CPA3}, if $K$ is sufficiently large, and if nonzero
$\sigma\in V_{i,j,k}$ satisfy
$
|\overline{\gamma_i}\sigma_i| > K|\overline{\gamma_j}\sigma_j|,\,
K|\overline{\gamma_k}\sigma_k|$,
then
\begin{gather}
\label{eqn:d1}
\mathfrak{s}_{\C}(C)(\sigma)\neq 0,\\
\label{eqn:d2}
\arg(\mathfrak{s}_{\C}(C)(\sigma)) \neq \arg(\sigma_i).
\end{gather}

\begin{definition}
\label{def:dominate}
Write
\[
i \searrow j,k
\]
if nonzero $\sigma\in V_{i,j,k}$ satisfy \eqref{eqn:d1} and
\eqref{eqn:d2}. 
\end{definition}

The following lemma will be used repeatedly.

\begin{lemma}
\label{lem:CLS}
If $r$ is sufficiently large, if $r_2$ is sufficiently large with
respect to $r_1$, and if $R$ is sufficiently large with respect to all
other choices, then the following holds: Let $i$ be a positive end and
let $j>j'$ be negative ends.  Let $v$ denote the central vertex for
$i$, $j$, and $j'$.  Then:
\begin{description}
\item{(a)}
If the path $P_{v,i}$ stays above the level $\rho=R_N$ and if
$i\le\overline{N}_+$, then $i\searrow j,j'$.
\item{(b)}
If the path $P_{v,j}$ stays below the level $\rho=-R_N$ (e.g.\ if the
path $P_{j,j'}$ does), then $j\searrow i, j'$.
\end{description}
\end{lemma}

\begin{proof}
We begin with a key estimate.  Let $\sigma\in V_{i,j,j'}$ be
normalized so that $|\sigma(z)|=1$ for some $z\in\Sigma$ with
$p(z)=v$.  Let $x_i\in\tau(\Sigma)$ denote the point on the edge
corresponding to the $i^{th}$ end for which $\rho(x_i)=R+r$.  Let
$x_j\in\tau(\Sigma)$ denote the point on the edge corresponding to the
$j^{th}$ end for which $\rho(x_j)=-(R+r)$, and define $x_{j'}$ likewise.
By Lemma~\ref{lem:specialWinding}, Corollary~\ref{cor:rsc}, and
Proposition~\ref{prop:approx1}, if $r$ is sufficiently large then
there is a constant $r'$ such that
\begin{equation}
\label{eqn:CLS}
\left|\log |\sigma_i| + \sum_{e\in
P_{v,x_i}^+}\ell(e)\left(\frac{\ceil{m(e)\theta}}{m(e)}-\theta\right) +
\sum_{e\in P_{v,x_i}^-} \ell(e) \left(\theta -
\frac{\floor{m(e)\theta}}{m(e)}\right)\right| \le r'.
\end{equation}
The estimate \eqref{eqn:CLS} also holds if $i$ is replaced by $j$ or
 $j'$.

(a) If the path $P_{v,i}$ stays above the level $\rho=R_N$, then the
estimate \eqref{eqn:CLS} implies that there is a constant $\kappa>0$
such that
\[
\log|\sigma_i| \ge -r\left(\frac{\ceil{a_i\theta}}{a_i} - \theta\right)      -
\kappa r_1 - r'.
\]
The analogue of \eqref{eqn:CLS} for $j$ and $j'$ implies that the
constant $\kappa$ can be chosen so that
\[
\log|\sigma_j|, \log|\sigma_{j'}| \le -\kappa R + r'.
\]
By the above two inequalities, assertion (a) holds provided that $R$
is sufficiently large with respect to all of the other choices.

(b) If the path $P_{v,j}$ stays below the level $\rho=-R_N$, then
\eqref{eqn:CLS} and its analogues for $j$ and $j'$ imply that there is
a constant $\kappa>0$ such that
\[
\begin{split}
\log|\sigma_i| & \le -\kappa R + r',\\
\log|\sigma_j| & \ge
-r\left(\theta - \frac{\floor{a_j\theta}}{a_j}\right) -\kappa r_1 -
 r',\\
\log|\sigma_{j'}| &\le -r\left(\theta -
 \frac{\floor{a_{j'}\theta}}{a_{j'}}\right) +r'.
\end{split}
\]
Recall that our ordering convention gives
$\floor{a_j\theta}/a_j\ge\floor{a_{j'}\theta}/a_{j'}$.  So by
\eqref{eqn:CPA4}, assertion (b) holds provided that $r_2$ is large
enough with respect to $r_1$ and $r'$, and $R$ is large enough with
respect to all other choices.
\end{proof}

The obvious symmetric analogue of Lemma~\ref{lem:CLS} with positive
and negative ends switched also holds.  Henceforth assume that the
constants are chosen so that the conclusions of Lemma~\ref{lem:CLS}
and its symmetric analogue hold.  (We will later need to choose $r_1$
large.)

\subsection{Processing the positive ends}
\label{sec:positive}

To begin the inductive process, we now define $\epsilon_k$ and $\mc{N}^k$
when $k<N_+$ and verify that the crucial properties (Ind1) and (Ind2)
hold in this case.

When $k<N_+$, we choose $\epsilon_k$ to be the positive end labeled by $k$.

In the definition of $\mc{N}^k$ and below, we will use the following
notation.  If $v$ is a vertex of a tree and $e$ is an edge incident to
$v$, let $A(v,e)$ denote the set of ends that are accessible via
paths starting from $v$ along the edge $e$.  Also, if there is a
unique downward path from $v$ to a negative leaf, then we denote the
corresponding negative leaf by $v^-$.  We generally refer to the
leaves of a tree $\tau(\Sigma)$ as ``ends'', and identify the ends
with their labels in $\{1,\ldots,N_+\}\cup\{-1,\ldots,-N_-\}$.

\begin{definition}
\label{def:Nk+}
For $k<N_+$, define $\mc{N}^k$ to be the set of $\Sigma\in\mc{M}^k$ that
satisfy the following conditions for all $i=1,\ldots,k$:
\begin{description}
\item{(a)}
$v_i$ is a trivalent splitting vertex with $\rho(v_i)=-R_i$.
\item{(b)}
$v_i$ has an outgoing edge $e_i^0$ such that $\{i+1,\ldots,N_+\}\subset
A(v_i,e_i^0)$.
\item{(c)}
If $j$ is a negative end and $j\notin A(v_i,e_i^0)$, then $j\ge v_i^-$.
\end{description}
Let $e_i^+$ denote the outgoing edge of $v_i$ other than
$e_i^0$, and let $e_i^-$ denote the incoming edge of $v_i$.
\end{definition}

To better the above definition, we now consider the following
additional structure associated to elements of $\mc{N}^k$.

\begin{definition}
Let $\Sigma\in\mc{N}^k$ with $k<N_+$.  For each $i=0,\ldots,k$
define a tree $\tau_i^+$ and a forest $\tau_i^-$ inductively
as follows.
\begin{itemize}
\item
$\tau_0^+=\tau(\Sigma)$.
\item
For $i=1,\ldots,k$, the tree $\tau_i^+$ is obtained from the tree
$\tau_{i-1}^+$ by cutting along the edge $e_i^0$ and keeping the half
that contains the positive ends $i+1,\ldots,N_+$.
\item
$\tau_i^-$ is the complement of $\tau_i^+$ in $\tau(\Sigma)$.
\end{itemize}
\end{definition}

\begin{lemma}
\label{lem:tau-}
Let $\Sigma\in\mc{N}^k$ with $k<N_+$.  Then for each $i=1,\ldots,k$:
\begin{description}
\item{(a)}
The forest $\tau_i^-$ contains the vertices $v_1,\ldots,v_i$, and no
other splitting vertices.
\item{(b)} The positive leaves of $\tau_i^-$ are the first $i$
positive ends, together with, for each component of $\tau_i^-$, a
positive leaf where the component of $\tau_i^-$ is attached to $\tau_i^+$.
\item{(c)}
An upward path in $\tau(\Sigma)$ from $e_i^+$ to a positive end must
terminate at the $i^{th}$ end.
\item{(d)}
If there is a downward path from the end $i$ to the end $j$, then
$j\ge v_i^-$.
\end{description}
\end{lemma}

\begin{proof}
Let $s$, $p$, and $m$ denote the numbers of splitting vertices,
positive leaves, and components respectively in $\tau_i^-$.  By
construction, $\tau_i^-$ contains the splitting vertices
$v_1,\ldots,v_i$, so $s\ge i$.  Also, the positive leaves of
$\tau_i^-$ consist of one positive leaf in each component where it
attaches to $\tau_i^+$, together with some subset of the first $i$
positive ends.  Thus $p\le m+i$.  But since $\tau_i^-$ contains no
loops, we have $p\ge m+s$.  Therefore $s=i$ and $p=i+m$, and these
facts prove parts (a) and (b) respectively of the lemma.

To prove (c), note that the path under consideration stays in
$\tau_i^-\setminus\tau_{i-1}^-$.  We are then done by part (b).

To prove part (d), note that the downward path from $i$ to $j$
intersects the upward path from $v_i$ to $i$.  By part (c), the latter
path does not contain $e_i^0$.  It then follows that $j\notin
A(v_i,e_i^0)$, so we are done by Definition~\ref{def:Nk+}(c).
\end{proof}

\begin{lemma}
\label{lem:Ind+}
If $k<N_+$, then statements (Ind1) and (Ind2) hold.
\end{lemma}

\begin{proof}
(Ind1) Suppose $\Sigma\in\partial\mc{N}^{k-1}$ satisfies the $k$-boundary
equation \eqref{eqn:boundaryk} with $\Lambda_1,\ldots,\Lambda_k>0$.
We need to show that $\Sigma\in\mc{N}^k$.  Since we already know that
$\Sigma\in\mc{N}^{k-1}$, we just need to verify conditions (a)--(c) in
Definition~\ref{def:Nk+} for $i=k$, and we also need to check that
$\Sigma\in\mc{M}^k$.  We proceed in four steps.

{\em Step 1.\/} We first show that every vertex $v\in\tau(\Sigma)$ has
$\rho(v) < R_N$.

Let $v$ be a vertex with $\rho(v)\ge R_N$.  We can assume that
$\rho(v)$ is maximal.  Suppose first that $v$ is a splitting vertex.
Then $v$ has (at least) two outgoing edges incident to positive ends
$i < i'$.  Since $\Sigma\in\mc{N}^{k-1}$, Lemma~\ref{lem:tau-} implies that
$i,i'\ge k$.  We can also find a downward path from $v$ to a negative
end $j$, so that $v$ is the central vertex for $i$, $i'$, and $j$. By
the symmetric analogue of Lemma~\ref{lem:CLS}(b), we have $i\searrow
i', j$.  Let $0\neq\sigma\in V_{i,i',j}$.  If $i>k$, then the
$k$-boundary equation asserts that $\mathfrak{s}_{\C}(\Sigma)(\sigma)=0$,
contradicting \eqref{eqn:d1}.  If $i=k$ then the $k$-boundary equation
gives $\mathfrak{s}_{\C}(\Sigma)(\sigma)=\Lambda_i\sigma_i$, contradicting
\eqref{eqn:d2}.

The remaining possibility is that $v$ is a joining vertex with one
outgoing edge incident to a positive end $i$.  Then similarly to the
proof of Lemma~\ref{lem:cpm}, the partition minimality assumption
implies that $\overline{N}_+=N_+$.  We can find downward paths from
$v$, starting along distinct edges, to negative ends $j$ and $j'$.
Then $v$ is the central vertex for $i$, $j$, and $j'$.
Lemma~\ref{lem:CLS}(a) implies that $i\searrow j,j'$.  Let
$0\neq\sigma\in V_{i,j,j'}$.  If $i>k$, then the $k$-boundary equation
gives $\mathfrak{s}_{\C}(\Sigma)(\sigma)=0$, contradicting \eqref{eqn:d1}.  If
$i\le k$, then the $k$-boundary equation gives
$\mathfrak{s}_{\C}(\Sigma)(\sigma)=\Lambda_i\sigma_i$, which contradicts
\eqref{eqn:d2}.

{\em Step 2.\/} We now show that any vertex $v$ with $\rho(v)=-R_k$ is
a trivalent splitting vertex.

Suppose first that $v$ has (at least) two incoming edges, and let
$j_1>j_2$ be negative ends reached by downward paths starting along
these two edges.  Lemma~\ref{lem:CLS}(b) then gives $j_1 \searrow
j_2,N_+$.  But if $0\neq\sigma\in V_{j_1,j_2,N_+}$, then the
$k$-boundary equation implies that $\mathfrak{s}_{\C}(\Sigma)(\sigma)=0$,
which is a contradiction.

So $v$ has only one incoming edge.  In particular $v$ is a splitting
vertex, so $v$ cannot be in the forest $\tau_{k-1}^-$ by
Lemma~\ref{lem:tau-}(a).  Thus any upward path starting at $v$ stays
in $\tau_{k-1}^+$, and hence by Lemma~\ref{lem:tau-}(b) terminates at
a positive end indexed by $k,\ldots,N_+$.  If $v$ has more than two
outgoing edges, then at least two of these outgoing edges lead to
positive ends $i_1,i_2 > k$.  If $\sigma\in V_{i_1,i_2,v^-}$, then the
$k$-boundary equation gives $\mathfrak{s}_{\C}(\Sigma)(\sigma)=0$, while the
symmetric analogue of Lemma~\ref{lem:CLS}(a) gives $v^- \searrow i_1,
i_2$, which is a contradiction.

{\em Step 3.\/} Let $v_k$ be a trivalent splitting vertex with
$\rho(v_k)=-R_k$.  We now show that $v_k$ is unique and satisfies
conditions (b) and (c) in Definition~\ref{def:Nk+}.

To prove (b), let $e_k^0$ and $e_k^+$ denote the outgoing edges of
$v_k$.  The sets $A(v_k,e_k^0)$ and $A(v_k,e_k^+)$ cannot both contain
positive ends that are greater than $k$, or else we obtain a
contradiction as in Step 2.  So without loss of generality,
$A(v_k,e_k^+)$ does not contain any positive ends indexed by $i > k$.
Since $\Sigma\in\mc{N}^{k-1}$, the incoming edge $e_k^-$ of $v_k$
either comes out of the forest $\tau_{k-1}^-$ or is incident to a
negative end.  Hence $A(v_k,e_k^-)$ does not contain any positive ends
indexed by $i>k$.  Therefore all of the positive ends indexed by $i>k$
must be contained in $A(v_k,e_k^0)$.

To prove condition (c) in Definition~\ref{def:Nk+}, suppose
$A(v_k,e_k^+) \cup A(v_k,e_k^-)$ contains a negative end indexed by $j$
with $v_k^->j$.  Since $N_+\in A(v_k,e_k^0)$, the central vertex for
$j$, $v_k^-$, and $N_+$ is on the downward path from $v_k$ to $v_k^-$.
Lemma~\ref{lem:CLS}(b) then gives $v_k^- \searrow j, N_+$.  But if
$0\neq \sigma\in V_{j,v_k^-,N_+}$ then the $k$-boundary equation gives
$\mathfrak{s}_{\C}(\Sigma)(\sigma)=0$, which is a contradiction.

To prove that $v_k$ is unique, note that Step 2 and
Lemma~\ref{lem:tau-}(a) imply that $v_k\in\tau_{k-1}^+$, so by
Lemma~\ref{lem:tau-}(b) there is a unique upward path $P$ starting
along $e_k^+$, and the path $P$ leads to the $k^{th}$ positive end.
Now suppose that $w$ is another trivalent splitting vertex with
$\rho(w)=-R_k$.  Then $w$ must also have an outgoing edge $e$ such
that $A(w,e)$ contains no positive ends indexed by $i>k$, there is a
unique upward path $P'$ starting along $e$, and the path $P'$ leads to
the $k^{th}$ positive end.  The two upward paths $P$ and $P'$ must
intersect.  By proceeding from $w$ along $P'$ to its intersection with
$P$, and then backwarrds along $P$ to $v_k$, we find that $N_+\in
A(w,e)$, which is a contradiction.

{\em Step 4.\/} To complete the proof that $\Sigma\in\mc{N}^k$, we must
check that $\Sigma\in\mc{M}^k$, i.e.\ that any vertex $w$ other than
$v_1,\ldots,v_k$ satisfies $\rho(w)\in[-R_{k+1},R_{k+1}]$.

We know from Steps 1--3 that $\rho(w)\in(-R_k,R_N)$.  Suppose to get a
contradiction that $\rho(w)\in(-R_k,-R_{k+1})$.  We can assume that
$\rho(w)$ is minimal.  If $w$ has more than one incoming edge, then we
get a contradiction as in Step 2.  So $w$ has (at least) two outgoing
edges.  By Lemma~\ref{lem:tau-}(a) we know that $w$ is in
$\tau_k^+$ (which is well-defined by Steps 1--3), so upward paths starting
along these outgoing edges lead to positive ends with labels in
$k+1,\ldots,N_+$.  This again gives a contradiction as in Step 2.

(Ind2) Suppose that $\Sigma\in\mc{N}^{k-1}$ solves the $k$-boundary
equation with $\Lambda_1,\ldots,\Lambda_k\ge 0$.  We need to show that
$\Lambda_1,\ldots,\Lambda_{k-1} \neq 0$.  Given
$i\in\{1,\ldots,k-1\}$, let $0\neq \sigma\in V_{i,N_+,v_i^-}$.
Observe that $v_i$ is the central vertex for $i$, $N_+$, and $v_i^-$.
Hence the symmetric analogue of Lemma~\ref{lem:CLS}(a) implies that
$v_i^-
\searrow i, N_+$, so $\mathfrak{s}_{\C}(\Sigma)(\sigma)\neq 0$.  However the
$k$-boundary equation gives $\mathfrak{s}_{\C}(\Sigma)(\sigma) =
\Lambda_i\sigma_i$, whence $\Lambda_i\neq 0$.
\end{proof}

\subsection{Processing the negative ends}

We now define $\epsilon_k$ and $\mc{N}^k$ for $k=N_+,\ldots,N-2$.
Here $\epsilon_k$ depends on the component of $\mc{N}^{N_+-1}$.  We
will then prove that conditions (Ind1) and (Ind2) continue to hold.

Given $\Sigma\in\mc{N}^{N_+-1}$, for each $i=1,\ldots,N_+$ let $E_i$
denote the set of negative ends that are accessible by downward paths
in the tree $\tau(\Sigma)$ starting from the $i^{th}$ positive end.
By Lemma~\ref{lem:tau-}(d), if $i<N_+$ then $v_i^-$ is the smallest
element of the set $E_i$.  For $i=1,\ldots,N_+$ define $E_i'\eqdef
E_i\setminus\{v_i^-\}$, where we interpret $v_{N_+}^-$ to be the
smallest element of the set $E_{N_+}$.  

\begin{lemma}
\label{lem:notethat}
For each $\Sigma\in\mc{N}^{N_+-1}$, the following hold:
\begin{description}
\item{(a)}
$
\{-1,\ldots,-N_-\} = E_{N_+} \sqcup \bigsqcup_{i=1}^{N_+-1}
E_i'.
$
\item{(b)}
$v_{N_+}^-=-N_-$.
\item{(c)}
$(E_1,\ldots,E_{N_+})$ depends only on the component of
$\mc{N}^{N_+-1}$ containing $\Sigma$.
\end{description}
\end{lemma}

\begin{proof}
(a) Let $j$ be a negative end; we need to show that there is a unique
positive end $i$ such that
\begin{description}
\item{(i)}
 $j\in E_i'$ if $i<N_+$, and $j\in E_{N_+}$
if $i=N_+$.
\end{description}
Note that condition (i) is equivalent to
\begin{description}
\item{(ii)}
The path $P$ from $j$ to $i$ is an upward path, and
\begin{description}
\item{(*)}
for all $k\in\{1,\ldots,N_+-1\}$, if the path $P$ meets the vertex
$v_k$, then the path $P$ continues along the edge $e_k^0$ in
Definition~\ref{def:Nk+}.
\end{description}
\end{description}
The reason is that by Lemma~\ref{lem:tau-}(c), condition (ii) fails if
and only if $j\notin E_i$ or there exists $k\in\{1,\ldots,N_+-1\}$
such that $j=v_k^-$ and $i=k$.  But there is a unique upward path $P$
starting at $j$ and satisfying condition (*), because every vertex
other than $v_1,\ldots,v_{N_+-1}$ is a joining vertex.

(b) We need to show that $-N_-\in E_{N_+}$.  If not, then part (a)
implies that $-N_-\in E_i'$ for some $i\in\{1,\ldots,N_+-1\}$.  But
then $-N_->v_i^-$, which is impossible.

(c) The sets $E_i$ can be characterized in terms of which ends are
accessible from which edges incident to the vertices
$v_1,\ldots,v_{N_+-1}$.  The latter information depends only on the
component of $\mc{N}^{N_+-1}$.
\end{proof}

\begin{definition}
For a given component of $\mc{N}^{N_+-1}$, define the sequence
$\epsilon_{N_+},\ldots,\epsilon_{N-2}$ by first listing the ends in
$E_1'$ in decreasing order, then listing the ends in $E_2'$ in
decreasing order, and so on up to $E_{N_+}'$.
\end{definition}

Thus the two remaining ends are the positive end $N_+$ and the
negative end $-N_-$; we denote these by $\epsilon_{N-1}$ and
$\epsilon_N$ respectively.

\begin{definition}
\label{def:Nk-}
For $k=N_+,\ldots,N-2$, define $\mc{N}^k$ to be the set of
$\Sigma\in\mc{N}^{k-1}\cap \mc{M}^k$ such that:
\begin{description}
\item{(a)}
$v_k$ is a trivalent joining vertex.
\item{(b)}
For one of the incoming edges of $v_k$, call it $e_k^-$, there is a
unique downward path starting along $e_k^-$, and this leads to the
negative end $\epsilon_k$.
\item{(c)}
$\rho(v_k)=+R_k$.
\end{description}
Let $e_k^0$ denote the incoming edge of $v_k$ other than $e_k^-$, and
let $e_k^+$ denote the outgoing edge of $v_k$.
\end{definition}

\begin{lemma}
\label{lem:Ind-}
For $k=N_+,\ldots,N-1$, statements (Ind1) and (Ind2) hold.
\end{lemma}

\begin{proof}
(Ind1) Recall that this is vacuous when $k=N-1$.  Now given
$k\in\{N_+,\ldots,N-2\}$, suppose $\Sigma\in\partial\mc{N}^{k-1}$
satisfies the $k$-boundary equation
\eqref{eqn:boundaryk} with $\Lambda_1,\ldots,\Lambda_{k-1}>0$ and
$\Lambda_k\ge 0$.  We need to show that $\Sigma\in\mc{N}^k$.  We proceed in
three steps.

{\em Step 1.\/} Let $v_k$ be a vertex with $\rho(v_k)=\pm R_k$.  We
now show that $v_k$ is unique and satisfies conditions (a) and (b) in
Definition~\ref{def:Nk-}.

To start, the tree $\tau(\Sigma)$ contains at most $N_+-1$ splitting
vertices, and these are accounted for by $v_1,\ldots,v_{N_+-1}$.
Since $k\ge N_+$, it follows that $v_k$ is a joining vertex with only
one outgoing edge.  Since $v_k$ is above all of the splitting
vertices, there is a unique upward path starting from $v_k$.  Since
$\Sigma\in\mc{N}^{k-1}$, a downward path starting at $v_k$ cannot lead
to an end in the set $\{\epsilon_1,\ldots,\epsilon_{k-1}\}$.

Now suppose that (a) or (b) fails.  Then there are downward paths
starting from $v_k$ along distinct incoming edges, leading to negative
ends $j > j'$ not in the set $\{\epsilon_1, \ldots,
\epsilon_k\}$.

To get a contradiction, suppose first that $\rho(v_k)=-R_k$.  In this
case Lemma~\ref{lem:CLS}(b) gives $j \searrow j', N_+$.  On the other
hand, if $0\neq\sigma\in V_{j,j',N_+}$ then the $k$-boundary equation
gives $\mathfrak{s}_{\C}(\Sigma)(\sigma)=0$, which contradicts
\eqref{eqn:d1}.

Suppose next that $\rho(v_k)=+R_k$.  Let $i$ denote the positive end
reached by the unique upward path starting from $v_k$.  Let
$0\neq\sigma\in V_{i,j,j'}$.  Note that $v_k$ is the central vertex for
$i$, $j$, and $j'$.  If $i<N_+$, then the $k$-boundary equation gives
$\mathfrak{s}_{\C}(\Sigma)(\sigma)=\Lambda_i\sigma_i$, while Lemma~\ref{lem:CLS}(a)
gives $i\searrow j, j'$.  This contradicts \eqref{eqn:d2}.  If
$i=N_+$, then the $k$-boundary equation gives
$\mathfrak{s}_{\C}(\Sigma)(\sigma)=0$.  However the partition minimality
assumption guarantees that $\overline{N}_+=N_+$ here, so
Lemma~\ref{lem:CLS}(a) applies again to give $i\searrow j,j'$, which
contradicts \eqref{eqn:d1}.  This completes the proof of (a) and (b).

To prove uniqueness of $v_k$, recall that there is a unique upward
path from $v_k$, and there is a unique downward path starting along
the incoming edge $e_k^-$ of $v_k$.  These paths lead respectively to
the positive end $i$ for which $\epsilon_k\in E_i'$, and to the
negative end $\epsilon_k$.  If $w$ is another vertex with these
properties, then the downward paths meet at some vertex other than
$v_k$ or $w$ (by uniqueness of these paths, since $v_k$ and $w$ are
joining vertices).  Then the downward paths and the upward paths
together contain a loop in $\tau(\Sigma)$, which is a contradiction.

{\em Step 2.\/} We now show that $v_k$ satisfies condition (c) in
Definition~\ref{def:Nk-}.  Suppose to the contrary that
$\rho(v_k)=-R_k$.  Choose a downward path from $v_k$ starting along
the incoming edge $e_k^0$; this leads to a negative end $j$
with $\epsilon_k>j$.  Then Lemma~\ref{lem:CLS}(b) gives
$\epsilon_k\searrow j, N_+$.  But if $0\neq\sigma\in
V_{N_+,\epsilon_k,j}$ then the $k$-boundary equation gives
$\mathfrak{s}_{\C}(\Sigma)(\sigma) = \Lambda_{k}\sigma_{\epsilon_k}$.  If
$\Lambda_k=0$ then this contradicts
\eqref{eqn:d1}, while if $\Lambda_k>0$ then this contradicts \eqref{eqn:d2}.

{\em Step 3.\/} To complete the proof that $\Sigma\in\mc{N}^k$, we must now
show that any vertex $w$ other than $v_1,\ldots,v_k$ has
$\rho(w)\in(-R_{k+1},R_{k+1})$.
The proof of this is essentially the same as the proof that $v_k$ is unique.

(Ind2) Given $k\in\{N_+,\ldots,N-1\}$, suppose that
$\Sigma\in\mc{N}^{k-1}$ solves the $k$-boundary equation with
$\Lambda_1,\ldots,\Lambda_k\ge 0$.  (When $k=N-1$, the hypothesis is
that $\Sigma\in\mc{N}^{N-2}$ satisfies the $(N-2)$-boundary equation
with $\Lambda_1,\ldots,\Lambda_{N-2}\ge 0$.)  Let
$j\in\{1,\ldots,k-1\}$; we must show that $\Lambda_j \neq 0$.  There are
two cases.

Case 1: $\epsilon_j\in E_{N_+}'$.  Since $\Sigma\in\mc{N}^j$,
Lemma~\ref{lem:CLS}(a) implies that $\epsilon_{N-1} \searrow
\epsilon_j,\epsilon_N$, so if $0\neq\sigma\in
V_{\epsilon_j,\epsilon_{N-1},\epsilon_N}$
then $\mathfrak{s}_{\C}(\Sigma)(\sigma)\neq 0$.  But the $k$-boundary equation
gives $\mathfrak{s}_{\C}(\Sigma)(\sigma)=\Lambda_j\sigma_{\epsilon_j}$ whence
$\Lambda_j\neq 0$.

Case 2: $\epsilon_j\notin E_{N_+}'$.  Then there exists
$i\in\{1,\ldots,N_+-1\}$ such that either $j=i$ or $\epsilon_j\in
E_i'$. Observe that
$v_i$ is the central vertex for $j$, $v_i^-$, and $N_+$.  The path
from $v_i$ to $v_i^-$ stays below the level $\rho=-R_N$, while the
paths from $v_i$ to $j$ and $N_+$ go above the level $\rho=+R_N$.  It
then follows from the decay estimate \eqref{eqn:CLS} that $v_i^-
\searrow j, N_+$.  Let $0\neq \sigma\in V_{j,v_i^-,N_+}$.
There are now two subcases.

Case 2a: $v_i^-\notin\{\epsilon_1,\ldots,\epsilon_k\}$. Then the
$k$-boundary equation asserts that
$\mathfrak{s}_{\C}(\Sigma)(\sigma)=\Lambda_j\sigma_{\epsilon_j}$, which together
with \eqref{eqn:d1} implies that $\Lambda_j\neq 0$.

Case 2b: $v_i^-=\epsilon_l$ with $l\in\{1,\ldots,k\}$.  If
$\Lambda_j=0$, then the $k$-boundary equation gives
$\mathfrak{s}_{\C}(\Sigma)(\sigma)=\Lambda_l\sigma_{\epsilon_l}$.  If
$\Lambda_l=0$ then this contradicts \eqref{eqn:d1}, while if
$\Lambda_l>0$ then this contradicts \eqref{eqn:d2}.
\end{proof}

By Lemmas~\ref{lem:Ind+} and \ref{lem:Ind-} and equations
\eqref{eqn:sLsC} and
\eqref{eqn:nsC}, we can apply Lemma~\ref{lem:strategy} inductively to
obtain
\begin{equation}
\label{eqn:nn2}
\#\mathfrak{s}_0^{-1}(0) = (-1)^{N+\frac{(N-2)(N-3)}{2}} \cdot n_{N-2}.
\end{equation}

\subsection{Rotation rates}
\label{sec:rotrat}

To prepare to compute $n_{N-2}$, we now digress to consider the
following question:  Let $\Sigma\in\mc{N}^{N-2}$, let $\sigma\in
V_{i,j,k}$ be a nonzero special cokernel element, and let $v$ be a
vertex of $\tau(\Sigma)$.   Approximately how does the argument of
$\sigma_i/\sigma_j$ change as we rotate the corresponding branch point
in the $S^1$ direction?

\begin{definition}
\label{def:rr}
If $\Sigma\in\mc{M}$, if $v$ is a vertex of $\tau(\Sigma)$, and if
 $\sigma\in V_{i,j,k}$ is nonzero, define the {\em rotation rate\/}
$
r(\sigma_i/\sigma_j,v)\in\Q
$
 as follows.
Let $e_i$ and $e_j$ denote the edges of $v$ that lead to the ends $i$
and $j$ respectively.  Then
\begin{equation}
\label{eqn:rr}
r(\sigma_i/\sigma_j,v) \eqdef \frac{\eta(\sigma,e_j)}{m(e_j)} -
\frac{\eta(\sigma,e_i)}{m(e_i)}.
\end{equation}
\end{definition}

Recall that $\eta(\sigma,e)$ denotes the winding number of $\sigma$ around
$e$, which is computed by Lemma~\ref{lem:specialWinding}.  Note that
if $v$ is not on the path $P_{i,j}$, then $r(\sigma_i/\sigma_j,v) =
0$. 

 The following lemma is a special case of Proposition~\ref{prop:RBP0},
 which is proved in \S\ref{sec:RR}.

\begin{lemma}
\label{lem:RBP}
For all $\varepsilon>0$, if the constant $r_1$ in \S\ref{sec:CO} is
sufficiently large, then the following holds.  Let
$\Sigma\in\mc{N}^{N-2}$, let $0\neq\sigma\in V_{i,j,k}$, and let $v$ be
a vertex of $\tau(\Sigma)$.  If one rotates the corresponding branch
point in the $S^1$ direction by angle $\varphi\in\R$, and if the
resulting change in the argument of $\sigma_i/\sigma_j$ is $r\in\R$, then
\[
|r - \varphi r(\sigma_i/\sigma_j,v)| < \varepsilon.
\]
\end{lemma}

\subsection{Beginning the computation of $n_{N-2}$}

The integer $n_{N-2}$ that we want to compute can be decomposed as a
sum as follows.  Recall that each component of $\mc{N}^{N-2}$
determines a data set $(E_1,\ldots,E_{N_+})$, where $E_i$ denotes the
set of negative ends that can be reached by downward paths starting at
the $i^{th}$ positive end.  Since the outgoing edges incident to
$v_1,\ldots,v_{N_+-1}$ have positive multiplicities, it follows from
Lemma~\ref{lem:notethat}(a) that
\[
E\eqdef (E_1,\ldots,E_{N_+}) \in \mc{E}(S).
\]
(See \S\ref{sec:EAT} for the definition of $\mc{E}(S)$.)  Given
$E\in\mc{E}(S)$, let $\mc{N}(E)$ denote the corresponding union of
components of $\mc{N}^{N-2}$.  We can then write
\begin{equation}
\label{eqn:nsum}
n_{N-2} = \sum_{E\in\mc{E}(S)} n(E),
\end{equation}
where $n(E)$ denotes the signed count of $\Sigma\in\mc{N}(E)$ solving
the equations \eqref{eqn:simplified'}.

Observe that if $\Sigma\in\mc{N}(E)$, then the associated trivalent
tree $\tau(\Sigma)$, with the function $\rho$ forgotten, is exactly
the tree $\phi(E)$ defined in \S\ref{sec:EAT}.

Given $E\in\mc{E}(S)$, we now derive a formula for $n(E)$. To state
the formula, let $k\in\{1,\ldots,N-2\}$, and consider a branched cover
$\Sigma\in\mc{N}(E)$ and a nonzero special cokernel element
\[
\sigma^{(k)}\in
V_{\epsilon_k,\epsilon_{N-1},\epsilon_N}.
\]
Here, unlike in \eqref{eqn:psi(i)}, we are not requiring
$\sigma^{(k)}_{\epsilon_k}=1$.  Now there is a ``dominant'' end
$d(k)\in\{\epsilon_k,\epsilon_{N-1},\epsilon_N\}$ whose contribution
to $\mathfrak{s}_{\C}(C)(\sigma^{(k)})$ is much larger than the other
contributions, in the sense of Definition~\ref{def:dominate}.  That
is:

\begin{lemma}
Given $\Sigma\in\mc{N}^{N-2}$ and $k\in\{1,\ldots,N-2\}$, let $v$
denote the central vertex for $\epsilon_k,\epsilon_{N-1},\epsilon_N$.
Then:
\begin{description}
\item{(a)}
If $\rho(v)>0$, then $\epsilon_{N-1}\searrow \epsilon_k,\epsilon_N$.
\item{(b)}
If $\rho(v)<0$, then $\epsilon_N\searrow \epsilon_k,\epsilon_{N-1}$.
\end{description}
\end{lemma}

\begin{proof}
By Lemma~\ref{lem:notethat}(b), there is a downward path $P$ from
$\epsilon_{N-1}=N_+$ to $\epsilon_N=N_-$.  The central vertex $v$ is
somewhere on the path $P$.  Suppose that $\rho(v)>0$.  The path
$P_{v,\epsilon_k}$ must dip below the level $\rho=-R_N$, because all
vertices with $\rho>0$ are joining.  It then follows as in the proof of
Lemma~\ref{lem:CLS}(a) that $\epsilon_{N-1}\searrow
\epsilon_k,\epsilon_N$.  This proves assertion (a), and assertion (b)
follows by a symmetric argument.
\end{proof}
\noindent
Define $d(k)\eqdef\epsilon_{N-1}$ in case (a) above, and
$d(k)\eqdef\epsilon_N$ in case (b).

Next, define a square matrix $A(E)$ over $\Q$ of size $N-2$ as
follows.  The rows of $A(E)$ correspond to the ends
$\epsilon_1,\ldots,\epsilon_{N-2}$.  The columns of $A(E)$
correspond to the vertices $v_1,\ldots,v_{N-2}$.  The entries of
$A(E)$ are defined by the rotation rates
\[
A(E)_{k,l} \eqdef r\left(\sigma^{(k)}_{d(k)} / \sigma^{(k)}_{\epsilon_k},
v_l\right).
\]
Let $\op{Edge}(E)$ and $\op{Vert}(E)$ denote the sets of edges and
internal vertices respectively in the tree $\phi(E)$.

\begin{lemma}
\label{lem:det1}
If $r_1$ is sufficiently large, then for each $E\in\mc{E}(S)$, we have
\begin{equation}
\label{eqn:det1}
n(E)=(-1)^{\frac{(N-2)(N-3)}{2}+(N_+-1)} \det(A(E))\prod_{e\in
\op{Edge}(E)}m(e).
\end{equation}
\end{lemma}

\begin{proof}
There is a natural action of $\R^{N-2}$ on $\mc{N}(E)$ that rotates
the $N-2$ branch points in the $S^1$ direction at speed $2\pi$.  The
kernel of this action is a nondegenerate lattice $\Lambda(E) \subset
\Z^{N-2}$.  In fact, the proof of Lemma~\ref{lem:covering} shows that
$\Lambda(E)$ is the kernel of the homomorphism
\[
\bigoplus_{\op{Vert}(E)}\Z \longrightarrow
\bigoplus_{e\in\op{Edge}(E)}\Z/m(e)
\]
that sends (the generator corresponding to) a vertex $v$ to the sum of
the outgoing edges of $v$ minus the sum of the incoming edges of $v$.
Thus we can identify
\begin{equation}
\label{eqn:CATE0}
\mc{N}(E) \simeq \bigsqcup_{\pi_0\mc{N}(E)} \R^{N-2}/\Lambda(E).
\end{equation}
By Lemma~\ref{lem:covering}, we have
\begin{equation}
\label{eqn:CATE1}
\det(\Lambda(E)) \cdot \left|\pi_0\mc{N}(E)\right| =
  \prod_{e\in \op{Edge}(E)}m(e).
\end{equation}

Now define a map
\[
\begin{split}
f:\mc{N}(E) & \longrightarrow (S^1)^{N-2},\\ \Sigma & \longmapsto
\left\{\op{arg}\left(
\sigma^{(k)}_{d(k)} / \sigma^{(k)}_{\epsilon_k}
\right)\right\}_{k=1}^{N-2}.
\end{split}
\]
By the domination condition \eqref{eqn:d2}, the map $f$ is homotopic
to the map sending
\[
\Sigma \longmapsto \left\{\op{arg}\left(
\mathfrak{s}_{\C}(\Sigma)(\sigma^{(k)}) / \sigma^{(k)}_{\epsilon_k}
\right)\right\}_{k=1}^{N-2}.
\]
Therefore the count of solutions to \eqref{eqn:simplified'} on $\mc{N}(E)$
is given by
\begin{equation}
\label{eqn:CATE2}
n(E) = 
\deg(f).
\end{equation}

On the other hand, by Lemma~\ref{lem:RBP}, if $r_1$ is sufficiently
large, then under the identification \eqref{eqn:CATE0} the map $f$ is
homotopic to the linear map $A(E)$ on each component of $\mc{N}(E)$.
Therefore
\begin{equation}
\label{eqn:CATE3}
\deg(f) = 
(-1)^{\frac{(N-2)(N-3)}{2} + (N_+-1)} \cdot  \det(\Lambda(E)) \cdot
  \left|\pi_0\mc{N}(E)\right| \cdot \det(A(E)).
\end{equation}
Here the sign arises from the orientation convention for
$\mc{N}^{N-2}$ in \S\ref{sec:CO}.

Combining equations \eqref{eqn:CATE1}, \eqref{eqn:CATE2}, and
\eqref{eqn:CATE3} proves the lemma.
\end{proof}

\subsection{Calculating the determinant}

Here is where we stand.  By equations \eqref{eqn:nn2},
\eqref{eqn:nsum}, and \eqref{eqn:det1}, we have
\begin{equation}
\label{eqn:HWS}
\#\mathfrak{s}_0^{-1}(0) = (-1)^{N_-+1} \sum_{E\in\mc{E}(S)} \det(A(E))\prod_{e\in
\op{Edge}(E)}m(e).
\end{equation}
By equation \eqref{eqn:HWS} and Lemma~\ref{lem:EATSum}, the following
lemma will finish off the proof of Proposition~\ref{prop:count}.  In
the statement of this lemma, recall from \S\ref{sec:ATEP} that
$P_{\phi(E)}$ denotes the canonical edge pairing on the admissible
tree $\phi(E)$, and $W_\theta(\phi(E),T_{\phi(E)})$ denotes the
associated positive integer weight.

\begin{lemma}
For each $E\in\mc{E}(S)$, we have
\[
\det(A(E))\prod_{e\in
\op{Edge}(E)}m(e) = (-1)^{N_--1} W_\theta(\phi(E),P_{\phi(E)}).
\]
\end{lemma}

\begin{proof}
It follows from the definitions that the canonical edge pairing
$P_{\phi(E)}$ on $\phi(E)$ is given by $e_{v_k}^\pm=e_k^\pm$, where
the edges $e_k^\pm$ are specified in Definitions~\ref{def:Nk+} and
\ref{def:Nk-}.

Now define a matrix $B$ as follows.  Let $A_l$ denote the $l^{th}$ row of
$A\eqdef A(E)$.  Then the rows of $B$ are given by the following
prescription.
\begin{itemize}
\item
If $i\in\{1,\ldots,N_+-1\}$, then
\[
B_i \eqdef \left\{\begin{array}{cl} A_i - A_k, & v_i^- =
\epsilon_k\neq\epsilon_N,\\
A_i, & v_i^- = \epsilon_N.
\end{array}\right.
\]
\item
If $k\in\{N,\ldots,N-2\}$, then (cf.\ Lemma~\ref{lem:notethat}(a),(b))
\[
B_k \eqdef \left\{\begin{array}{cl} A_k - A_i, & \epsilon_k\in E_i',
\; i<N_+,\\
A_k, & \epsilon_k\in E_{N_+}'.
\end{array}\right.
\]
\end{itemize}
By equation \eqref{eqn:weight}, to prove the lemma it suffices to
prove (i)--(iii) below:
\begin{description}
\item{(i)} $\det(A)=\det(B)$.
\item{(ii)} $B$ is lower triangular, for a suitable reordering of
$\{1,\ldots,N-2\}$.
\item{(iii)} The $l^{th}$ diagonal entry $B_{l,l}$ of $B$ is given by
\end{description}
\begin{equation}
\label{eqn:Bll}
\frac{m(e_l^-)\ceil{m(e_l^+)\theta} -
m(e_l^+)\floor{m(e_l^-)\theta}}{m(e_l^+)m(e_l^-)} = \left\{\begin{array}{cl}
B_{l,l}, & l=1,\ldots,N_+-1,\\
-B_{l,l}, & l=N_+,\ldots,N-2.
\end{array}\right.
\end{equation}

{\em Proof of (i):\/} The matrix $B$ is obtained from $A$ by
performing the following row operations for $i=1,\ldots,N_+-1$ in
order:
\begin{itemize}
\item
For each $k$ such that $\epsilon_k\in E_i'$, subtract the $i^{th}$ row
from the $k^{th}$ row.
\item
If $v_i^-=\epsilon_k$ with $k\neq N$, then subtract the $k^{th}$ row
(which has not yet been modified since $\epsilon_k\in E_j'$ for some
$j>i$) from the $i^{th}$ row.
\end{itemize}

{\em Proof of (ii):\/} We claim that the matrix $B$ is lower triangular
if one lists the numbers $1,\ldots,N-2$, which index the rows and
columns of $B$, in the order
\begin{equation}
\label{eqn:LTO}
\{k\mid \epsilon_k\in E_1'\},1,\ldots,\{k\mid \epsilon_k\in
E_{N_+-1}'\},N_+-1,
\{k\mid \epsilon_k\in E_{N_+}'\}.
\end{equation}
Here the set $\{k\mid \epsilon_k\in E_i'\}$ is listed in
increasing order of $k$ for each $i$.

To prove lower triangularity, we first investigate the $k^{th}$ row of
$B$ when $\epsilon_k\in E_i'$ and $i<N_+$.
By the definition of $A$,
\[
\begin{split}
A_{i,l} & = r(\sigma^{(i)}_{d(i)}/\sigma^{(i)}_i,v_l),\\
A_{k,l} & = r(\sigma^{(k)}_{d(k)}/\sigma^{(k)}_{\epsilon_k},v_l).
\end{split}
\]
We now calculate these rotation rates using Definition~\ref{def:rr}
and Lemma~\ref{lem:specialWinding}.  First note that
the dominant ends $d(i)$ and $d(\epsilon_k)$ are equal.  The reason is
that the central vertex for $i, \epsilon_{N-1},
\epsilon_N$ is the same as the central vertex for $\epsilon_k,
\epsilon_{N-1}, \epsilon_N$, because the paths from $\epsilon_k$ or
$i$ to $P_{\epsilon_{N-1},\epsilon_N}$ both pass through $v_k$.

More precisely, the path $P_{i,d(i)}$ passes first through the
vertices $v_j$ for $\epsilon_j\in E_i'$ in increasing order, then
through the vertex $v_i$.  If $v_i^-\neq\epsilon_N$, then at $v_i$ the
path $P_{i,d(i)}$ turns (at least temporarily) upward and passes
through some additional internal vertices which are all in $\tau_i^+$;
otherwise the path $P_{i,d(i)}$ stays downward, and any additional
internal vertices on this path are in $\tau_{i-1}^-$.  Likewise, the
path $P_{\epsilon_k,d(k)}$ possibly first passes through some vertices
in $\tau_{i-1}^-$, then hits the vertex $v_k$, and then agrees with
the rest of the path $P_{i,d(i)}$.

Since the central vertex for $i, \epsilon_{N-1},
\epsilon_N$ is the same as the central vertex for $\epsilon_k,
\epsilon_{N-1}, \epsilon_N$, it follows
by Definition~\ref{def:rr} and Lemma~\ref{lem:specialWinding} that
$A_{i,l}=A_{k,l}$ whenever the paths $P_{i,d(i)}$ and
$P_{\epsilon_k,d(\epsilon_k)}$ either both avoid $v_l$, or both pass
through $v_l$ along the same ordered pair of edges.  By the above
description of these two paths, this can fail only if
$v_l\in\tau_{i-1}^-$, or $\epsilon_l\in E_i'$ with $l\le k$.  It
follows that the $k^{th}$ row of $B$ has the required form for lower
triangularity with respect to the ordering \eqref{eqn:LTO}.  Similar
arguments show that all other rows of $B$ have the required form.

\medskip

{\em Proof of (iii):\/} We now prove equation \eqref{eqn:Bll} in several
cases.  In these calculations recall that $e_l^0$ denotes the edge of
$v_l$ that is neither $e_l^+$ nor $e_l^-$.

Suppose first that $k\in\{N_+,\ldots,N-2\}$.  Then equation
\eqref{eqn:Bll} for $l=k$ asserts that
\begin{equation}
\label{eqn:Bkk}
B_{k,k} = -\frac{\ceil{m(e_k^+)\theta}}{m(e_k^+)} +
\frac{\floor{m(e_k^-)\theta}}{m(e_k^-)}.
\end{equation}
If $\epsilon_k\in E_i'$ with $i<N_+$, then
by the definition of $A$ and Lemma~\ref{lem:specialWinding}, and using
the descriptions of the paths $P_{i,d(i)}$ and $P_{\epsilon_k,d(k)}$
from the proof of part (ii), we obtain
\[
\begin{split}
A_{k,k} &= \frac{\floor{m(e_k^-)\theta}}{m(e_k^-)} -
\frac{\ceil{m(e_k^0)\theta}}{m(e_k^0)},\\
A_{i,k} &= \frac{\ceil{m(e_k^+)\theta}}{m(e_k^+)} -
\frac{\ceil{m(e_k^0)\theta}}{m(e_k^0)}.
\end{split}
\]
Subtracting these two equations gives \eqref{eqn:Bkk}.  If
$\epsilon_k\in E_{N_+}'$, then \eqref{eqn:Bkk} holds since
$B_{k,k}=A_{k,k}$ and $d(k)=\epsilon_{N-1}$.

Finally, suppose $i\in\{1,\ldots,N_+-1\}$.  Then equation
\eqref{eqn:Bll} for $l=i$ is
\begin{equation}
\label{eqn:Bii}
B_{i,i} = \frac{\ceil{m(e_i^+)\theta}}{m(e_i^+)} -
\frac{\floor{m(e_i^-)\theta}}{m(e_i^-)}.
\end{equation}
By Lemma~\ref{lem:specialWinding} and the definition of $A$, if 
$v_i^-=\epsilon_N$ then $A_{i,i}$ is given by the right hand side of
\eqref{eqn:Bii}.  On the other hand, if
$v_i^- =\epsilon_k \neq \epsilon_N$, then
\[
A_{i,i} = 
\frac{\ceil{m(e_i^+)\theta}}{m(e_i^+)} -
\frac{\floor{m(e_i^0)\theta}}{m(e_i^0)}.
\]
Also, since $\epsilon_k\in E_j'$ for some $j>i$, we have
\[
A_{k,i} = 
\frac{\floor{m(e_i^-)\theta}}{m(e_i^-)} -
\frac{\floor{m(e_i^0)\theta}}{m(e_i^0)}. 
\]
The above calculations imply \eqref{eqn:Bii}.
\end{proof}

\begin{remark}
\label{rem:EPN}
One might try to give a more direct proof of
Proposition~\ref{prop:count} as follows.  If $\mathfrak{s}_\C(\Sigma)=0$,
then ``generically'' the tree $\tau(\Sigma)$ is trivalent, and given a
nonzero special cokernel element $\sigma\in V_{i,j,k}$, in the equation
\[
\mathfrak{s}_\C(\Sigma)(\sigma) = \pm \overline{\gamma_i}\sigma_i
\pm\overline{\gamma_j}\sigma_j \pm \overline{\gamma_k}\sigma_k=0,
\]
one term is much smaller than the other two.  The two larger terms
specify two distinguished edges incident to the central vertex $v$ for
$i$, $j$, and $k$.  One can check that these two edges depend only on
$v$ and define an edge pairing on $\tau(\Sigma)$, modulo the choice of
which distinguished edge is $e_v^+$ and which is $e_v^-$.  Moreover,
similarly to the above calculations, the count of solutions with this
tree and edge pairing is given by plus or minus the weight in
Definition~\ref{def:weight}.  Thus one finds that
$\#\mathfrak{s}_\C^{-1}(0)$ is naturally given by a sum over certain trees
with edge pairings of their corresponding weights.  However the sum
that arises is sometimes different than the sum over admissible trees
in Lemma~\ref{lem:IFT}, and the combinatorics of this approach seems difficult.
\end{remark}

\section{Detailed analysis of the obstruction bundle}
\label{sec:approx}

In this section, as in \S\ref{sec:obstructionBundle}, fix positive
integers $a_1,\ldots,a_{N_+}$ and $a_{-1},\ldots,a_{-N_-}$ satisfying
\eqref{eqn:M}, fix an admissible almost complex structure $J$ on $\R\times Y$,
and fix an embedded elliptic Reeb orbit $\alpha$ with monodromy angle
$\theta\in\R\setminus\Q$ satisfying \eqref{eqn:kappa1}.  Let
$\mc{M}\eqdef
\mc{M}(a_1,\ldots,a_{N_+}\mid a_{-1},\ldots,a_{-N_-})$ denote the
moduli space of branched covers of the cylinder $\R\times S^1$ from
Definition~\ref{def:M}, and given $\Sigma\in\mc{M}$ recall the
operator $D_\Sigma$ defined in \S\ref{sec:DSigma}.  As usual, identify
an element of $\Coker(D_{\Sigma})$ with a smooth, square integrable
$(0,1)$ form $\sigma$ on $\Sigma$ satisfying $D_{\Sigma}^*\sigma=0$,
and away from the ramification points use $d\overline{z}$ to identify
$\sigma$ with a complex function.

In this section we give the previously deferred proof of
Proposition~\ref{prop:approx1}, which describes the approximate
behavior of nonvanishing cokernel elements away from the ramification
points.  We also prove a result on the approximate behavior of
nonvanishing cokernel elements near isolated clusters of ramification
points.  The latter result is stated in
\S\ref{sec:ICR}, and the proofs of both results are given in \S\ref{sec:PAR}.
 In \S\ref{sec:RSProof} and \S\ref{sec:RR} we use these results to give the
previously deferred proofs of Proposition~\ref{prop:RS} and
Lemma~\ref{lem:RBP}.

\subsection{Isolated clusters of ramification points}
\label{sec:ICR}

We now state a result which asserts, roughly, that the behavior of a
nonvanishing cokernel element near an isolated cluster of ramification
points does not depend much on the nature of the distant ramification
points.

We need the following preliminary definitions.  Let
$\pi:\Sigma\to\R\times S^1$ be a branched cover in $\mc{M}$.  Recall
from \S\ref{sec:obstructionBundle} that $\Sigma$ determines a tree
$\tau(\Sigma)$ with a metric and a map $p:\Sigma\to\tau(\Sigma)$.

\begin{definition}
\label{def:cluster}
A nonempty set $Z$ of ramification points in $\Sigma$ is a {\em
cluster\/} if there is a connected set $B\subset\tau(\Sigma)$, such
that a ramification point $z\in\Sigma$ is in $Z$ if and only if
$p(z)\in B$.  In this case let ${\Sigma}_Z$ denote the
branched cover of $\R\times S^1$ obtained by attaching half-infinite
cylinders to the boundary circles of $p^{-1}(B)$.  The {\em
diameter\/} of $Z$ is the diameter of the set $p(Z)$ in
$\tau(\Sigma)$.  For a real number $R>0$, the cluster $Z$ is {\em
$R$-isolated\/} if every vertex in $p(Z)$ has distance at least $R$
from all vertices of $\tau(\Sigma)$ not in $p(Z)$.
\end{definition}

Note that there is a canonical identification between a cluster of
ramification points $Z$ and the set of ramification points in
$\Sigma_Z$.  Also, if $Z$ is $R$-isolated, then there is a canonical
identification between the set of points in $\Sigma$ within distance
$R$ of a ramification point in $Z$, and the set of points in
$\Sigma_Z$ within distance $R$ of a ramification point in $\Sigma_Z$.

\begin{definition}
A $(0,1)$-form $\sigma$ on $\Sigma$ has {\em exponential growth\/} if
there exists a constant $c$ such that $|\sigma| \le
c\exp(c|{\pi}^*s|)$ at every point in $\Sigma$.  Define
$\widetilde{\Coker}(D_\Sigma)$ to be the space of $(0,1)$-forms with
exponential growth on $\Sigma$ that are annihilated by $D_\Sigma^*$.
\end{definition}

\begin{proposition}
\label{prop:approx2}
Given $r,\varepsilon_0>0$, there exists $R>1/\varepsilon_0$ such that
the following holds.  Let $\Sigma\in\mc{M}$, let $Z$ be an
$R$-isolated cluster of ramification points in $\Sigma$ of diameter
$\le 2r$, and let $\sigma\in\Coker(D_\Sigma)$ be nonvanishing.  Then
there exists a nonvanishing
${\sigma}_Z\in\widetilde{\Coker}(D_{\Sigma_Z})$ such that
$|\sigma-\sigma_Z|\le
\varepsilon_0|\sigma|$ at all points in $\Sigma$ within distance
$1/\varepsilon_0$ of a ramification point in $Z$.
\end{proposition}

To say more about the forms
$\sigma_Z\in\widetilde{\Coker}(D_{\Sigma_Z})$ that can arise, first
note the following basic lemma:

\begin{lemma}
\label{lem:AK}
Fix a positive integer $m$ and an integer $\eta$.  Then there exists
$\kappa>0$ with the following property.  Let $\sigma$ be a complex
function on $[0,\infty)\times\R/2\pi m\Z$ which is annihilated by
$\partial_s-L_m$, is nonvanishing with winding number $\eta$, and has
exponential growth.  Then there is a normalized ($L^2$-norm $1$)
eigenfunction $\gamma$ of $L_m$ with eigenvalue $E_\gamma$ and winding
number $\eta$, and constants $\sigma_\gamma\neq 0$ and $c_\sigma$,
such that
\[
|\sigma(s,t) - \sigma_\gamma\exp(E_\gamma s)\gamma(t)| \le
 c_\sigma\exp((E_\gamma-\kappa)s)
\]
for all $(s,t)$ with $s\ge 0$.
\end{lemma}

\begin{proof}
This follows by writing $\sigma(s,t) =
\sum_\gamma\sigma_\gamma\exp(E_\gamma s)\gamma(t)$, where the sum is
over an orthonormal basis of $L^2(\R/2\pi m\Z;\R^2)$ consisting of
eigenfunctions $\gamma$ of $L_m$ with eigenvalues $E_\gamma$.
\end{proof}

Now suppose $\sigma\in\widetilde{\Coker}(D_\Sigma)$ is nonvanishing.
As in \S\ref{sec:DSigma}, let $\eta_i^+$ and $\eta_j^-$ denote the
winding numbers of $\sigma$ around the positive and negative ends of
$\Sigma$.  By equation \eqref{eqn:ACZ}, we must have
\begin{equation}
\label{eqn:nvw}
\sum_{i=1}^{N_+}\eta_i^+ - \sum_{j=-1}^{-N_-}\eta_j^- = N-2.
\end{equation}
Given integers $\eta_i^+$ and $\eta_j^-$ satisfying \eqref{eqn:nvw},
let $V(\eta_1^+,\ldots,\eta_{N_+}^+\mid
\eta_{-1}^-,\ldots,\eta_{-N_-}^-)$ denote the vector space of
$\sigma\in\widetilde{\Coker}(D_\Sigma)$ such that either $\sigma=0$,
or $\sigma$ is nonvanishing with winding numbers $\eta_i^+$ and
$\eta_j^-$.  Calculations similar to those in
Lemmas~\ref{lem:windingBounds} and
\ref{lem:SCE}, using Lemma~\ref{lem:AK}, show that
\begin{equation}
\label{eqn:dimRV}
\dim_\R V(\eta_1^+,\ldots,\eta_{N_+}^+\mid
\eta_{-1}^-,\ldots,\eta_{-N_-}^-) = 2.
\end{equation}
In particular, there is a vector bundle
\begin{equation}
\label{eqn:vvb}
\mc{V}(\eta_1^+,\ldots,\eta_{N_+}^+\mid
\eta_{-1}^-,\ldots,\eta_{-N_-}^-)\longrightarrow\mc{M},
\end{equation}
whose fiber over $\Sigma\in\mc{M}$ is the vector space
$V(\eta_1^+,\ldots,\eta_{N_+}^+\mid
\eta_{-1}^-,\ldots,\eta_{-N_-}^-)$ associated to $\Sigma$.

\begin{remark}
Also in connection with Lemma~\ref{lem:AK}, one of the difficulties in
proving Proposition~\ref{prop:approx1} is that there is no a
priori upper bound on the ratio $c_\sigma/|\sigma_\gamma|$.  For
example, fix $\theta\in(0,1)$ and take $S(t)=-\theta$ and $m=1$.  Let
$a$ be a nonnegative real number and take $\sigma = e^{-\theta
s}\exp(ae^{-s+it})$ on $[0,\infty)\times S^1$.  This $\sigma$ is
nonvanishing and square integrable with $\eta=0$, so
Lemma~\ref{lem:AK} gives $\gamma(t)=1\sqrt{2\pi}$.  It is then easy to
check that the smallest possible value of $c_\sigma/|\sigma_\gamma|$
limits to $\infty$ as $a\to\infty$.

The following even worse situation can occur over a compact cylinder.
Again fix $\theta\in(0,1)$ and take $S(t)=-\theta$ and $m=1$.  Let $a$
be a nonnegative real number and define $\sigma:
\R\times S^1\to\C$ by
\[
\sigma(s,t) \eqdef e^{-\theta s} \exp\left(\frac{ia}{2}\left(e^{s-it}
+ e^{-s+it}\right)\right),
\]
This is nonvanishing, is annihilated by $\partial_s-L$, and has
winding number $\eta=0$.  Even so, there exists $a$ such that
$\sigma(0,\cdot)$ has no constant term in its Fourier series, so that
$\Pi_W\sigma(0,\cdot)=0$.
\end{remark}

\subsection{Proof of the approximation results}
\label{sec:PAR}

We now prove Propositions~\ref{prop:approx1} and \ref{prop:approx2}
together.  If either of these propositions fails, then we can find
constants $\varepsilon_0,r>0$, and a sequence of pairs
$\{(\Sigma_k,\sigma_k)\}_{k=1,2,\ldots}$ such that the conclusions of
the propositions do not all hold for
$(\Sigma,\sigma)=(\Sigma_k,\sigma_k)$ when $R<k$.  Hence to prove
Propositions~\ref{prop:approx1} and \ref{prop:approx2}, it is enough
to prove the following statement:
\begin{itemize}
\item
Consider a sequence $\{(\Sigma_k,\sigma_k)\}_{k=1,2,\ldots}$ where
$\pi_k:\Sigma_k\to\R\times S^1$ is a branched cover in $\mc{M}$ and
$\sigma_k\in\Coker(D_{\Sigma_k})$ is nonvanishing.  Let $r$ be given.
Then we can pass to a subsequence (again indexed by $k=1,2,\ldots$)
such that for all $\varepsilon_0>0$,
there exists $R$ such that the conclusions of
Proposition~\ref{prop:approx1} and
\ref{prop:approx2} hold for $(\Sigma,\sigma)=(\Sigma_k,\sigma_k)$
whenever $k$ is sufficiently large.
\end{itemize}
We now prove the above statement in 7 steps.  The strategy is to pass
to a subsequence with appropriate convergence properties, and then use
estimates on the limit to produce $R$ from $\varepsilon_0$.

\medskip

{\em Step 1.\/}
We begin by passing to a subsequence so that the
sequence $\{\Sigma_k\}$ has certain convergence properties.

By Lemma~\ref{lem:convergence}, we can pass to a
subsequence so that the sequence $\{[\Sigma_k]\}$ in $\mc{M}/\R$
converges, in the sense of Definition~\ref{def:convergence}, to an
element $(T;[\Sigma_{*1}],\ldots,[\Sigma_{*p}])
\in\overline{\mc{M}/\R}$.  Fix $\Sigma_{kj}$ and $\Phi_{k,e}$ as in
Definition~\ref{def:convergence}, and carry over the other notation
from Definition~\ref{def:convergence}.  By passing to a further
subsequence and increasing $r$ if necessary, we may assume that:
\begin{itemize}
\item
$\Sigma_{kj}$ is a component of the $\pi_k$-inverse image of a
subcylinder in $\R\times S^1$ of length $2k$.
\item
If $z\in\Sigma_{kj}$ is a ramification point, then $\pi_k(z)$ has
distance $\le r$ from the center of the subcylinder
$\pi_k(\Sigma_{kj})$.
\end{itemize}

Since
$\lim_{k\to\infty}T_{-s_{kj}}(\widehat{\Sigma}_{kj})=\Sigma_{*j}$, we
can, possibly after passing to a further subsequence, choose
diffeomorphisms of the domains
$\Psi_{kj}:{\Sigma}_{*j}\to\widehat{\Sigma}_{kj}$ such that:
\begin{itemize}
\item
$T_{s_{kj}}\circ\pi_{*j}\circ\Psi_{kj}^{-1}$ agrees with the projection
$\widehat{\Sigma}_{kj}\to\R\times S^1$ at all points in
$\widehat{\Sigma}_{kj}$ that have distance $\ge 1$ from the
ramification points in $\widehat{\Sigma}_{kj}$.
\item
Let $\Psi_{kj}^{0,1}$ denote the composition of the pullback
$\Psi_{kj}^*:T^*_\C\widehat{\Sigma}_{kj}\to T^*_\C \Sigma_{*j}$ with
orthogonal projection $T^*_\C\Sigma_{*j}\to T^{0,1}\Sigma_{*j}$.  Then the
sequence of differential operators $\left\{\Psi_{kj}^{0,1}\circ
D_{\widehat{\Sigma}_{kj}}^*\right\}_{k=1,2,\ldots}$ converges to
$D_{\Sigma_{*j}}^*$.
\end{itemize}

The following notation will be used below.  Choose a $k$-independent
number $a>2r+1$. Restrict attention to $k\ge a$.  Let $V_{kj}\subset
\Sigma_{kj}$ denote $\pi_k^{-1}$ of the set of points with $s_{kj}-a< s <
s_{kj}+a$.  Let $U_{kj}\subset V_{kj}$ denote $\pi_k^{-1}$ of the set
of points with $s_{kj}-a+1<s<s_{kj}+a-1$.  Let $\pi_{*j}$ denote the
projection $\Sigma_{*j}\to\R\times S^1$, let $V_{*j}\subset
\Sigma_{*j}$ denote $\pi_{*j}^{-1}$ of the set of points with
 $|s|<a$, and let $U_{*j}\subset V_{*j}$ denote $\pi_{*j}^{-1}$ of the
set of points with $|s|<a-1$.

\medskip

{\em Step 2.\/}  We now pass to a further subsequence so that the
sequence $\{\sigma_k\}$ has certain convergence properties.

To start, normalize the $\sigma_k$'s to have $L^2$ norm $1$.  Define
$\theta_{kj}$ to be the $L^2$ norm of $\sigma_k$ over $V_{kj}$.

\begin{lemma}
\label{lem:asub}
For each $j$, there is a smooth $(0,1)$-form $\sigma_{*j}$ on $V_{*j}$
which is annihilated by $D_{\Sigma_{*j}}^*$ such that:
\begin{description}
\item{(a)}
The sequence of $(0,1)$-forms $
\left\{\Psi_{kj}^{0,1} \left(\theta_{kj}^{-1}
\sigma_k|_{V_{kj}}\right)
\right\}_{k=a,a+1,\ldots}
$ has a subsequence that converges in the $C^\infty$
topology\footnote{Here and below, `convergence in the $C^\infty$
topology' means convergence in the $C^n$ topology on any compact set
for any integer $n$.} to $\sigma_{*j}$.
\item{(b)}
$\sigma_{*j}$ is nonvanishing.
\end{description}
\end{lemma}

\begin{proof}
A standard compactness argument using a priori elliptic estimates
finds a subsequence of the sequence in (a) converging to a smooth
$(0,1)$-form $\sigma_{*j}$ on $V_{*j}$ that is annihilated by
$D_{\Sigma_{*j}}^*$.
The $(0,1)$-form $\sigma_{*j}$ is nonvanishing provided that it is not
identically zero, because it is the $C^0$ limit of a sequence of
nonvanishing $(0,1)$-forms, and any zero of $\sigma_{*j}$ must have
negative multiplicity.  Thus it remains only to prove that
$\sigma_{*j}$ is not identically zero.

Suppose to the contrary that $\sigma_{*j}=0$.  By elliptic estimates,
this assumption implies that
for any neighborhood $N$ of the boundary of $V_{*j}$,
\begin{equation}
\label{eqn:olddagger}
\lim_{k\to\infty}\int_{V_{*j}\setminus N}
\left|\Psi_{kj}^{0,1}(\theta_{kj}^{-1}\sigma_k)\right|^2=0.
\end{equation}

Now pass to a subsequence so that for each edge $e$ of the tree $T$
incident to the $j^{th}$ internal vertex, the $L^2$ norm of
$\Psi_{kj}^{0,1}(\sigma_k/\theta_{kj})$ over the component of
$V_{*j}\setminus U_{*j}$ corresponding to $e$ converges as
$k\to\infty$ some $c_{j,e}\ge 0$.  By \eqref{eqn:olddagger} with
$N=V_{*j}\setminus U_{*j}$, we must have $\sum_ec_{j,e}=1$.  Hence
there is an edge $e$ of the tree $T$ adjacent to the $j^{th}$ internal
vertex with $c_{j,e}>0$.

We will now show that if $j$ is an internal vertex with
$\sigma_{*j}=0$, and if $e$ is an edge of $T$ incident to $j$ with
$c_{j,e}>0$, then:
\begin{description}
\item{(i)}
$e$ is an internal edge.
\item{(ii)}
Let $j'\neq j$ denote the other internal vertex of $T$ incident to
$e$; then $\sigma_{*j'}=0$.
\item{(iii)}
If $e'$ is an edge of $T$ incident to $j'$ with $c_{j',e'}>0$, then
$e\neq e'$.
\end{description}
By induction using (i) and (ii), we can find an infinite sequence of
internal vertices $j_0=j, j_1=j', j_2,\ldots$ of $T$, and an infinite
sequence of edges $e_0=e, e_1=e', e_2,\ldots$ such that
$\sigma_{*j_i}=0$; the edge $e_i$ is incident to $j_i$ and $j_{i+1}$;
and $c_{j_i,e_i}>0$.  Then property (iii) implies $e_i\neq e_{i+1}$.
Since $T$ is a tree, this will give the desired contradiction.

Proof of (i): For each $k$, let $\mc{E}_k$
denote the component cylinder of $\Sigma_k\setminus
\bigcup_{j'}U_{kj'}$ corresponding to $e$.
Without loss of generality, $s\le 1$ on $\mc{E}_k$, with $s=1$
denoting the boundary circle of $U_{kj}$ and $s=0$ the boundary circle
of $V_{jk}$.  By \eqref{eqn:olddagger} again,
\begin{equation}
\label{eqn:boundaryV}
\begin{split}
\lim_{k\to\infty} \theta_{kj}^{-2}\int_{1/2}^1
\|\sigma_k|_{s=\tau}\|^2d\tau & = 0,\\
\lim_{k\to\infty}\theta_{kj}^{-2}\int_{0}^{1/2}\|\sigma_k|_{s=\tau}\|^2d\tau
& =
c_{j,e}^2.
\end{split}
\end{equation}
To be more explicit, expand $\sigma_k|_{\mc{E}_k}$, regarded as a
complex function, in terms of eigenfunctions of $L_m$ as
\begin{equation}
\label{eqn:expandsigma}
\sigma_k|_{\mc{E}_k}(s,t) = \sum_\gamma\sigma_{k\gamma}\exp(E_\gamma
s)\gamma(t).
\end{equation}
It follows from
\eqref{eqn:boundaryV} and \eqref{eqn:expandsigma} that for every real
number $\Lambda$, there exists $c_\Lambda>0$ such that for all
$\varepsilon>0$, if $k$ is sufficiently large then
\begin{equation}
\label{eqn:EL}
\sum_{E_\gamma>\Lambda}
|\sigma_{k\gamma}|^2 \le \varepsilon \theta_{kj}^2,
\quad\quad\quad\quad
\sum_{E_\gamma\le\Lambda} |\sigma_{k\gamma}|^2 \ge
c_\Lambda(1-\varepsilon)c_{j,e}^2\theta_{kj}^2.
\end{equation}
Taking $\Lambda<0$ in the right most inequality shows that $\mc{E}_k$
is compact.  This is because if $\mc{E}_k$ is not compact, then square
integrability of $\sigma_k$ requires that $\sigma_{k\gamma}=0$ when
$E_\gamma<0$.

Proof of (ii): Let $s_k<0$ denote the value of $s$ on the boundary
circle of $\mc{E}_k$ in $V_{kj'}$.  Let $m$ denote the multiplicity of
the edge $e$, and let $\Pi_{\le \Lambda}$ denote the orthogonal
projection in $L^2(\R/2\pi m\Z;\R^2)$ onto the span of the
eigenfunctions of $L_m$ with eigenvalue $\le \Lambda$.  Likewise let
$\Pi_{>\Lambda}$ denote the projection onto the sum of the eigenspaces
with eigenvalues $>\Lambda$.  It follows from \eqref{eqn:EL} that if
$k$ is sufficiently large then
\begin{equation}
\label{eqn:DTV}
\frac{\|\Pi_{>\Lambda} \sigma_k|_{s=s_k}\|_2}
{\|\Pi_{\le\Lambda} \sigma_k|_{s=s_k}\|_2}
=
\frac{\sum_{E_\gamma > \Lambda}|\sigma_{k\gamma}|^2\exp(2E_\gamma
s_k)}
{\sum_{E_\gamma \le \Lambda}|\sigma_{k\gamma}|^2\exp(2E_\gamma s_k)}
\le \exp(\kappa_{\Lambda}s_k)
\end{equation}
where $\kappa_\Lambda$ is a positive constant.  By the convergence in
(a) to $\sigma_{*j'}$, it follows from \eqref{eqn:DTV} that
$\sigma_{*j'}=0$ on the boundary circle of $V_{*j'}\setminus U_{*j'}$
corresponding to the edge $e$, and hence on all of $V_{*j'}$, because
$\Lambda$ can be taken arbitrarily negative and $s_k\to-\infty$ as
$k\to\infty$.

Proof of (iii): The $L^2$ norm of $\theta_{kj'}^{-1}\sigma_k$ over
$\mc{E}_k\cap V_{kj'}$ must converge to zero, because otherwise the
analogue of \eqref{eqn:EL} for $j'$, in which the inequalities on the
eigenvalues are reversed, would contradict \eqref{eqn:DTV}.
\end{proof}

Now pass to a subsequence such that the convergence in
Lemma~\ref{lem:asub} holds for each $j$.  This convergence (or an
argument independent of Lemma~\ref{lem:asub} using winding bounds)
allows us to pass to a further subsequence such that for each edge $e$
of the tree $T$, the winding number of $\sigma_k$ around the component
of $\Sigma_k\setminus\bigcup_k U_{kj}$ corresponding to $e$ does not
depend on $k$.

\medskip

{\em Step 3.\/} We now show that $\sigma_{*j}$ has an extension over
$\Sigma_{*j}$ with various nice properties.  More properties of
$\sigma_{*j}$ will be established later in Lemma~\ref{lem:extA}.

\begin{lemma}
\label{lem:extend}
For each $j$, the $(0,1)$-form $\sigma_{*j}$ extends to a smooth
$(0,1)$-form $\sigma_{*j}$ on $\Sigma_{*j}$ which is annihilated by
$D_{\Sigma_{*j}}^*$, and is such that:
\begin{description}
\item{(a)}
Let $e$ be an external edge of $T$ incident to the $j^{th}$ internal
vertex, and let $\mc{E}_k$ denote the corresponding noncompact
component of $\Sigma_k\setminus\bigcup_{j'}U_{kj'}$.  On the
corresponding component of $\Sigma_{*j}\setminus V_{*j}$, the
$(0,1)$-form $\sigma_{*j}$ is square integrable, and the limit in the
$C^\infty$ topology of the sequence
$\{\Psi_{kj}^{0,1}(\theta_{kj}^{-1}\sigma_k|_{\mc{E}_k\cap
\Sigma_{kj}})\}_{k=a,a+1,\ldots}$.  In particular, $\sigma_{*j}$ is
nonvanishing here.
\item{(b)}
$\sigma_{*j}$ has exponential growth on all of $\Sigma_{*j}$.
\end{description}
\end{lemma}

\begin{proof}
(a) We extend $\sigma_{*j}$ over the end in question as follows.
Without loss of generality, $s\le -a+1$ on $\mc{E}_k$.  Expand
$\sigma_k$ on $\mc{E}_k$ by the formula \eqref{eqn:expandsigma}, and
on the corresponding end of $\Sigma_{*j}$ where $-a<s\le-a+1$ write
\begin{equation}
\label{eqn:expandsigma*}
\sigma_{*j} = \sum_{\gamma}
\sigma_{*j\gamma} \exp(E_\gamma s)\gamma(t).
\end{equation}
By the convergence in
Lemma~\ref{lem:asub}(a) at $s=-a+1$, we have
\begin{equation}
\label{eqn:coco}
\lim_{k\to\infty} \sum_\gamma
\exp(2E_\gamma(-a+1))
|\sigma_{*j\gamma}-\theta_{kj}^{-1}\sigma_{k\gamma}|^2=0.
\end{equation}
If $E_\gamma<0$, then square integrability of $\sigma_k$ implies that
$\sigma_{k\gamma}=0$, and hence $\sigma_{*j\gamma}=0$ also by
\eqref{eqn:coco}.  Consequently \eqref{eqn:expandsigma*} defines an
extension of $\sigma_{*j}$ over the component of $\Sigma_{*j}\setminus
V_{*j}$ corresponding to $\mc{E}_k$ that has all of the required
properties.

(b) We now extend each $\sigma_{*j}$ over the rest of $\Sigma_{*j}$.
 Let $e$ be an internal edge of $T$ and let $\mc{E}_k$ denote the
 corresponding compact component of $\Sigma_k
 \setminus\bigcup_{j}U_{kj}$.  Let $j$ and $j'$ denote the upper and
 lower vertices of $e$, and suppose without loss of generality that
 $s_k \le s \le -a+1$ on $\mc{E}_k$.  Expand $\sigma_k$ on
 ${\mc{E}_k}$ where $-a\le s \le -a+1$ as in \eqref{eqn:expandsigma},
 and expand $\sigma_{*j}$ where $-a\le s \le -a+1$ as in
 \eqref{eqn:expandsigma*}.  Meanwhile, expand $\sigma_{*j'}$ where
 $a-1\le s < a$ on the component of $\Sigma_{*j'}\setminus U_{*j'}$
 that corresponds to $\mc{E}_k$ as in \eqref{eqn:expandsigma*} but
 with $j'$ replacing $j$.  By the convergence in
 Lemma~\ref{lem:asub}(a) at $s=-a+1$ and at $s=s_k$, we have
\begin{gather}
\label{eqn:coco2}
\lim_{k\to\infty}\sum_\gamma
\exp(2E_\gamma(-a+1))|\sigma_{*j\gamma} - \theta_{kj}^{-1}
\sigma_{k\gamma}|^2 = 0,\\
\label{eqn:coco3}
\lim_{k\to\infty} \sum_\gamma 
|\exp(E_\gamma(a-1))
\sigma_{*j'}-\theta_{kj'}^{-1}\exp(E_\gamma s_k)
\sigma_{k\gamma}|^2
= 0.
\end{gather}
It follows that if $\sigma_{*j\gamma}\neq 0$ and
$\sigma_{*j'\gamma'}\neq 0$ then $E_{\gamma}\ge E_{\gamma'}$, because
otherwise
\[
\frac{\sigma_{*j\gamma'}\sigma_{*j'\gamma} }{\sigma_{*j\gamma}
\sigma_{*j'\gamma'}}
= \lim_{k\to\infty} \exp((E_{\gamma'}-E_{\gamma})(-s_k+a-1))
\]
is infinite since $\lim_{k\to\infty}s_k=-\infty$.

We know from Lemma~\ref{lem:asub}(b) that $\sigma_{*j}$ and
$\sigma_{*j'}$ are nonzero, so there exist $\gamma,\gamma'$ with
$\sigma_{*j\gamma}$ and $\sigma_{*j'\gamma'}$ nonzero.  Hence there
is a smallest eigenvalue $E_+$ such that $E_+=E_{\gamma}$ with
$\sigma_{*j\gamma}\neq 0$, and a largest eigenvalue $E_-$ such
that $E_-=E_{\gamma'}$ with $\sigma_{*j'\gamma'}\neq 0$.  Hence
\eqref{eqn:expandsigma*} defines an extension of $\sigma_{*j}$ over
the negative end of $\Sigma_{*j}$ corresponding to $e$, and this
extension is a smooth $(0,1)$-form with exponential growth annihilated
by $D_{\Sigma_{*j}}$.  Likewise, $\sigma_{*j'}$ extends over the
positive end of $\Sigma_{*j'}$ corresponding to $e$ as a smooth
$(0,1)$-form with exponential growth annihilated by
$D_{\Sigma_{*j'}}$.
\end{proof}

{\em Step 4.\/} We now show that for any $\varepsilon_0>0$, there
exists $R$ such that for all $k$, the conclusions of
Proposition~\ref{prop:approx1} hold for
$(\Sigma,\sigma)=(\Sigma_k,\sigma_k)$ whenever $p(z)$ is on an
external edge $e$ of the tree $\tau(\Sigma)$.

The external edge $e$ of $\tau(\Sigma)$ corresponds to an external
edge of $T$ which we also denote by $e$.  Let $m$ denote the
multiplicity of $e$ and let $j$ denote the internal vertex of $T$
incident to $e$.  By symmetry, we may assume that the leaf incident to
$e$ is negative.  Let $\mc{E}_k$ denote the corresponding noncompact
component of $\Sigma_k\setminus\bigcup_{j'}U_{kj'}$.  As usual, there
is no loss of generality in assuming that $s\le -a+1$ on $\mc{E}_k$.
Let $s_k^++1>-a+1$ denote the $s$ value of the closest ramification
point in $\Sigma_k$ to the $s=-a+1$ circle in $\mc{E}_k$.  The
cylinder $\mc{E}_k$ then extends as $\overline{\mc{E}}_k =
(-\infty,s_k^++1)\times\R/2\pi m\Z$.  The convergence of the sequence
of branched covers $T_{-s_{kj}}(\Sigma_{kj})$ from Step 1 implies that
the sequence $\{s_k^+\}$ converges to a number $s_+$ with $|s_++1|\le
r$.

Expand $\sigma_k$ and $\sigma_{*j}$ on $\overline{\mc{E}_k}$ as in
\eqref{eqn:expandsigma} and
\eqref{eqn:expandsigma*}.
By the convergence in Lemma~\ref{lem:asub}(a) at $s=s_+$, we have
\begin{equation}
\label{eqn:coco4}
\lim_{k\to\infty} \sum_\gamma \exp(2E_\gamma s_k^+)
\left|\sigma_{*j\gamma} -
\theta_{kj}^{-1}\sigma_{k\gamma} \right|^2 = 0.
\end{equation}

Suppose that $\sigma_k$ has winding number $\eta$ on
$\overline{\mc{E}_k}$ for all $k$.  Let $\gamma_+$ and $\gamma_-$ be
orthonormal eigenfunctions of $L_m$ with winding number $\eta$ and
eigenvalues $E_+\ge E_-$.  By Lemma~\ref{lem:AK}, the following
hold for each $k$:
\begin{itemize}
\item
At least one of the coefficients $\sigma_{k\gamma_-}$,
$\sigma_{k\gamma_+}$ is nonzero.
\item
If $\gamma$ is an eigenfunction of $L_m$ with $E_\gamma < E_-$, then
$\sigma_{k\gamma}=0$.
\end{itemize}
By Lemmas~\ref{lem:asub}(a) and \ref{lem:extend}(a), the function
$\sigma_{*j}$ also has winding number $\eta$ on $\overline{\mc{E}_k}$,
so the above two properties also hold for the coefficients
$\sigma_{*j\gamma}$.

It now follows from \eqref{eqn:coco4} that there is a $k$-independent
number $c$ with
\begin{equation}
\label{eqn:KIN}
\sum_{E_\gamma > E_+} \exp(2E_\gamma s_k^+)|\sigma_{k\gamma}|^2 <
c\left(\exp(2E_- s_k^+)|\sigma_{k\gamma_-}|^2 + \exp(2E_+ s_k^+)
|\sigma_{k\gamma_+}|^2\right).
\end{equation}
Let $\kappa>0$ denote the difference between $E_+$ and the next
largest eigenvalue.  Then it follows from \eqref{eqn:KIN} and elliptic
regularity for the operator $D_\Sigma^*$ that there is
a $k$-independent constant $c$ with
\begin{gather*}
\left|
\sigma_k(s,t) - \sigma_{k\gamma_-}\exp(E_-s)\gamma_-(t) -
 \sigma_{k\gamma_+}\exp(E_+s)\gamma_+(t)
\right|
\quad\quad\quad\quad\quad\quad\quad
\\
\quad\quad\quad\quad
< c \cdot 
\exp(\kappa(s-s_k^+))
\left|
\sigma_{k\gamma_-}\exp(E_-s)\gamma_-(t) +
\sigma_{k\gamma_+}\exp(E_+s)\gamma_+(t)
\right|
\end{gather*}
for all $(s,t)\in\overline{\mc{E}_k}$ with $s\le s_k^+-1$.  Given
$\varepsilon_0>0$, choose $R\ge 2$ sufficiently large that
\begin{equation}
\label{eqn:R'}
c\exp(-\kappa (R-1))<\varepsilon_0.
\end{equation}
Then the conclusions of Proposition~\ref{prop:approx1} follow when
$(\Sigma,\sigma)=(\Sigma_k,\sigma_k)$ and the point $p(z)$ lies in the
external edge $e$ of $\tau(\Sigma)$.

\medskip

{\em Step 5.\/} We now show that for any $\varepsilon_0>0$, there
exists $R$ such that for all $k$, the conclusions of
Proposition~\ref{prop:approx1} hold for
$(\Sigma,\sigma)=(\Sigma_k,\sigma_k)$ whenever $p(z)$ is on an
internal edge of the tree $\tau(\Sigma)$.

To start, we can assume that $R>2r+1$.  This ensures that we only have
to consider $z$ in a compact component $\mc{E}_k$ of
$\Sigma_k\setminus \cup_{j''}U_{kj''}$ corresponding to an internal
edge $e$ of $T$.  Let $j'$ and $j$ denote the lower and upper vertices
respectively of $e$.  Let $s_k^++1 > -a+1$ denote the $s$ value of the
nearest ramification point in $\Sigma_{kj}$ to the $s=-a+1$ circle in
$\mc{E}_k$, and let $s_k^--1<s_k$ denote the $s$ value of the nearest
critical point in $\Sigma_{kj'}$ to the $s=s_k$ circle in $\mc{E}_k$.
Thus $\mc{E}_k$ extends to a cylinder $\overline{\mc{E}_k}\simeq
(s_k^--1,s_k^++1)\times\R/2\pi m(e)\Z$.  As in Step 4, the sequence
$\{s_k^+\}$ converges to a number $s_+$, while $\lim_{k\to\infty}s_k^-
= -\infty$.  Define $E_+$ and $E_-$ as in the proof of
Lemma~\ref{lem:extend}(b), and let $\gamma_\pm$ be a normalized
eigenfunction with eigenvalue $E_\pm$.  We assume in what follows that
if $E_+=E_-$, then the corresponding eigenspace is one dimensional;
the argument in the case when the dimension is two has no substantive
differences.

Similarly to \eqref{eqn:KIN}, there are $k$-independent numbers $c_+$
and $c_-$ such that
\begin{equation}
\label{eqn:c+c-}
\begin{split}
\sum_{E_\gamma > E_+}\exp(2E_\gamma s_k^+)|\sigma_{k\gamma}|^2 & <
c_+\exp(2E_+s_k^+)|\sigma_{k\gamma_+}|^2,\\
\sum_{E_\gamma < E_-} \exp(2E_\gamma s_k^-)|\sigma_{k\gamma}|^2 & <
c_-\exp(2E_-s_k^-)|\sigma_{k\gamma_-}|^2.
\end{split}
\end{equation}
It follows that there are $k$-independent numbers
$c,\kappa>0$ such that
\begin{equation}
\label{eqn:CKI}
\begin{split}
\bigg|\sum_{E_\gamma > E_+}\sigma_{k\gamma}\exp(E_\gamma s)\gamma(t)\bigg| &<
 c|\sigma_{k\gamma_+}| \exp(E_+s) \exp(\kappa(s-s_k^+)),\\
\bigg|\sum_{E_\gamma < E_-}\sigma_{k\gamma} \exp(E_\gamma
s)\gamma(t)\bigg|
 &< c|\sigma_{k\gamma_-}| \exp(E_-s)
\exp(\kappa(s_k^- - s))
\end{split}
\end{equation}
whenever $s_k^- \le s \le s_k^+$.  

Suppose that $\sigma_k$ has winding number $\eta$ on $\mc{E}_k$ for
all $k$.  We now show that the eigenfunction $\gamma_+$ has winding
number $\eta$.  If $\gamma$ is a normalized eigenfunction with
$E_\gamma<E_+$, then since $\sigma_{*j\gamma_+}\neq 0$ and
$\sigma_{*j\gamma}=0$, it follows from \eqref{eqn:coco4} that
\[
\lim_{k\to\infty}
\frac{\sigma_{k\gamma}}{\sigma_{k\gamma_+}} 
=
0.
\]
Combining this limit for $E_-\le E_\gamma < E_+$ with the inequalities
\eqref{eqn:CKI}, we deduce that for any $\varepsilon>0$, if $k$ is
sufficiently large then
\begin{equation}
\label{eqn:sameDegree}
\frac{|\sigma_k(s,t) -
\sigma_{k\gamma_+}\exp(E_+s)\gamma_+(t)|}{|\sigma_{k\gamma_+}|\exp(E_+s)}
<
c\exp(\kappa(s-s_k^+)) + \varepsilon \exp((E_--E_+)s)
\end{equation}
whenever $s_k^- \le s \le s_k^+$.  By taking $s$ sufficiently small
and then taking $\varepsilon$ sufficiently small (both of which
we can do by taking $k$ sufficiently large), we can make the right hand
side of \eqref{eqn:sameDegree} less than $\min_t|\gamma_+(t)|$.  Hence
there exist $k$ and $s$ such that that $\gamma_+$ has the same winding
number as $\sigma_k(s,\cdot)$, and of course the latter winding number
is $\eta$.

Likewise, $\gamma_-$ has winding number $\eta$.  In particular, there
are no eigenvalues between $E_-$ and $E_+$.  There are now two cases
to consider regarding $E_-$ and $E_+$.

Suppose first that $E_+>E_-$.  Then the inequalities \eqref{eqn:c+c-}
imply that
\begin{gather}
\label{eqn:approx1c'}
\left|\sigma_k(s,t) - \sigma_{k\gamma_+}\exp(E_+s)\gamma_+(t) -
\sigma_{k\gamma_-}\exp(E_-s)\gamma_-(t)
\right|
<
 \\
\nonumber
c\left|
\exp(\kappa(s-s_k^+))\sigma_{k\gamma_+}
\exp(E_+s)\gamma_+(t) + 
\exp(\kappa(s_k^--s))\sigma_{k\gamma_-}\exp(E_-s)\gamma_-(t)\right|
\end{gather}
whenever $s_k^- +1 \le s \le s_k^+ - 1$.  Given $\varepsilon_0>0$, choose
$R\ge 2$ sufficiently large that \eqref{eqn:R'} holds.  Then
\eqref{eqn:approx1c'} implies the conclusions of
Proposition~\ref{prop:approx1} when
$(\Sigma,\sigma)=(\Sigma_k,\sigma_k)$ and $p(z)$ is in the edge of
$\tau(\Sigma)$ corresponding to $e$.

Suppose next that $E_+=E_-$.  Recall that we are assuming that the
corresponding eigenspace is one dimensional, so that
$\gamma_+=\gamma_-$.  Then \eqref{eqn:approx1c'} holds with
the $\gamma_-$ term on the left hand side deleted.  So given
$\varepsilon_0>0$, it is enough choose $R\ge 2$ sufficiently large that
$
c\exp(-\kappa(R-1)) <
\varepsilon_0/2.
$

This completes the proof of Proposition~\ref{prop:approx1}.
\qed

\medskip

{\em Step 6.\/}
We now prove an addendum to Lemma~\ref{lem:extend}.

\begin{lemma}
\label{lem:extA}
For each $j$, and for each internal edge $e$ of $T$ incident to the
$j^{th}$ internal vertex, the following two points hold:
\begin{description}
\item{(i)}
$\sigma_{*j}$ is nonvanishing on the component of $\Sigma_{*j}\setminus
V_{*j}$ corresponding to $e$.
\item{(ii)}
The sequence $\{\Psi_{kj}^{0,1}(\sigma_k) -
\theta_{kj}\sigma_{*j}\}_{k=a,a+1,\ldots}$ converges to zero in the
$C^\infty$ topology on the end in $\Sigma_{*j}$ that corresponds to $e$.
\end{description}
\end{lemma}

\begin{proof}
Let $m$ denote the multiplicity of $e$ and write
$\widetilde{S^1}\eqdef\R/2\pi m\Z$.  Without loss of generality, the
component $\mc{E}_k$ of $\Sigma_k\setminus \bigcup_{j'}U_{kj'}$
corresponding to $e$ is identified with
$[s_k,-a+1]\times\widetilde{S^1}$.  On the corresponding end of
$\Sigma_{*j}$ where $-a\le s \le -a+1$, expand $\sigma_{*j}$ as in
\eqref{eqn:expandsigma*}.  Recall from Step 3 that there is a smallest
eigenvalue $E_+$ of $L_m$ such that $E_+=E_{\gamma_+}$ with
$\sigma_{*j\gamma_+}\neq 0$; and in particular the expansion
\eqref{eqn:expandsigma*} is valid for all $s\le -a+1$.  It follows
that if $-s$ is large, then the winding number of
$\sigma_{*j}(s,\cdot)$ around $\widetilde{S^1}$ equals the winding
number of $\gamma_+$.  By Step 5 and Lemma~\ref{lem:asub}(a), the
latter is the winding number of $\sigma_{*j}(s,\cdot)$ when $s>-a$.
Since all zeroes of $\sigma_{*j}$ have negative degree, we conclude
that $\sigma_{*j}$ is nonvanishing on $(-\infty,-a] \times
\widetilde{S^1}$.  This proves (i).

To prove (ii), it is enough to show that given $s\le -a$,
\begin{equation}
\label{eqn:ii}
\lim_{k\to\infty}\sum_\gamma \exp(2E_\gamma s)|\theta_{kj}\sigma_{*j\gamma}
-\sigma_{k\gamma}|^2 =0.
\end{equation}
Here we have expanded $\sigma_k$ on $\mc{E}_k$ as in \eqref{eqn:expandsigma}.
Since $\theta_{kj}\le 1$, it follows from the convergence in
Lemma~\ref{lem:asub}(a) that (ii) holds when $s=-a+1$.
It is then enough to show that given $s\le -a$,
\begin{equation}
\label{eqn:iiETS}
\lim_{k\to\infty}\sum_{E_\gamma < E_+} \exp(2E_\gamma
s)|\sigma_{k\gamma}|^2=0.
\end{equation}
By \eqref{eqn:c+c-}, there is 
a $k$-independent constant $c$ such that
\begin{equation}
\label{eqn:TLI}
\sum_{E_\gamma < E_+} \exp(2E_\gamma
s)|\sigma_{k\gamma}|^2 < c\exp(2E_-s)|\sigma_{k\gamma_-}|^2
\end{equation}
whenever $s_k\le s \le -a+1$.  For any given $s\le -a$, if $k$ is
sufficiently large then $s_k\le s$ so that \eqref{eqn:TLI} is
applicable.  The inequality \eqref{eqn:TLI} then implies
\eqref{eqn:iiETS} because $\lim_{k\to\infty} \sigma_{k\gamma_-}=0$.
\end{proof}

\medskip

{\em Step 7.\/}
Let $\varepsilon_0>0$ be given;  we now show that there
exists $R$ such that the conclusions of
Proposition~\ref{prop:approx2} hold for $(\Sigma,\sigma)=
(\Sigma_k,\sigma_k)$ whenever $k$ is sufficiently large.

In fact, we can take $R\eqdef r+1$, where $r$ was fixed in Step 1.  To
see why, let $Z$ be a cluster of ramification points in $\Sigma_k$
satisfying the assumptions of Proposition~\ref{prop:approx2}.
Then $Z$ contains all the ramification points in $\Sigma_{kj}$ for
some $j$, while our assumption that $k\ge a >2r+1$ implies that $Z$
contains no other ramification points.  Thus
$\Sigma_{kZ}=\widehat{\Sigma}_{kj}$.  By Lemmas~\ref{lem:asub},
\ref{lem:extend}, and \ref{lem:extA}, there is a nonvanishing
$(0,1)$-form $\sigma_{*j}\in\widetilde{\Coker}(D_{\Sigma_{*j}})$ such
that if $k$ is sufficiently large, then
\[
\left|\sigma_k - (\Psi_{kj}^{0,1})^{-1}(\theta_{kj}\sigma_{*j})\right| <
\frac{\varepsilon_0}{2}|\sigma_k|
\]
at all points in $\Sigma_k$ within distance $1/\varepsilon_0$ of a
ramification point in $Z$.  By the conditions on $\Psi_{kj}$, and
using the vector bundle structure on \eqref{eqn:vvb}, if $k$
is sufficiently large then we can also find a nonvanishing $(0,1)$-form
$\sigma_{kZ}\in\widetilde{\Coker}(D_{\Sigma_{kZ}})$ such that
\[
\left|
(\Psi_{kj}^{0,1})^{-1}(\theta_{kj}\sigma_{*j})
- \sigma_{kZ}
\right|
<
\frac{\varepsilon_0}{2}|\sigma_k|
\]
at all points in $\Sigma_k$ within distance $1/\varepsilon_0$ of a
ramification point in $Z$.  Combining the above two inequalities shows
that the conclusions of Proposition~\ref{prop:approx2} hold for
$(\Sigma,\sigma)= (\Sigma_k,\sigma_k)$ whenever $k$ is sufficiently
large.

This completes the proof of Proposition~\ref{prop:approx2}.
\qed

\subsection{Proof of the relative size estimate}
\label{sec:RSProof}

We now prove Proposition~\ref{prop:RS}.  The proof has four steps.

\medskip

{\em Step 1.\/} We begin by using Proposition~\ref{prop:approx1} to
derive an estimate for the change in $|\sigma|$ along a cylinder away
from the ramification points.
Let $\Sigma\in\mc{M}$, let $\sigma\in\Coker(D_\Sigma)$ be
nonvanishing, and let $e$ be an edge of the tree $\tau(\Sigma)$ of
multiplicity $m$.  Let $\mc{E}$ denote the cylinder in $\Sigma$
corresponding to $e$, and identify $\mc{E}$ with an interval cross
$\R/2\pi m\Z$ as usual.  Then on $\mc{E}$ we can write
\begin{equation}
\label{eqn:writePiW}
\Pi_W\sigma(s,t) = \exp(E_-(\sigma,e)s)\gamma_-(t) +
\exp(E_+(\sigma,e)s)\gamma_+(t),
\end{equation}
where $\gamma_\pm$ are orthogonal eigenfunctions of $L_m$ with
eigenvalues $E_{\pm}(\sigma,e)$, and at least one of $\gamma_\pm$ is
nonzero.  It follows from the above equation that
\[
\log\|\Pi_W\sigma(s,\cdot)\| = \max\big\{E_-(\sigma,e)s +
\log\|\gamma_-\|, E_+(\sigma,e)s + \log\|\gamma_+\|\big\} + \text{Error}
\]
where $\|\cdot\|$ denotes the $L^2$ norm on $\R/2\pi m\Z$, and $0 \le
\text{Error} \le \frac{\log 2}{2}$.
Consequently, if $s'<s$ then
\begin{equation}
\label{eqn:logPiW1}
\log\frac{\|\Pi_W\sigma(s,\cdot)\|}{\|\Pi_W\sigma(s',\cdot)\|} \in
\bigg[
E_-(\sigma,e)(s-s') - \log 2,
E_+(\sigma,e)(s-s') + \log 2
\bigg].
\end{equation}

Next, observe that there are only finitely many possible values of the
winding number $\eta(\sigma,e)$ when $\Sigma\in\mc{M}$,
$\sigma\in\Coker(D_\Sigma)$ is nonvanishing, and $e$ is an edge of
$\tau(\Sigma)$.  This follows from the winding bounds
\eqref{eqn:windbound} together with equations \eqref{eqn:ACZ} and
\eqref{eqn:formZeroes}.   Hence there is a constant $\varepsilon_0>0$
such that for any $\Sigma\in\mc{M}$ and nonvanishing
$\sigma\in\Coker(D_\Sigma)$, on the cylinder $\mc{E}$ in $\Sigma$ corresponding
to an edge $e$ of $\tau(\Sigma)$, the inequality \eqref{eqn:approx1}
implies that
\begin{equation}
\label{eqn:logPiW2}
\big|\log|\sigma(s,t)| - \log\|\Pi_W\sigma(s,\cdot)\|\big| <
1/\varepsilon_0.
\end{equation}
Let $R$ denote the constant provided by Proposition~\ref{prop:approx1} for
this $\varepsilon_0$.  Then by the inequalities \eqref{eqn:logPiW1} and
\eqref{eqn:logPiW2}, we conclude that there is a constant $c>0$ such
that on a cylinder $\mc{E}$ corresponding to an edge $e$, if $s'<s$
have distance at least $R$ from the endpoints of the corresponding
interval, then
\begin{equation}
\label{eqn:logPiW}
\log|\sigma(s,\cdot)| -
\log|\sigma(s',\cdot)| \in
\big[
E_-(\sigma,e)(s-s') - c,
E_+(\sigma,e)(s-s') + c.
\big]
\end{equation}

{\em Step 2.\/} We now inductively define certain constants
$d_k,r_k,r_k'$ for $k=1,\ldots,N$, whose significance will become
clear in subsequent steps.

To start, define $d_1\eqdef 0$.

Next, supposing that $d_k$ has been defined, we want to choose $r_k>R$
and $r_k'$ with the following property: Let $\Sigma\in\mc{M}$ and let
$\sigma\in\Coker(D_\Sigma)$ be nonvanishing.  Let
$B\subset\tau(\Sigma)$ be a compact connected set such that:
\begin{description}
\item{(i)}
The vertices in $B$ correspond to an $r_k$-isolated cluster of
ramification points with diameter $\le d_k$ and total ramification index $k$.
\item{(ii)}
Each boundary point in $B$ has distance exactly $R$ from the nearest
vertex in $B$.
\end{description}
Let $z_1,z_2\in \Sigma$ with $p(z_1), p(z_2)\in B$.  Then
\begin{equation}
\label{eqn:BE}
\big|\log |\sigma(z_1)| - \log|\sigma(z_2)|\big| \le r_k'.
\end{equation}
The existence of such $r_k$ and $r_k'$ follows by applying
Proposition~\ref{prop:approx2} with $r=d_k/2$ and
$\varepsilon_0<1/\max(R,d_k)$, and using compactness of the
projectivization of the vector bundle $\mc{V}(\cdots\mid\cdots)$ in
\eqref{eqn:vvb} over $\mc{M}_r(\cdots\mid\cdots)$.  This last
compactness follows from Lemmas~\ref{lem:mrc} and \ref{lem:LFI}.

Finally, let $k>1$, and suppose that $d_{i},r_{i},r_{i}'$ have been
defined for all $i<k$.  Then 
\begin{equation}
\label{eqn:dk}
d_k \eqdef \max_{i+j=k,\;\; 0<i,j<k}(d_i+d_j+\max\{r_i,r_j\}).
\end{equation}

{\em Step 3.\/} We claim now that for any $\Sigma\in\mc{M}$, the
internal vertices in the tree $\tau(\Sigma)$ can be partitioned into disjoint
subsets $V_1,\ldots,V_l$ such that for each $i=1,\ldots,l$, the
following two properties hold.  Let $k_i$ denote the total
ramification index of the ramification points in $\Sigma$
corresponding to vertices in $V_i$.
\begin{description} \item{(a)} The set $V_i$ has diameter at most
$d_{k_i}$ in $\tau(\Sigma)$, and is contained in a connected set $B$
which does not intersect any $V_j$ with $i\neq j$.
\item{(b)}
Let $e$ be an edge incident to vertices in $V_i$ and $V_j$ with $i\neq
j$.  Then the length of $e$ is greater than $\max\{r_{k_i},r_{k_j}\}$.
\end{description}
We construct a partition satisfying (a) and (b) by induction as
follows.  Start with the partition into sets of cardinality one.  Then
(a) automatically holds.  For the induction step, suppose we have a
partition satisfying (a) but not (b).  Then there exists an edge $e$
incident to vertices in $V_i$ and $V_j$ with $i\neq j$ whose
length is at most $ \max\{r_{k_i}, r_{k_j}\}$.  We now modify
the partition by merging the subsets $V_i$ and $V_j$ into a single
subset.  Condition (a) still holds because
\[
\op{diam}(V_i\cup V_j) = \op{diam}(V_i \cup V_j \cup e) \le d_{k_i} +
d_{k_j} + \max\{r_{k_i}, r_{k_j}\}  \le d_{k_i+k_j} 
\]
by \eqref{eqn:dk}.  Since there are only finitely many vertices,
repeating this step must eventually yield a partition satisfying both
(a) and (b).

{\em Step 4.\/} We now complete the proof of
Proposition~\ref{prop:RS}.  By conditions (a) and (b) in Step 3, we
can find compact connected subsets $B_1,\ldots,B_l$ of $\tau(\Sigma)$,
together containing all of the internal vertices, such that each $B_i$
satisfies conditions (i) and (ii) in Step 2 with $B=B_i$ and $k=k_i$.
Then to prove the estimate \eqref{eqn:RS}, divide the path $P_{x,y}$
into segments, each of which is either contained in one of the $B_i$'s
or outside the interiors of all of the $B_i$'s.  Use
\eqref{eqn:BE} to estimate the change in $\log|\sigma|$ along segments
of the former type, and use
\eqref{eqn:logPiW} to estimate the change in
$\log|\sigma|$ along segments of the latter type.
\qed

\subsection{How moving a ramification point affects the cokernel}
\label{sec:RR}

This subsection proves Proposition~\ref{prop:RBP0} below, which
describes how the cokernel of $D_\Sigma$ changes as one modifies
$\Sigma$ by moving a ramification point.  Lemma~\ref{lem:RBP} is a
special case of Proposition~\ref{prop:RBP0}.

To state Proposition~\ref{prop:RBP0}, assume that $S(t)=\theta$, so
that the operator $D_\Sigma$ is $\C$-linear for each
$\Sigma\in\mc{M}$.  Fix integers $\eta_i^+$ for $i=1,\ldots,N_+$ and
$\eta_j^-$ for $j=-1,\ldots,-N_-$ satisfying
\eqref{eqn:nvw}.  Assume also that
\begin{equation}
\label{eqn:decayforce}
\eta_i^+ \ge \ceil{a_i\theta}, \quad \quad \eta_j^- \le
\floor{a_j\theta}.
\end{equation}
Then the vector space $V\eqdef V(\eta_1^+,\ldots,\eta_{N_+}^+ \mid
\eta_{-1}^-,\ldots,\eta_{-N_-}^-)$ from \eqref{eqn:dimRV} is a
complex linear subspace of $\Coker(D_\Sigma)$ with complex dimension $1$.

Given $\Sigma\in\mc{M}$ and an edge $e$ of the tree $\tau(\Sigma)$,
let $\mc{E}$ denote the cylinder in $\Sigma$ corresponding to $e$.
Given $0\neq\sigma\in V$, let $\eta(e)\eqdef\eta(\sigma,e)$ denote the
winding number of $\sigma$ around $\mc{E}$; by equation
\eqref{eqn:formZeroes}, this depends only on the numbers $\eta_i^+$
and $\eta_j^-$.  Also recall the notation $W(e)$ from
Definition~\ref{def:PiW}; here we have $\dim_\C W(e)=1$.  If we choose
an identification of $\mc{E}$ with an interval cross $\R/2\pi m(e)\Z$
commuting with the projections to $\R\times S^1$, then as in
\eqref{eqn:writePiW}, on $\mc{E}$ we can write
\begin{equation}
\label{eqn:gammae1}
\Pi_W\sigma(s,t) =
\exp\left(\left(\theta-\frac{\eta(e)}{m(e)}\right)s\right)
\sigma_e(t)
\end{equation}
where $\sigma_e\in W(e)$ is given by
\begin{equation}
\label{eqn:gammae2}
\sigma_e(t) = a_e \exp\left(\frac{\eta(e)}{m(e)} it \right)
\end{equation}
for some $a_e\in\C^\times$.

Since $\dim_\C V=1$, the eigenfunction $\sigma_e$ determines
$\sigma$, which in turn determines $\sigma_{e'}$ for any other edge
$e'$ of $\tau(\Sigma)$.  Thus for every pair of edges $e,e'$, the map
sending $\sigma_e$ to $\sigma_{e'}$ is an isomorphism
\begin{equation}
\label{eqn:Phiee}
\Phi_{e,e'}(\Sigma)\in\Hom(W(e),W(e')) = W(e)^* \tensor W(e')
\end{equation}
which depends only on the branched cover $\Sigma\in\mc{M}$ and on our fixed
integers $\eta_i^+$ and $\eta_j^-$.  We now want to study how
$\Phi_{e,e'}$ changes as we rotate the ramification points in the $t$
direction.

If $v$ is an internal vertex, define a rational number $r(e,e';v)$ as
follows: Let $f$ and $f'$ denote the edges incident to $v$ that lead
from $v$ to $e$ and $e'$ respectively, and define
\[
r(e,e';v) \eqdef \frac{\eta(f)}{m(f)} - \frac{\eta(f')}{m(f')}.
\]

\begin{proposition}
\label{prop:RBP0}
For all $\varepsilon>0$ there exists $r>0$ such that the following
holds.  Fix integers $\eta_i^+$ and $\eta_j^-$ satisfying
\eqref{eqn:nvw} and \eqref{eqn:decayforce}.  Fix $\Sigma\in\mc{M}$
such that $\tau(\Sigma)$ is trivalent and each edge of $\tau(\Sigma)$
has length $\ge r$.  Let $v$ be an internal vertex of
$\tau(\Sigma)$, and let $\Sigma'$ be obtained from $\Sigma$
by rotating the ramification point corresponding to $v$ by angle
$\varphi\in\R$ in the $t$ direction. Then
\begin{equation}
\label{eqn:PP'}
\Phi_{e,e'}(\Sigma') = (1+O(\varepsilon))\exp(i \varphi
r(e,e';v))
\Phi_{e,e'}(\Sigma).
\end{equation}
\end{proposition}

Here and below, `$O(\varepsilon)$' denotes a complex number $z$ with
$|z|<\varepsilon$.

\begin{proof}
It follows from the definitions that $\Phi_{e,e''}(\Sigma) =
\Phi_{e',e''}(\Sigma)\circ\Phi_{e,e'}(\Sigma)$ and
$r(e,e'';v)=r(e,e';v)+r(e',e'';v)$.  Hence by induction, it suffices
to prove the lemma when $e$ and $e'$ are both incident to the same
vertex $w$.  We do so in three steps.

{\em Step 1.\/} As in Definition~\ref{def:cluster}, let
$\widehat{\Sigma}\eqdef\Sigma_{\{w\}}$ denote the thrice-punctured
sphere obtained by attaching cylindrical ends to a neighborhood in
$\Sigma$ of the component of the constant $s$ locus corresponding to
$w$.  Let $\widehat{V}$ denote the space of
$\widehat{\sigma}\in\widetilde{\Coker}(D_{\widehat{\Sigma}})$ such
that if $\widehat{\sigma}\neq 0$, then for each edge $e$ of $\tau(\Sigma)$
incident to $w$, $\widehat{\sigma}$ has winding number $\eta(e)$ around the
corresponding cylinder in $\widehat{\Sigma}$.  By equation
\eqref{eqn:formZeroes}, the winding numbers in the definition of $\widehat{V}$
satisfy the appropriate version of equation \eqref{eqn:nvw}, so that
$\dim_\C \widehat{V}=1$.  Thus there is a well-defined element
\[
\Phi_{e,e'}(\widehat{\Sigma}) \in W(e)^* \tensor W(e')
\]
as in \eqref{eqn:Phiee}.  Define $\widehat{\Sigma}'$, $\widehat{V}'$,
and $\Phi_{e,e'}(\widehat{\Sigma}')$ analogously from $\Sigma'$.
Propositions~\ref{prop:approx1} and
\ref{prop:approx2} imply that for any $\varepsilon>0$, if $r$ is
sufficiently large then
\begin{equation}
\begin{split}
\label{eqn:zetazetahat}
\Phi_{e,e'}(\Sigma) &= (1+O(\varepsilon))\Phi_{e,e'}(\widehat{\Sigma}),\\
\Phi_{e,e'}(\Sigma') &= (1+O(\varepsilon))
\Phi_{e,e'}(\widehat{\Sigma}').
\end{split}
\end{equation}

{\em Step 2.\/} We now prove the lemma when $v\neq w$.  Here it
follows from the definition that $r(e,e';v)=0$.  On the other hand,
$\Phi_{e,e'}(\widehat{\Sigma})=\Phi_{e,e'}(\widehat{\Sigma}')$, because in
passing from $\Sigma$ to $\widehat{\Sigma}$ or from ${\Sigma}'$ to
$\widehat{\Sigma}'$, the location of the ramification point
corresponding to $w$ is forgotten.  Thus
\eqref{eqn:PP'} follows from \eqref{eqn:zetazetahat}.

{\em Step 3.\/} We now prove the lemma when $v=w$.  Here
\begin{equation}
\label{eqn:reev}
r(e,e';v) = \frac{\eta(e)}{m(e)} - \frac{\eta(e')}{m(e')}.
\end{equation}
Observe that there is an isomorphism
$\widehat{\Sigma}\to\widehat{\Sigma}'$ covering the automorphism of
$\R\times S^1$ that sends $(s,t)\mapsto (s,t+\varphi)$.  Given an
element $\widehat{\sigma}\in \widehat{V}$, we can push it forward via
this isomorphism to obtain an element
$\widehat{\sigma}'\in\widehat{V}'$.  It follows from
\eqref{eqn:gammae1} and \eqref{eqn:gammae2} that
\[
\widehat{\sigma}'_e =
\exp\left(- i\varphi\frac{\eta(e)}{m(e)}t\right)\widehat{\sigma}_e,
\]
and likewise for $e'$.  Therefore
\begin{equation}
\label{eqn:phieehat}
\Phi_{e,e'}(\widehat{\Sigma}') = \exp\left(i\varphi
\left(\frac{\eta(e)}{m(e)} -
\frac{\eta(e')}{m(e')}\right) t \right) \Phi_{e,e'}(\widehat{\Sigma}).
\end{equation}
We are now done by \eqref{eqn:zetazetahat},
\eqref{eqn:reev}, and \eqref{eqn:phieehat}.
\end{proof}

\section{Application to embedded contact homology}
\label{sec:ech}

As in \S\ref{sec:GP}, let $Y$ be a closed oriented $3$-manifold with a
contact form $\lambda$ whose Reeb orbits are nondegenerate, and let
$J$ be an admissible almost complex structure on $\R\times Y$.  Out of
these data one can define the embedded contact homology (ECH), which
is the homology of a chain complex whose differential $\partial$
counts certain (mostly) embedded $J$-holomorphic curves in $\R\times
Y$.  The significance of ECH is that as explained in \cite[\S1.1]{t3},
it is conjecturally isomorphic to versions of the Ozsv\'{a}th-Szab\'{o} and
Seiberg-Witten Floer homologies defined in \cite{ozsz,krmr}.  However,
most of the foundations of ECH have not yet been established.

In this section, we apply the gluing formula of Theorem~\ref{thm:main}
in a special (but nontrivial) case to prove that the ECH differential
$\partial$ satisfies $\partial^2=0$.  Essentially the same argument
shows that the differential in the periodic Floer homology of mapping
tori \cite{pfh3} also has square zero.


After some combinatorial preliminaries in
\S\ref{sec:pinpout}, the definition of the ECH differential $\partial$
is reviewed in \S\ref{sec:ECHReview}.  The proof that $\partial^2=0$
is given in \S\ref{sec:d2}, using a gluing coefficient calculation
which is carried out in \S\ref{sec:GCC1} and \S\ref{sec:GCC2}.  This
section uses only \S\ref{sec:SGT} (if one accepts the statement of
Theorem~\ref{thm:main}), and is not used elsewhere in the paper.

\subsection{Incoming and outgoing partitions}
\label{sec:pinpout}

To prove that $\partial^2=0$, we need to know the
multiplicities of the ends of the curves that are counted by
$\partial$.  The description of these multiplicities in
\S\ref{sec:ECHReview} requires the following preliminary combinatorial
definitions.

Fix an irrational number $\theta$.  For each nonnegative integer $M$,
we now define two distinguished partitions of $M$, called the
``incoming partition'' and the ``outgoing partition'', and denoted
here by $\pin_\theta(M)$ and $\pout_\theta(M)$ respectively.

\begin{definition}
\label{def:pin}
\cite[\S4]{pfh2}
Define the {\em incoming partition\/} $\pin_\theta(M)$ as follows.
Let $S_\theta$ denote the set of all positive integers $a$ such that
\[
\frac{\ceil{a\theta}}{a} < \frac{\ceil{a'\theta}}{a'}, \quad\quad
\forall a'\in\{1,\ldots,a-1\}.
\]
Let $a\eqdef\max (S_\theta \cap \{1,\ldots,M\})$.  Define
$\pin_\theta(0)=\emptyset$, and inductively define\footnote{ If
$P=(a_1,\ldots,a_k)$ is a partition of $M$ and $Q=(b_1,\ldots,b_l)$ is
a partition of $N$, define a partition $P\cup Q$ of $M+N$ by
\[
P\cup Q \eqdef (a_1,\ldots,a_k,b_1,\ldots,b_l).
\]}
\[
\pin_\theta(M) \eqdef (a) \cup \pin_\theta(M-a).
\]
Define the {\em outgoing partition\/}
\[
\pout_\theta(M) \eqdef \pin_{-\theta}(M).
\]
\end{definition}

We will make frequent use of the following alternate description of
the incoming and outgoing partitions. 

\begin{definition}
Let $\lin_\theta(M)$ denote the lowest convex polygonal path in the
plane that starts at $(0,0)$, ends at $(M,\ceil{M\theta})$, stays
above the line $y=\theta x$, and has corners at lattice points.  That
is, the boundary of the convex hull of the set of lattice points
$(x,y)\in\Z^2$ such that $0\le x \le M$ and $y\ge \theta x$ consists
of the ray $(x=0,y\ge 0)$, the path $\lin_\theta(M)$, and the ray
$(x=M,y\ge\ceil{M\theta})$.
\end{definition}

\begin{lemma}
\label{lem:interpretPartition}
The integers in the incoming partition $\pin_\theta(M)$ are the
horizontal displacements of the segments of the path
$\lin_\theta(M)$ between lattice points.
\end{lemma}

\begin{proof}
Let $(x_1,y_1)\in \Z^2$ denote the first lattice point on the path
$\lin_\theta(M)$, after the initial endpoint $(0,0)$.  Since
$\lin_\theta(M)$ is convex, there are no lattice points in the open
region bounded by the lines $y=\theta x$, $y=(y_1/x_1)x$, and $x=M$.
Hence
\begin{equation}
\label{eqn:firstEdge}
\frac{\ceil{x_1\theta}}{x_1} = \frac{y_1}{x_1} \le
\frac{\ceil{x\theta}}{x}, \quad \forall x=1,\ldots, M.
\end{equation}
Also, since the vector $(x_1,y_1)\in \Z^2$ is indivisible, equality can
hold in \eqref{eqn:firstEdge} only when $x\ge x_1$.  It follows that
$x_1=\max(S_\theta\cap\{1,\ldots,M\})$.  Moreover, the rest of the path
$\lin_\theta(M)$ is the translation of the path
$\lin_\theta(M-x_1)$ by $(x_1,y_1)$.  The lemma follows by
induction.
\end{proof}

\begin{figure}
\begin{center}
\scalebox{0.5}
{\includegraphics{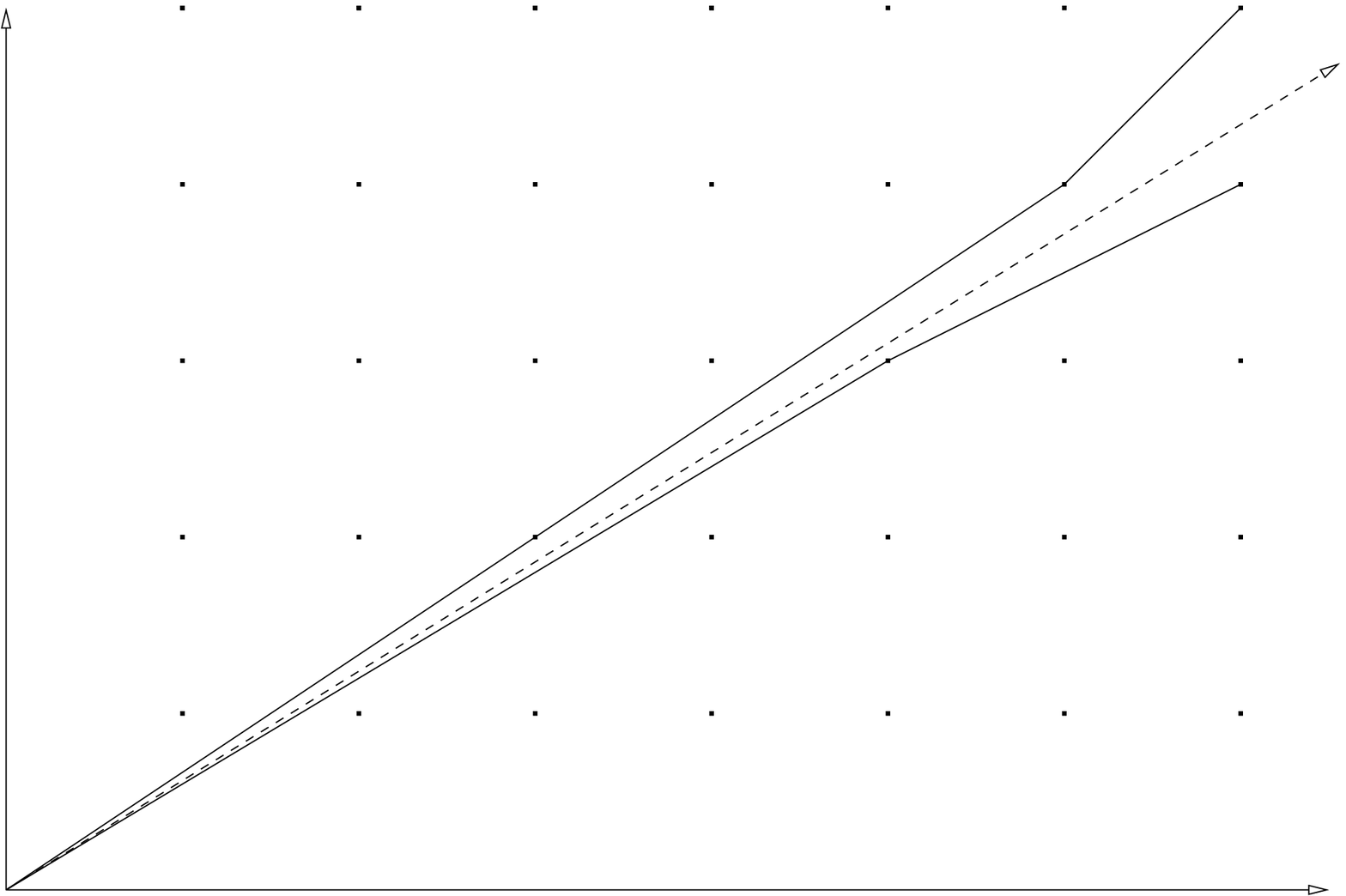}}
\end{center}
\caption{Here $3/5 < \theta < 2/3$, the dashed line is $y=\theta x$, the upper
path is $\lin_\theta(7)$, and the lower path is $\lout_\theta(7)$.
In particular, $\pin_\theta(7)=(3,3,1)$ and
$\pout_\theta(7)=(5,2)$.  }
\label{fig:pinpout}
\end{figure}

Likewise, let $\lout_\theta(M)$ denote the highest concave polygonal
path in the plane which starts at $(0,0)$, ends at
$(M,\floor{M\theta})$, stays below the line $y=\theta x$, and has
corners at lattice points.  Then by
Lemma~\ref{lem:interpretPartition}, the integers in the outgoing
partition $\pout_\theta(M)$ are the horizontal displacements of the
segments of the path $\lout_\theta(M)$ between lattice points.  An
example is shown in Figure~\ref{fig:pinpout}.

The following basic facts about the incoming and outgoing partitions
will be needed later.

\begin{lemma}
\label{lem:pkappa}
\begin{description}
\item{(a)}
If $\pin_\theta(M)=(a_1,\ldots,a_k)$, then
$\sum_{i=1}^k\ceil{a_i\theta} = \ceil{M\theta}$.
\item{(b)}
  If
$\pout_\theta(M)=(b_1,\ldots,b_l)$ then $\sum_{j=1}^l
\floor{b_j\theta} = \floor{M\theta}$.
\end{description}
\end{lemma}

\begin{proof}
This is an immediate consequence of the descriptions of $\pin_\theta(M)$ and
$\pout_\theta(M)$ in terms of the paths $\lin_\theta(M)$ and
$\lout_\theta(M)$.
\end{proof}

\begin{lemma}
\label{lem:PO}
Under the partial order $\ge_\theta$ on the set of partitions of $M$
(see Definition~\ref{def:newpo}), $\pin_\theta(M)$ is maximal and
$\pout_\theta(M)$ is minimal.
\end{lemma}


\begin{proof}
It is an exercise using either Definition~\ref{def:newpo} or
Lemma~\ref{lem:potd} to show that $P\ge_\theta Q$ if and only if one
can get from $P$ to $Q$ by a sequence of the following operations:
\begin{itemize}
\item
Replace $a_1,a_2$ by $a_1+a_2$ where
$\ceil{a_1\theta}+\ceil{a_2\theta}=\ceil{(a_1+a_2)\theta}$.
\item
Replace $a_1+a_2$ by $a_1,a_2$ where
$\ceil{(a_1+a_2)\theta}=\ceil{a_1\theta}+\ceil{a_2\theta}-1$.
\end{itemize}

Now to prove the lemma, by symmetry it is enough to show that
 $\pin_\theta(M)$ is maximal.
Suppose to the contrary that $Q\ge_\theta \pin_\theta(M)$ and $Q\neq
\pin_\theta(M)$.  Then at least one of the following situations occurs:
\begin{description}
\item{(i)}
$\pin_\theta(M)$ contains $m_1$ and $m_2$ with
$\ceil{(m_1+m_2)\theta}=\ceil{m_1\theta}+\ceil{m_2\theta}-1$.
\item{(ii)}
$\pin_\theta(M)$ contains $m_1+m_2$ where
$\ceil{m_1\theta}+\ceil{m_2\theta}=\ceil{(m_1+m_2)\theta}$.
\end{description}

In case (i), write $\pin_\theta(M)=(m_1,m_2,\ldots,m_k)$.
Then by Lemma~\ref{lem:pkappa}(a),
\[
\ceil{(m_1+m_2)\theta} + \sum_{i=3}^k\ceil{m_i\theta} =
\sum_{i=1}^k\ceil{m_i\theta}-1 = \floor{M\theta}.
\]
But this is impossible, since the left side is greater than $M\theta$
and the right side is smaller than $M\theta$.

In case (ii), since $m_1+m_2\in S_\theta$, by
Lemma~\ref{lem:interpretPartition} the path $\lin_\theta(m_1+m_2)$ is
just a line segment from the origin to the point
$(m_1+m_2,\ceil{(m_1+m_2)\theta})$, and this line segment has no
lattice points in its interior.  By the definition of
$\lin_\theta(m_1+m_2)$, the slope of this line segment must be
strictly less than that of the vectors $(m_1,\ceil{m_1\theta})$ and
$(m_2,\ceil{m_2\theta})$.  This contradicts (ii).
\end{proof}

\begin{definition}
The {\em standard ordering convention\/} for the incoming and outgoing
partitions is to write
$\pin_\theta(M)=(a_1,\ldots,a_{k})$
and
$\pout_\theta(M)=(b_1,\ldots,b_{l})$
where
\[
a_{i} \ge a_{i+1}, \quad \quad b_{j} \ge b_{j+1}.
\]
Note that by Definition~\ref{def:pin}, this is equivalent to
\begin{equation}
\label{eqn:OC}
\frac{\ceil{a_i\theta}}{a_i} \le
\frac{\ceil{a_{i+1}\theta}}{a_{i+1}}, \quad \quad
  \frac{\floor{b_j\theta}}{b_j} \ge \frac{\floor{b_{j+1}\theta}}{b_{j+1}}.
\end{equation}
\end{definition}

\subsection{The ECH differential $\partial$}
\label{sec:ECHReview}

We now briefly review the definition of the differential $\partial$ in
embedded contact homology, in preparation for showing that
$\partial^2=0$.

A $J$-holomorphic curve may have several ends at covers of an embedded
Reeb orbit $\gamma$, with various covering multiplicities.  The ECH
chain complex only keep tracks of the sum of these multiplicities.
For this purpose we make the following definitions.

\begin{definition}
An {\em orbit set\/} is a finite set of pairs
$\alpha=\{(\alpha_i,m_i)\}$ where the $\alpha_i$'s are distinct
embedded Reeb orbits and the $m_i$'s are positive integers%
\footnote{This is
different from the notation in Definition~\ref{def:MSFT}, where
$\alpha$ and $\beta$ are ordered lists of Reeb orbits which might be
multiply covered.}.  Define $[\alpha]\eqdef\sum_im_i[\alpha_i]\in
H_1(Y)$.  The orbit set $\alpha$ is {\em admissible\/} if $m_i=1$
whenever $\alpha_i$ is hyperbolic.  If $\beta=\{(\beta_j,n_j)\}$ is
another orbit set with $[\alpha]=[\beta]$, define
$H_2(Y,\alpha,\beta)$ to be the set of relative homology classes of
2-chains $Z$ in $Y$ with
\[
\partial Z = \sum_im_i\alpha_i - \sum_j n_j\beta_j.
\]
Thus $H_2(Y,\alpha,\beta)$ is an affine space over $H_2(Y)$.
\end{definition}

\begin{definition}
If $\alpha=\{(\alpha_i,m_i)\}$ and $\beta=\{(\beta_j,n_j)\}$ are orbit
sets with $[\alpha]=[\beta]$, let $\mc{M}^J(\alpha,\beta)$ denote the
moduli space of $J$-holomorphic curves $u$ with positive ends at
covers of $\alpha_i$ with total multiplicity $m_i$, negative ends at
covers of $\beta_j$ with total multiplicity $n_j$, and no other ends.
In contrast to Definition~\ref{def:MSFT}, the ends of $u$ are not
ordered or asymptotically marked.  Note that the projection of each
$u\in\mc{M}^J(\alpha,\beta)$ to $Y$ has a well-defined relative
homology class $[u]\in H_2(Y,\alpha,\beta)$.  For $Z\in
H_2(Y,\alpha,\beta)$ we then define
\[
\mc{M}^J(\alpha,\beta,Z) \eqdef \{u\in \mc{M}^J(\alpha,\beta) \mid
[u]=Z\}.
\]
\end{definition}

\begin{definition}
Given a homology class $\Gamma\in H_1(Y)$, the ECH chain complex
$C_*(Y,\lambda;\Gamma)$ is a free $\Z$-module with one generator for
each admissible orbit set $\alpha$ with $[\alpha]=\Gamma$.
\end{definition}

To fix the signs in the differential below,  for each admissible orbit
set we need to choose an ordering of its positive hyperbolic orbits%
\footnote{
Alternately one can define the chain complex to be generated by
admissible orbit sets in which the positive hyperbolic orbits are
ordered, modulo the relation that reordering the positive hyperbolic
orbits in a generator multiplies the generator by the sign of the
reordering permutation.
}.
To simplify the discussion below, let us do this by fixing some
ordering of all the embedded positive hyperbolic Reeb orbits in $Y$.

The relative index on this chain complex is defined as follows.
(This should be contrasted with Definition~\ref{def:ind}.)

\begin{definition}
(cf.\ \cite[\S1]{pfh2})
If $\alpha$ and $\beta$ are orbit sets with $[\alpha]=[\beta]$, and if
$Z\in H_2(Y,\alpha,\beta)$, define the {\em ECH index\/}
\[
I(\alpha,\beta,Z) \eqdef c_1(\xi|_Z,\tau) + Q_\tau(Z) +
\sum_i\sum_{k=1}^{m_i} \op{CZ}_\tau(\alpha_i^k) -
\sum_j\sum_{k=1}^{n_j} \op{CZ}_\tau(\beta_j^k).
\]
Here $\tau$ is a trivialization of $\xi$ over the $\alpha_i$'s and
$\beta_j$'s.  As in \S\ref{sec:GP}, $c_1$ denotes the relative first
Chern class of $\xi$ over a surface representing $Z$, and
$\op{CZ}_\tau(\gamma^k)$ denotes the Conley-Zehnder index of the
$k^{th}$ iterate of $\gamma$.  Also, $Q_\tau$ is the ``relative
self-intersection pairing'' defined in \cite[\S2]{pfh2}.  If
$u\in\mc{M}^J(\alpha,\beta,Z)$, write $I(u)\eqdef I(\alpha,\beta,Z)$.
\end{definition}

It is shown in \cite{pfh2} that $I$ depends only on the orbit sets
$\alpha$ and $\beta$ and on the relative homology class $Z$.  Also,
$I$ is additive in the following sense: if $\gamma$ is another orbit
set with $[\beta]=[\gamma]$ and if $W\in H_2(Y,\beta,\gamma)$, then
there is a well-defined relative homology class $Z+W\in
H_2(Y,\alpha,\gamma)$, and we have
\[
I(\alpha,\gamma,Z+W) = I(\alpha,\beta,Z) + I(\beta,\gamma,W).
\]
The key nontrivial property of the ECH index $I$ is that it gives an
upper bound on the Fredholm index $\op{ind}$ from
Definition~\ref{def:ind}.  Moreover, curves that realize this upper
bound are highly restricted.  To give the precise statements, we need
the following definitions.

\begin{definition}
(cf.\ \cite[\S4]{pfh2})
If $\gamma$ is an embedded Reeb orbit and $M$ is a positive integer, define two
partitions of $M$, the {\em incoming partition\/} $\pin_\gamma(M)$ and
the {\em outgoing partition\/} $\pout_\gamma(M)$, as follows.
\begin{itemize}
\item
If $\gamma$ is positive hyperbolic, then
\begin{equation}
\label{eqn:PHP}
\pin_\gamma(M) \eqdef \pout_\gamma(M) \eqdef (1,\ldots,1).
\end{equation}
\item
If $\gamma$ is negative hyperbolic, then
\begin{equation}
\label{eqn:NHP}
\pin_\gamma(M) \eqdef \pout_\gamma(M) \eqdef \left\{\begin{array}{cl}
(2,\ldots,2), & \mbox{if $M$ is even},\\
(2,\ldots,2,1), & \mbox{if $M$ is odd}.
\end{array}\right.
\end{equation}
\item
If $\gamma$ is elliptic with monodromy angle $\theta$, then (see
\S\ref{sec:pinpout})
\[
\pin_\gamma(M) \eqdef \pin_\theta(M), \quad \quad \pout_\gamma(M)
\eqdef \pout_\theta(M).
\]
\end{itemize}
The {\em standard ordering convention\/} for $\pin_\gamma(M)$ or
$\pout_\gamma(M)$ is to list the entries in nonincreasing order.
\end{definition}

\begin{notation}
Any $J$-holomorphic curve $u\in\mc{M}^J(\alpha,\beta)$ can be uniquely
written as $u=u_0\cup u_1$, where $u_0$ and $u_1$ are unions of
components of $u$, each component of $u_0$ maps to an $\R$-invariant
cylinder, and no component of $u_1$ does.  Given an embedded Reeb
orbit $\gamma$, let $n_\gamma$ denote the total multiplicity of covers
of $\R\times\gamma$ in $u_0$.  Let $m_\gamma^+$ denote the total
multiplicity of all positive ends of $u$ at covers of $\gamma$, and
let $P_\gamma^+$ denote the partition of $m_\gamma^+ - n_\gamma$
consisting of the multiplicities of the positive ends of $u_1$ at
covers of $\gamma$.  Define $m_\gamma^-$ and $P_\gamma^-$ analogously
for the negative ends.
\end{notation}

\begin{definition}
\label{def:admissible}
$u=u_0\cup u_1\in\mc{M}^J(\alpha,\beta)$ is {\em admissible\/} if:
\begin{description}
\item{(a)}
$u_1$ is embedded and does not intersect $u_0$.
\item{(b)}
For each embedded Reeb orbit $\gamma$, under the standard ordering convention:
\begin{description}
\item{$\bullet$}
$P_\gamma^+$ is an initial segment of
$\pout_\gamma(m_\gamma^+)$.
\item{$\bullet$}
$P_\gamma^-$ is an initial segment of $\pin_\gamma(m_\gamma^-)$.
\end{description}
\end{description}
\end{definition}

We can now state the key index inequality.

\begin{proposition}
\label{prop:II}
Let $u=u_0\cup u_1\in\mc{M}^J(\alpha,\beta)$ and suppose that $u_1$
is not multiply covered.  Then:
\begin{description}
\item{(a)}
$\op{ind}(u_1) \le I(u_1)-2\delta(u_1)$, with equality only if
for each embedded Reeb orbit $\gamma$:
\begin{description}
\item{$\bullet$}
$P_\gamma^+ = \pout_\gamma(m_\gamma^+ - n_\gamma)$.
\item{$\bullet$}
$P_\gamma^- = \pin_\gamma(m_\gamma^- - n_\gamma)$.
\end{description}
\item{(b)}
$I(u_1) \le I(u) - 2 \#(u_0 \cap u_1)$, with equality only if the
following hold for each embedded Reeb orbit $\gamma$, under the
standard ordering convention:
\begin{description}
\item{$\bullet$}
$\pout_\gamma(m_\gamma^+ - n_\gamma)$ is an initial
segment of $\pout_\gamma(m_\gamma^+)$.
\item{$\bullet$}
$\pin_\gamma(m_\gamma^- - n_\gamma)$ is an initial segment of
$\pin_\gamma(m_\gamma^-)$.
\end{description}
\end{description}
\end{proposition}

Here $\delta(u_1)$ is a count of the singularities of $u_1$ with
positive integer weights; in particular $\delta(u_1)\ge 0$, with
equality if and only if $u_1$ is embedded.  Also $\#(u_0\cap u_1)$ is
the algebraic intersection number; by intersection positivity, each
intersection point counts positively.

\begin{proof}
Part (a) is proved in \cite[Eq. (18) and Prop. 6.1]{pfh2}, and part
(b) is proved in \cite[Prop. 7.1]{pfh2}, except for two issues.
First, these results are proved in \cite{pfh2} in a slightly different
setting where $Y$ is a mapping torus and an analytical simplifying
assumption (``local linearity'') is made.  The asymptotic analysis
needed to transfer these results to the present setting is carried out
in \cite{siefring}.  Second, the necessary condition for equality in
part (b) is different from the one given in \cite[Prop. 7.1]{pfh2}.
However these two conditions are equivalent by
Lemma~\ref{lem:tfae}(a),(d) below.
\end{proof}

We can now classify the curves with small ECH index for generic $J$.

\begin{proposition}
\label{prop:LIC}
Suppose that $J$ is generic and $u=u_0\cup
u_1\in\mc{M}^J(\alpha,\beta)$.  Then:
\begin{description}
\item{(a)}
$I(u)\ge 0$.
\item{(b)}
If $I(u)=0$, then $u_1=\emptyset$.
\item{(c)}
If $I(u)=1$, then $u$ is admissible and
$\op{ind}(u_1)=1$.
\item{(d)}
If $I(u)=2$ and $\alpha$ and $\beta$ are admissible, then $u$ is
admissible and $\op{ind}(u_1)=2$.
\end{description}
\end{proposition}

\begin{proof}
(This is based on \cite[Lem. 9.5]{pfh2} with simplifications from
\cite[Cor.\ 11.5]{t3}.)  The image of $u_1$ is the union of $k$
irreducible components $v_1,\ldots,v_k$, covered by $u$ with positive
integer multiplicities $d_1,\ldots,d_k$.  Since $J$ is generic,
$\op{ind}(v_i)\ge 1$ for each $i$.

Let $u_1'$ be the union of $d_i$ translates of $v_i$ for $i=1,\ldots,k$.
Then $\op{ind}(u_1')=\sum_{i=1}^kd_i\op{ind}(v_i)$ by definition, and
$I(u_1')=I(u_1)$ since $u_1'$ and $u_1$ go between the same orbit sets and
have the same relative homology class.  So by
Proposition~\ref{prop:II}(a) applied to $u_1'$ and
Proposition~\ref{prop:II}(b) applied to $u_0\cup u_1'$, we obtain
\begin{equation}
\label{eqn:MCII}
\sum_{i=1}^k d_i \op{ind}(v_i) \le I(u) - 2\delta(u_1') - 2\#(u_0\cap u_1),
\end{equation}
with equality only if condition (b) in Definition~\ref{def:admissible}
holds.  Parts (a)--(c) follow immediately from \eqref{eqn:MCII}.

To prove part (d), note that if $I(u)=2$ then $k>0$, because a union
of $\R$-invariant cylinders has $I=0$.  Furthermore the left hand side
of \eqref{eqn:MCII} must equal $2$, because by \cite[Lem.\ 9.4]{pfh2},
if $\alpha$ and $\beta$ are admissible then $\op{ind}(u)$ and $I(u)$
have the same parity. Now there are three possibilities: (i) $k=2$ and
$d_1=d_2=1$; (ii) $k=1$ and $d_1=1$; (iii) $k=1$ and $d_1=2$.  In
cases (i) and (ii) we are done.

To complete the proof we now rule out case (iii).  In
this case we must have $\op{ind}(v_1)=1$.  However, since $\alpha$ and
$\beta$ are admissible, and since $d_1>1$, all Reeb orbits in $\alpha$
and $\beta$ are elliptic.  Since elliptic orbits have odd
Conley-Zehnder index, it follows from the definition of $\op{ind}$ and
the formula for the Euler characteristic of a surface that
$\op{ind}(v_1)$ is even, a contradiction.
\end{proof}

The differential $\partial$ in ECH counts $I=1$ curves in
$\mc{M}^J(\alpha,\beta)/\R$ where $\alpha$ and $\beta$ are admissible
orbit sets.  Such curves may contain multiple covers of the
$\R$-invariant cylinder $\R\times\gamma$ when $\gamma$ is an elliptic
embedded Reeb orbit.  The differential $\partial$ only keeps track of
the total multiplicity of such coverings for each $\gamma$.  We now
give the precise definition of $\partial$, in notation which will be
convenient for the proof that $\partial^2=0$.

\begin{definition}
Let $\alpha$ and $\beta$ be orbit sets.
Define $\mc{M}_1^J(\alpha,\beta,Z)$ to be the set of curves
$u\in\mc{M}^J(\alpha,\beta,Z)$ such that if $\gamma$ is an elliptic
Reeb orbit, then $u$ does not contain $\R\times\gamma$ or any cover
thereof.
\end{definition}

\begin{notation}
If $\alpha$ and $\beta$ are orbit sets, define a ``product'' orbit set
$\alpha\beta$ by adding the multiplicities of all embedded Reeb orbits
involved.  (The index and differential are not well-behaved with
respect to this ``multiplication''.)  Write $\alpha|\beta$ if $\beta$
is divisible by $\alpha$ in this sense, in which case denote the
quotient by $\beta/\alpha$.  Call an orbit set ``elliptic'' if all of
its Reeb orbits are elliptic.
\end{notation}

\begin{definition}
Given a generic $J$ and a system of coherent orientations, define the
ECH differential
\[
\partial: C_*(Y,\lambda,\Gamma) \longrightarrow
C_{*-1}(Y,\lambda,\Gamma)
\]
as follows.  If $\alpha$ and $\beta$ are admissible orbit sets with
$[\alpha]=[\beta]=\Gamma$, then the coefficient of $\beta$ in
$\partial\alpha$ is
\[
\langle\partial\alpha,\beta\rangle \eqdef \sum_{\substack{Z\in
H_2(Y,\alpha,\beta)\\
I(\alpha,\beta,Z)=1}}
\left(
\sum_{\substack{\mbox{\scriptsize $\gamma$ elliptic orbit set}\\
\gamma|\alpha,\beta}}
\#\frac{\mc{M}_1^J\left(\alpha/\gamma, \beta/\gamma,Z-[\R\times\gamma]
\right)}{\R}
\right) \cdot \beta
\]
Here the symbol `$\#$' indicates the signed count.
\end{definition}

To see why $\langle\partial\alpha,\beta\rangle$ is well-defined, first
note that the set being counted is zero-dimensional, because by
Proposition~\ref{prop:LIC}(c), if $I(\alpha,\beta,Z)=1$ and
\[
u_1\in
\mc{M}_1^J\left(\alpha/\gamma, \beta/\gamma,Z-[\R\times\gamma]
\right)
\]
then $\op{ind}(u_1)=1$.  By Remark~\ref{rem:hypor}, the sign
$\epsilon(u_1)$ is well-defined, because admissibility of $\alpha$ and
$\beta$ ensures that $u_1$ does not have an end at a double cover of a
negative hyperbolic orbit or more than one end at a positive
hyperbolic orbit, and an ordering of all positive hyperbolic orbits
has been chosen.  Finiteness of the count results from the following
compactness lemma.

\begin{lemma}
\label{lem:compactness1}
If $\alpha$ and $\beta$ are (not necessarily admissible) orbit sets
and $J$ is generic, then the set
\[
\bigsqcup_{\substack{Z\in
H_2(Y,\alpha,\beta)\\
I(\alpha,\beta,Z)=1}} \mc{M}_1^J(\alpha,\beta,Z)/\R
\]
is compact (and therefore finite).
\end{lemma}

\begin{proof}
(Cf.\ \cite[\S9]{pfh2}.)  Let $\{u_n\}$ be a sequence of curves in
$\mc{M}_1^J(\alpha,\beta)/\R$ with $I(u_n)=1$.  By Stokes' theorem,
the ``energy'' of $u_n$, namely the integral of $u_n^*d\lambda$ over the
domain of $u_n$, is
\[
\int_{u_n}d\lambda = \int_{\alpha}\lambda - \int_{\beta}\lambda,
\]
which does not depend on $n$.  So by Gromov compactness as in
\cite[Lem.\ 9.8]{pfh2}, we can pass to a subsequence so that $\{u_n\}$
converges in the sense of \cite{behwz} to a (possibly) broken curve
with $\op{ind}=1$.

By Proposition~\ref{prop:LIC}(a) and the additivity of the ECH index,
one level of the broken curve has $I=1$ and all other levels have
$I=0$.  By Proposition~\ref{prop:LIC}(b), Lemma~\ref{lem:BCCIndex},
and the additivity of $\op{ind}$, the $I=0$ levels also have
$\op{ind}=0$.  Then the top level cannot have $I=0$, or else by
Lemmas~\ref{lem:BCCIndex} and \ref{lem:PO} it would have $\op{ind}\ge
1$.  Likewise the bottom level cannot have $I=0$.  Hence there is only
one level.

The limiting curve cannot contain a cover of $\R\times\gamma$ with
$\gamma$ elliptic, because the $u_n$'s contain no such covers, and any
$J$-holomorphic curve in the same moduli space component as a cover of
$\R\times\gamma$ is itself a cover of $\R\times\gamma$ because it has
energy zero.
\end{proof}

\subsection{Proof that $\partial^2=0$}
\label{sec:d2}

\begin{theorem}
\label{thm:d2}
If $J$ is generic, then the ECH differential $\partial$ satisfies
$\partial^2=0$.
\end{theorem}

The proof of Theorem~\ref{thm:d2} follows the standard strategy of
analyzing ends of moduli spaces of $I=2$ curves, and consists of a
compactness argument and a gluing argument.  The following are the
kinds of pairs of curves that we will need to glue.

\begin{definition}
Let $\alpha_+$ and $\alpha_-$ be admissible orbit sets.  An {\em ECH
gluing pair\/} is a pair of curves $u_+\in\mc{M}^J(\alpha_+,\beta)$
and $u_-\in\mc{M}^J(\beta,\alpha_-)$ such that:
\begin{description}
\item{(a)}
$I(u_+)=I(u_-)=1$.
\item{(b)}
For each embedded elliptic Reeb orbit $\gamma$:
\begin{description}
\item{(i)} All covers of $\R\times\gamma$ in $u_+$ and $u_-$ are unbranched.
\item{(ii)} If $u_+$ (resp.\ $u_-$) contains covers of $\R\times\gamma$ with
total multiplicity $n_\gamma^+$ (resp.\ $n_\gamma^-$), then the individual
multiplicities comprise the outgoing partition $\pout_\gamma(n_\gamma^+)$
(resp.\ the incoming partition $\pin_\gamma(n_\gamma^-)$).
\item{(iii)} $u_+$ and $u_-$ do not both contain covers of
$\R\times\gamma$.
\end{description}
\end{description}
\end{definition}

To glue ECH gluing pairs, we will apply Theorem~\ref{thm:main}, for
which purpose we will need the following calculation of gluing
coefficients.  If $P$ is a partition in which the
positive integer $n$ appears $r(n)$ times, define
\[
P! \eqdef \prod_{n=1}^\infty n^{r(n)} \cdot r(n)!.
\]
In particular, if $P$ is the (empty) partition of $0$, then $P!=1$.

\begin{proposition}
\label{prop:PTC}
Given integers $0 \le M_+,M_-\le M$, define
\[
S\eqdef (\pin_\theta(M_+);\pout_\theta(M-M_+) \mid
\pout_\theta(M_-) ; \pin_\theta(M-M_-)).
\]
Suppose that under the standard ordering convention,
\begin{equation}
\label{eqn:star}
\begin{array}{l}
\text{$\pin_\theta(M_+)$ is an
initial segment of $\pin_\theta(M)$, and}\\
\text{$\pout_\theta(M_-)$ is an
initial segment of $\pout_\theta(M)$.}
\end{array}
\end{equation}
Then:
\begin{description}
\item{(a)}
If $M_-, M_+ < M$, then for every $\theta$-decomposition of $S$, see
Definition~\ref{def:TD}, there exists $\nu$ with
$I_\nu=J_\nu=\emptyset$ and $|I_\nu'|=|J_\nu'|=1$.
\item{(b)}
If $M_-=M$ then
$c_\theta(S) = \pout_\theta(M-M_+)!$.
\item{(c)}
If $M_+=M$ then
$c_\theta(S) = \pin_\theta(M-M_-)!$.
\end{description}
\end{proposition}

The proof of this proposition is deferred to \S\ref{sec:GCC1} and
\S\ref{sec:GCC2}.  We can now carry out the compactness part of the
proof that $\partial^2=0$ and see how ECH gluing pairs arise.

\begin{lemma}
\label{lem:compactness2}
Assume that $J$ is generic and let $\alpha_+$ and $\alpha_-$ be
admissible orbit sets.  Let $\{u_n\}$ be a sequence of curves in
$\mc{M}_1^J(\alpha_+,\alpha_-)/\R$ such that $I(u_n)=2$.  Then after
passing to a subsequence, $\{u_n\}$ converges in the sense of
\cite{behwz} either to a curve in
$\mc{M}^J_1(\alpha_+,\alpha_-)/\R$, or to a broken curve
$(u_+,\tau_1,\ldots,\tau_k,u_-)$ for some $k\ge 0$, such that each
$\tau_i$ maps to a union of $\R$-invariant cylinders and $(u_+,u_-)$
is an ECH gluing pair.
\end{lemma}

\begin{proof}
As in the proof of Lemma~\ref{lem:compactness1}, we can pass to a
subsequence so that $\{u_n\}$ converges to a (possibly) broken curve
with $\op{ind}=2$, in which each level has $I\ge 0$, and the ECH
indices of the levels sum to $2$.  The top level must have $I>0$;
otherwise, since $\alpha_+$ is admissible, by
Proposition~\ref{prop:LIC}(b) and Lemma~\ref{lem:PO} it would have
$\op{ind}\ge 2$, contradicting additivity of $\op{ind}$ for the broken
curve.  Likewise the bottom level has $I>0$.

Suppose there are at least two levels.  Then it follows that the
limiting broken curve has the form $(u_+,\tau_1,\ldots,\tau_k,u_-)$
where $I(u_+)=I(u_-)=1$ and $I(\tau_i)=0$ for all $i$.  By
Proposition~\ref{prop:LIC}(a), each $\tau_i$ maps to a union of
$\R$-invariant cylinders.

To complete the proof we must verify condition (b) in the definition
of ECH gluing pair.  Let $\gamma$ be an embedded elliptic Reeb orbit.
By Proposition~\ref{prop:LIC}(c),(d), the $u_n$'s and $u_\pm$ are
admissible.  It then follows from Definition~\ref{def:pin} (cf.\
Lemma~\ref{lem:tfae}(a) below) that the multiplicities of the positive
ends of $u_+$ (resp.\ negative ends of $u_-$) at covers of $\gamma$
must comprise the outgoing (resp.\ incoming) partition of $n_\gamma^+$
(resp.\ $n_\gamma^-$).  Assertion (i) now follows from
Lemma~\ref{lem:PO} and additivity of $\op{ind}$ as before.  Assertion
(ii) then follows from the above description of the multiplicities of the
positive ends of $u_+$ and negative ends of $u_-$. 

To prove assertion (iii), let $\theta$ denote the monodromy angle of
$\gamma$, and let $m_\gamma$ denote the total multiplicity of the negative
ends of $u_+$ at covers of $\gamma$.  We can glue
$\tau_1,\ldots,\tau_k$ to obtain an index zero branched cover
$\pi:\Sigma\to\R\times\gamma$, where each positive end of $\Sigma$ is
paired with a negative end of $u_+$, and each negative end of $\Sigma$
is paired with a positive end of $u_-$.  The multiplicities of the
ends of the components of $\Sigma$ determine a $\theta$-decomposition
of
\[
S \eqdef (\pin_\theta(m_\gamma-n_\gamma^+);\pout_\theta(n_\gamma^+) \mid
\pout_\theta(m_\gamma-n_\gamma^-) ; \pin_\theta(n_\gamma^-)).
\]
If assertion (iii) is false, then Proposition~\ref{prop:PTC}(a)
implies that $\Sigma$ has a cylinder component which is attached to
$\R$-invariant cylinders in $u_+$ and $u_-$.  This contradicts the
fact that the $u_n$'s have no components mapping to $\R\times\gamma$.
\end{proof}

We now apply Theorem~\ref{thm:main} to deduce the gluing lemma that
will be needed in the proof that $\partial^2=0$.  Note that by
Lemma~\ref{lem:PO}, an ECH gluing pair becomes a gluing pair as in
Definition~\ref{def:gluingPair} after orderings and asymptotic
markings of the ends of $u_\pm$ are chosen.

\begin{lemma}
\label{lem:d22}
Assume $J$ is generic and let $(u_+,u_-)$ be an ECH gluing pair.
If orderings and asymptotic markings of the ends of $u_\pm$ are
chosen, then:
\begin{description}
\item{(a)}
If $\beta$ is not admissible then $\#G(u_+,u_-)=0$.
\item{(b)}
If $\beta$ is admissible then
\begin{equation}
\label{eqn:factorials}
\#G(u_+,u_-) = \epsilon(u_+)\epsilon(u_-)
\prod_{\mbox{\scriptsize $\gamma$ elliptic embedded Reeb orbit}}
\pout_\gamma(n_\gamma^+)! 
\pin_\gamma(n_\gamma^-)!.
\end{equation}
\end{description}
\end{lemma}

\begin{proof}
For each embedded Reeb orbit $\gamma$, let $m_\gamma$ denote the total
multiplicity of negative ends of $u_+$ at covers of $\gamma$.  By
Theorem~\ref{thm:main}, it is enough to show:
\begin{description}
\item{(c)} If $\gamma$ is hyperbolic and $m_\gamma=1$ then
$c_\gamma(u_+,u_-)=1$.
\item{(d)} If $\gamma$ is hyperbolic and $m_\gamma > 1$ then
$c_\gamma(u_+,u_-)=0$.
\item{(e)} If $\gamma$ is elliptic then
$c_\gamma(u_+,u_-) = \pout_\gamma(n_\gamma^+)!
\pin_\gamma(n_\gamma^-)!$.
\end{description}

Assertion (c) follows immediately from Definition~\ref{def:HGC}.

To prove (d), suppose $\gamma$ is hyperbolic and $m_\gamma>1$.  Recall
from Proposition~\ref{prop:LIC}(c) that $u_+$ and $u_-$ are
admissible.  If $\gamma$ is positive hyperbolic, this means that all
ends of $u_+$ and $u_-$ at (covers of) $\gamma$ have multiplicity
$1$.  It then follows immediately from Definition~\ref{def:HGC}(b)
that $c_\gamma(u_+,u_-)=0$.  If $\gamma$ is negative hyperbolic, then
admissibility implies that $u_+$ has at least one negative end at a
double cover of $\gamma$.  Then $c_\gamma(u_+,u_-)=0$ by
Definition~\ref{def:HGC}(c).

Assertion (e) follows immediately from
Proposition~\ref{prop:PTC}(b),(c), thanks to condition (b) in the
definition of ECH gluing pair, and the admissibility of $u_+$ and
$u_-$.
\end{proof}

\begin{proof}[Proof of Theorem~\ref{thm:d2}.]
Let $\alpha_+$ and $\alpha_-$ be admissible orbit sets.
We will prove that $\langle\partial^2\alpha_+,\alpha_-\rangle = 0$ in
two steps.

{\em Step 1.\/}  We first show that
\begin{equation}
\label{eqn:d2step1}
\begin{split}
\sum_{\mbox{\scriptsize $\beta$ admissible}}
\sum_{\substack{\gamma_+ | \beta,\alpha_+ \\ \gamma_- | \beta, \alpha_-}}
\sum_{Z_+,Z_-}
& \#\frac{\mc{M}_1^J(\alpha_+/\gamma_+,
\beta/\gamma_+, Z_+-[\R\times\gamma_+])}{\R} \\
& \cdot \#\frac{\mc{M}_1^J(\beta/\gamma_-,
\alpha_-/\gamma_-, Z_--[\R\times\gamma_-])}{\R} = 0.
\end{split}
\end{equation}
Here $\gamma_+$ and $\gamma_-$ are elliptic orbit sets with no common
factor, while $Z_+\in H_2(Y,\alpha_+,\beta)$ and $Z_-\in
H_2(Y,\beta,\alpha_-)$ satisfy $I(\alpha_+,\beta,Z_+)=I(\beta,\alpha_-,Z_-)=1$.

To prove \eqref{eqn:d2step1}, we study the ends of the one-dimensional
manifold
\[
\mc{M}\eqdef \bigsqcup_{\substack{Z\in H_2(Y,\alpha_+,\alpha_-) \\
I(\alpha_+,\alpha_-,Z)=2}} \frac{\mc{M}_1^J(\alpha_+,\alpha_-,Z)}{\R}.
\]
If $(u_+,u_-)$ is an ECH gluing pair, in which
$u_+\in\mc{M}^J(\alpha_+,\beta)$ and $u_-\in\mc{M}^J(\beta,\alpha_-)$
for some orbit set $\beta$, let
$V(u_+,u_-)\subset\mc{M}_1^J(\alpha_+,\alpha_-)/\R$ be an open set
like the open set $U$
in Definition~\ref{def:countG}, but where the curves do not have
asymptotic markings or orderings of the ends.  Define
\[
\overline{\mc{M}} \eqdef \mc{M} \setminus \bigsqcup_{(u_+,u_-)}V(u_+,u_-).
\]
By Lemma~\ref{lem:compactness2}, $\overline{\mc{M}}$ is compact.  Thus
the signed count of boundary points is
\[
0 = \#\partial\overline{\mc{M}} = \sum_{(u_+,u_-)}
-\#\partial\overline{V(u_+,u_-)}.
\]
To understand this sum, let $v_\pm$ denote the $\op{ind}=1$ component
of $u_\pm$.  Then
\begin{equation}
\label{eqn:v}
\begin{split}
v_+ & \in
\mc{M}_1^J(\alpha_+/\gamma_+, \beta/\gamma_+, Z_+-[\R\times\gamma_+])/\R,\\
 v_- &\in \mc{M}_1^J(\beta/\gamma_-, \alpha_-/\gamma_-,
 Z_--[\R\times\gamma_-])/\R,
\end{split}
\end{equation}
where $\gamma_\pm$ and $Z_\pm$ are as above.  Thus
\begin{equation}
\label{eqn:dM1}
\#\partial\overline{\mc{M}} = \sum_\beta  \;
\sum_{\substack{\gamma_+ | \alpha_+, \beta \\ \gamma_- | \beta,
\alpha_-}}
\;\sum_{Z_+,Z_-}  \;\sum_{\mbox{\scriptsize $v_+, v_-$ as in
\eqref{eqn:v}}}
-\#\partial\overline{V(u_+,u_-)}.
\end{equation}

By Lemma~\ref{lem:d22},
\begin{equation}
\label{eqn:dM2}
-\#\partial\overline{V(u_+,u_-)} = \left\{\begin{array}{cl} 0, &
\mbox{if $\beta$ is not admissible},\\
\epsilon(v_+)\epsilon(v_-), & \mbox{if $\beta$ is admissible.}
\end{array}\right.
\end{equation}
Let us clarify the signs and factorials here.  First, the signs
 $\epsilon(v_+)$ and $\epsilon(v_-)$ are well-defined when $\beta$ is
 admissible.  If $\beta$ is not admissible, then to apply
 Lemma~\ref{lem:d22} one needs to choose some orderings and asymptotic
 markings of the ends of $u_\pm$.  However,
 $-\#\partial\overline{V(u_+,u_-)}$ is defined independently of this
 choice.  Second, the factorials in \eqref{eqn:factorials} have
 disappeared in
\eqref{eqn:dM2}, because the count
$\#G(u_+,u_-)=-\#\partial\overline{U}$ distinguishes curves in
$\partial\overline{U}$ that have different asymptotic markings and
orderings of the ends but represent the same element of
$\partial\overline{V(u_+,u_-)}$.  More precisely, given $v\in \partial
\overline{V(u_+,u_-)}$, the corresponding curves
$u\in\partial\overline{U}$ differ from each other by the following operations:
\begin{itemize}
\item
changing the asymptotic marking of a positive (resp.\ negative) end of
$u$ that corresponds to an $\R$-invariant component of $u_+$ (resp.\
$u_-$).
\item
switching the ordering of two positive (resp.\ negative) ends of $u$
that correspond to identical $\R$-invariant components of $u_+$ (resp.\
$u_-$).
\end{itemize}
Since $\gamma_+$ and $\gamma_-$ have no common factor, it follows that
$\#(\partial\overline{U})$ equals $\#\partial\overline{V(u_+,u_-)}$ times the
product of factorials in equation \eqref{eqn:factorials}.

Since $\#\partial\overline{\mc{M}}=0$, equations \eqref{eqn:dM1}
and \eqref{eqn:dM2} imply \eqref{eqn:d2step1}.

{\em Step 2.\/} By definition, the coefficient of $\alpha_-$ in
$\partial^2\alpha_+$ is given by
\[
\begin{split}
\langle \partial^2\alpha_+, \alpha_-\rangle &= \sum_{\mbox{\scriptsize
$\beta$ admissible}} \langle
\partial\alpha_+, \beta\rangle \langle \partial\beta,
\alpha_-\rangle\\
&= \sum_{\mbox{\scriptsize $\beta$ admiss.}}\;
\sum_{\substack{\gamma_+ | \beta,\alpha_+ \\ \gamma_- | \beta,\alpha_-}}
\#\frac{\mc{M}_1^J(\alpha_+/\gamma_+, \beta/\gamma_+)}{\R} \cdot
\#\frac{\mc{M}_1^J(\beta/\gamma_-, \alpha_-/\gamma_-)}{\R}.
\end{split}
\]
In the second line, $\gamma_+$ and $\gamma_-$ are elliptic orbit sets.
We are also implicitly summing over relative homology classes with
$I=1$, which are suppressed here in order to simplify the notation.
To process the above sum, let $\gamma_0$ denote the greatest common
divisor of $\gamma_+$ and $\gamma_-$.  Then after dividing $\gamma_+$
and $\gamma_-$ by $\gamma_0$, the above sum becomes
\[
\sum_{\mbox{\scriptsize $\beta$ admissible}}\;
\sum_{\substack{\gamma_0|\alpha_+,\beta,\alpha_-\\
\gamma_0\gamma_+ | \beta,\alpha_+\\
\gamma_0\gamma_- | \beta,\alpha_-}}
\#\frac{\mc{M}_1^J(\alpha_+/\gamma_0\gamma_+,
\beta/\gamma_0\gamma_+)}{\R}
 \cdot
\#\frac{\mc{M}_1^J(\beta/\gamma_0\gamma_-, \alpha_-/\gamma_0\gamma_-)}{\R}.
\]
Here $\gamma_0$, $\gamma_+$, and $\gamma_-$ are elliptic orbit sets
such that $\gamma_+$ and $\gamma_-$ have no common factor.  Now we can
sum over $\gamma_0$ first and divide $\beta$ by $\gamma_0$ to obtain
\[
\sum_{\gamma_0|\alpha_+,\alpha_-}
\sum_{\mbox{\scriptsize $\beta$
admiss.}}\;
\sum_{\substack{\gamma_+ | \beta,\alpha_+/\gamma_0 \\
\gamma_- | \beta, \alpha_-/\gamma_0}}
\#\frac{\mc{M}_1^J(\alpha_+/\gamma_0\gamma_+,
\beta/\gamma_+)}{\R}
 \cdot
\#\frac{\mc{M}_1^J(\beta/\gamma_-, \alpha_-/\gamma_0\gamma_-)}{\R}.
\]
Again, $\gamma_0$, $\gamma_+$, and $\gamma_-$ are elliptic orbit sets
such that $\gamma_+$ and $\gamma_-$ have no common factor.  For each
$\gamma_0$, by equation \eqref{eqn:d2step1} applied to
$\alpha_\pm/\gamma_0$, the above sum over $\beta$ equals zero.  This
completes the proof that $\partial^2=0$.
\end{proof}

\begin{remark}
There is also a ``twisted'' version of ECH, with coefficients in the
group ring over $H_2(Y)$ (or a quotient thereof), which keeps track of
the relative homology classes of the $J$-holomorphic curves, see
\cite[\S11.2]{t3}.  The same argument with a bit more notation shows
that $\partial^2=0$ for the twisted chain complex as well.
\end{remark}

\subsection{Calculation of ECH gluing coefficients, first half}
\label{sec:GCC1}

To prepare for the proof of Proposition~\ref{prop:PTC}, we now
establish a special case:

\begin{proposition}
\label{prop:pinpout}
For any irrational number $\theta$ and positive integer $M$,
\[
c_\theta(\pin_\theta(M) \mid \pout_\theta(M)) = 1.
\]
\end{proposition}

The proof of Proposition~\ref{prop:pinpout} uses induction on $M$.
The key to carrying out the induction is the following lemma.

\begin{lemma}
\label{lem:indkey}
Write
$\pin_\theta(M)=(a_1,\ldots,a_{k})$
and
$\pout_\theta(M)=(b_1,\ldots,b_{l})$,
with the standard ordering convention \eqref{eqn:OC}.  Then:
\begin{description}
\item{(a)} There is a unique subset
$
I=\{i_1<\cdots <i_m\}\subset \{1,\ldots,l\}
$
such that
\begin{equation}
\label{eqn:ICondition'}
\sum_{j=1}^{m-1} b_{i_j} < a_1 \le \sum_{j=1}^m b_{i_j},
\end{equation}
and moreover $I=\{1,\ldots,m\}$ for some $m$.
\item{(b)}
With $m$ as above, if $1\le n\le m$, then
\begin{equation}
\label{eqn:step2}
\delta_\theta\left(a_1-\sum_{j=1}^{n-1}b_j\, ; \, b_n\right) = 1.
\end{equation}
\item{(c)} $\pin_\theta(M-a_1) = (a_2,\ldots,a_k)$ with the standard
ordering convention.
\item{(d)}
Let $\overline{b} \eqdef
\sum_{j=1}^m b_j - a_1$.  Then with the standard ordering convention,
\begin{equation}
\label{eqn:wnst}
\pout_\theta(M-a_1) = \left(\overline{b},
b_{m+1},\ldots,b_l\right).
\end{equation}
\end{description}
\end{lemma}

\begin{proof}
We begin with some preliminary remarks.  Note that with the ordering
convention \eqref{eqn:OC}, the lattice points on the path
$\lin_\theta(M)$ are the points
\[
\sum_{i=1}^n \bpm a_i \\ \ceil{a_i\theta} \epm,\quad
n=0,\ldots,k,
\]
while the lattice points on the path $\lout_\theta(M)$ are
\[
\sum_{j=1}^n \bpm b_j \\ \floor{b_j\theta} \epm,\quad
n=0,\ldots,l.
\]
Let $\Delta_\theta(M)$ denote the open region in the plane consisting
of points $(x,y)$ that (i) have $0\le x \le M$, (ii) are strictly
below the path $\lin_\theta(M)$, and (iii) are strictly above the path
$\lout_\theta(M)$.  A key observation, which we will use repeatedly
below, is that by construction {\em the region $\Delta_\theta(M)$
contains no lattice points\/}.
Note also that by Lemma~\ref{lem:pkappa}, we have
\begin{equation}
\label{eqn:kappapinpout}
\kappa_\theta(a_1,\ldots,a_k\mid b_1,\ldots,b_l)=1.
\end{equation}

Proof of (a): It will suffice to show that for each $n=1,2,\ldots,$
\begin{equation}
\label{eqn:ISufficient}
\sum_{j=1}^{n-1} b_j < a_1 \Longrightarrow M-b_n < a_1.
\end{equation}
To see that
\eqref{eqn:ISufficient} suffices, suppose that $I$ satisfies
\eqref{eqn:ICondition'}, and let $n$ be the smallest positive integer
that is not in $I$.  Suppose to get a contradiction that $I$ contains
an integer larger than $n$.  Then the first inequality in
\eqref{eqn:ICondition'} implies that $\sum_{j=1}^{n-1} b_j < a_1$.  So
by \eqref{eqn:ISufficient} we have $M-b_n < a_1$.  Since $n\notin I$,
the second inequality in \eqref{eqn:ICondition'} is then impossible.

To prove \eqref{eqn:ISufficient}, suppose to the contrary that
\begin{equation}
\label{eqn:contrary}
\sum_{j=1}^{n-1}b_j < a_1,\quad\quad a_1\le M-b_n.
\end{equation}
Consider the lattice point in the plane
\[
\bpm x \\ y \epm
 \eqdef \bpm a_1 \\ \ceil{a_1\theta} \epm +
 \bpm b_n \\ \floor{b_n\theta} \epm.
\]
To get a contradiction we will show that
$(x,y)\in\Delta_\theta(M)$.

(i) To start, the second inequality in \eqref{eqn:contrary} implies
that $x\le M$.

(ii) Next, $(x,y)$ is strictly below the path
$\lin_\theta(M)$, because the vector $(a_1,\ceil{a_1\theta})$ is
on the path $\lin_\theta(M)$, while the vector
$(b_n,\floor{b_n\theta})$ points to the right and has slope less than
that of all subsequent edges on the path $\lin_\theta(M)$.
Indeed, $(b_n,\floor{b_n\theta})$ has slope less than $\theta$, while
all of the edges in the path $\lin_\theta(M)$ have slope greater
than $\theta$.

(iii) To see that $(x,y)$ is strictly above the path
$\lout_\theta(M)$, rewrite $(x,y)$ as a sum of two vectors as
follows:
\[
\bpm x \\ y \epm = \sum_{j=1}^n \bpm b_j \\ \floor{b_j\theta} \epm
+
\left(\bpm a_1 \\ \ceil{a_1\theta} \epm - \sum_{j=1}^{n-1} \bpm b_j \\
\floor{b_j\theta} \epm\right).
\]
Then the first vector is on the path $\lout_\theta(M)$, while the
second vector points to the right (by the first inequality in
\eqref{eqn:contrary}), and has slope greater than that of all
subsequent edges in the path $\lout_\theta(M)$ (because it has
slope greater than $\ceil{a_1\theta}/a_1>\theta$).

Proof of (b): For $n=1,\ldots,m$, let $T_n$ denote the triangle with
vertices
\[
\sum_{j=1}^{n-1} \bpm b_j \\ \ceil{b_j\theta}\epm,\quad \sum_{j=1}^n
\bpm b_j \\ \ceil{b_j\theta}\epm, \quad \bpm a_1 \\ \ceil{a_1\theta} \epm.
\]
Then the interior of $T_n$ is in $\Delta_\theta(M)$, and
hence contains no lattice points, and the interiors of the edges of
$T_n$ also contain no lattice points, by the definition of the
incoming and outgoing partitions.  It follows that $T_n$ has area
$1/2$, i.e.\
\begin{equation}
\label{eqn:area1/2}
\det\begin{pmatrix} b_n & a_1-\sum_{j=1}^{n-1} b_j \\
\floor{b_n\theta} & \ceil{a_1\theta} - \sum_{j=1}^{n-1}
\floor{b_j\theta}
\end{pmatrix}
= 1.
\end{equation}

Next, in the notation of Definition~\ref{def:indTheta}, we have
\[
\begin{split}
\op{ind}_\theta(a_1,\ldots,a_k \mid b_1,\ldots,b_l) & = 
\op{ind}_\theta\left(a_1 \mid
b_1,\ldots,b_{n-1},a_1-\sum_{j=1}^{n-1}b_j\right) + \\
& \;\; \;\; +
\op{ind}_\theta\left(a_1-\sum_{j=1}^{n-1}b_j, a_2,\ldots,a_k
 \mid b_n,\ldots,b_l\right).
\end{split}
\]
We know by Lemma~\ref{lem:pkappa} that the left side of this equation
equals zero, and the two terms on the right are nonnegative.  In
particular, the first term on the right must equal zero, so
\begin{equation}
\label{eqn:ceil}
\ceil{a_1\theta} -
\sum_{j=1}^{n-1}\floor{b_j\theta}
=
\ceil{\left(a_1-\sum_{j=1}^{n-1}b_j\right)\theta}.
\end{equation}
Equations \eqref{eqn:area1/2} and \eqref{eqn:ceil} imply equation
\eqref{eqn:step2}.

Proof of (c): This follows immediately from the definition of the
incoming partition and the ordering convention \eqref{eqn:OC}.

Proof (d): We begin with some preliminary calculations.  Suppose that $k>1$.
Then $\overline{b}>0$ by \eqref{eqn:kappapinpout}.
Next observe that
\[
\begin{split}
\op{ind}_\theta(a_1,\ldots,a_k \mid b_1,\ldots,b_l) & = 
\op{ind}_\theta\left(a_2, \ldots, a_k \mid
\overline{b},b_{m+1}\ldots,b_l\right) + \\
& \;\; \;\; +
\op{ind}_\theta\left(a_1, \overline{b} \mid b_1,\ldots,b_m\right).
\end{split}
\]
Similarly to \eqref{eqn:ceil}, the second term on the right must
vanish and so
\begin{equation}
\label{eqn:floor}
\floor{\overline{b}\theta}
=
\sum_{j=1}^m\floor{b_j\theta} - \ceil{a_1\theta}.
\end{equation}

We now prove \eqref{eqn:wnst} up to reordering.  Consider the
polygonal path $\Lambda$ whose initial vertex is
$(a_1,\ceil{a_1\theta})$, and whose subsequent vertices are the sums
$\sum_{j=1}^n(b_j,\floor{b_j\theta})$ for $n=m,\ldots,l$.  Note that
the interior of the initial edge of $\Lambda$ contains no lattice
points, because it is inside the region $\Delta_\theta(M)$.  It then
suffices to show that $\Lambda$ is the path $\lout_\theta(M-a_1)$
translated by $(a_1,\ceil{a_1\theta})$.

To prove this, first note that by equation \eqref{eqn:floor}, the
first edge of $\Lambda$ has slope
$\floor{\overline{b}\theta}/\overline{b}<\theta$, and hence all edges
of $\Lambda$ have slope less than $\theta$.  Also, by
\eqref{eqn:slopeCondition} the path $\Lambda$ is concave.
Second, the
final endpoint of the path $\Lambda$ is
\[
\bpm M \\ \floor{M\theta} \epm =
\bpm a_1 \\ \ceil{a_1\theta} \epm
 + \bpm M-a_1 \\ \floor{(M-a_1)\theta}\epm,
\]
by Lemma~\ref{lem:pkappa}(a) applied to $M$ and $M-a_1$ with the help of
part (c).  Third, there are no lattice points above the path $\Lambda$
and below the translate by $(a_1,\ceil{a_1\theta})$ of the line
$y=\theta x$, because any such lattice point would lie in the region
$\Delta_\theta(M)$.  This completes the proof of \eqref{eqn:wnst} up
to reordering.

To show that \eqref{eqn:wnst} respects the standard ordering
convention, it is enough to show that if $m<l$ then
\begin{equation}
\label{eqn:slopeCondition} 
\frac{\floor{\overline{b}\theta}}{\overline{b}}
\ge \frac{\floor{b_{m+1}\theta}}{b_{m+1}}.
\end{equation}
If $m<l$ and equation \eqref{eqn:slopeCondition} fails, consider the
lattice point
\begin{equation}
\label{eqn:xy1}
\bpm x \\ y \epm \eqdef \bpm a_1 \\ \ceil{a_1\theta} \epm + \bpm
b_{m+1} \\ \floor{b_{m+1}\theta} \epm.
\end{equation}
To get a contradiction, we will show that $(x,y)\in\Delta_\theta(M)$.
(i) First observe that $x=a_1+b_{m+1} < \sum_{j=1}^{m+1}b_j\le M$.
(ii) As in the proof of (a), it follows from \eqref{eqn:xy1} that
$(x,y)$ is strictly below the path $\lin_\theta(M)$.  (iii) By
\eqref{eqn:floor},
\begin{equation}
\label{eqn:xy2}
\bpm x \\ y \epm = \sum_{j=1}^{m+1} \bpm b_j \\ \floor{b_j\theta} \epm -
\bpm \overline{b} \\ \floor{\overline{b}\theta} \epm.
\end{equation}
By our assumption that \eqref{eqn:slopeCondition} fails, the vector
$(\overline{b},\floor{\overline{b}\theta})$ has strictly smaller slope
than the vectors $(b_j,\floor{b_j})$ for $j=1,\ldots,m+1$, so $(x,y)$
is strictly above the path $\lout_\theta(M)$.
\end{proof}

\begin{proof}[Proof of Proposition~\ref{prop:pinpout}.]
By equation \eqref{eqn:kappapinpout} and the definition of $c_\theta$, we have
\begin{equation}
\label{eqn:cthetaftheta}
c_\theta(\pin_\theta(M) \mid \pout_\theta(M)) =
f_\theta(a_1,\ldots,a_k\mid b_1,\ldots,b_l).
\end{equation}
By Lemma~\ref{lem:indkey}(a) and the
definition of $f_\theta$, we have
\[
\begin{split}
f_\theta(a_1,\ldots,a_k \mid b_1,\ldots,b_l) & =
f_\theta\left(a_2,\ldots,a_k \mid
\overline{b},b_{m+1},\ldots,b_l\right) \cdot \\ & \; \; \cdot
\prod_{n=1}^m \delta_\theta\left(a_1-\sum_{j=1}^{n-1}b_j\, , \, b_n\right).
\end{split}
\]
Then by Lemma~\ref{lem:indkey}(b),
\begin{equation}
\label{eqn:PET}
f_\theta(a_1,\ldots,a_k \mid b_1,\ldots,b_l) =
f_\theta\left(a_2,\ldots,a_k \mid
\overline{b},b_{m+1},\ldots,b_l\right).
\end{equation}
By Lemma~\ref{lem:indkey}(c),(d),
\begin{equation}
\label{eqn:PIS}
f_\theta\left(a_2,\ldots,a_k \mid
\overline{b},b_{m+1},\ldots,b_l\right)
=
c_\theta(\pin_\theta(M-a_1) \mid \pout_\theta(M-a_1)).
\end{equation}
Proposition~\ref{prop:pinpout} follows from \eqref{eqn:cthetaftheta},
\eqref{eqn:PET}, and \eqref{eqn:PIS} by induction on $k$.
\end{proof}


\subsection{Calculation of ECH gluing coefficients, second half}
\label{sec:GCC2}

We now prove Proposition~\ref{prop:PTC}.  We begin by clarifying the
hypothesis \eqref{eqn:star} in the statement of the proposition.  If $\Lambda_1$
and $\Lambda_2$ are two paths in the plane, let $\Lambda_1\Lambda_2$
denote the concatenated path that first traverses $\Lambda_1$ and then
traverses the appropriate translate of $\Lambda_2$.

\begin{lemma}
\label{lem:tfae}
For $0\le M'\le M$, the following are equivalent:
\begin{description}
\item{(a)}
$\pin_\theta(M'+n) = \pin_\theta(M') \cup \pin_\theta(n)\;\;$
for all $n=1,\ldots,M-M'$.
\item{(b)}
$
\ceil{(M'+n)\theta} = \ceil{M'\theta} + \ceil{n\theta} \quad\quad$ for all
$n=1,\ldots, M-M'$.
\item{(c)}
$\lin_\theta(M)=\lin_\theta(M')\lin_\theta(M-M')$.
\item{(d)}
Under the standard ordering convention, $\pin_\theta(M')$ is an
initial segment of $\pin_\theta(M)$.
\end{description}
\end{lemma}

\begin{proof}
(a) $\Rightarrow$ (b): For a given $n$, if $\pin_\theta(M'+n) =
\pin_\theta(M') \cup \pin_\theta(n)$, then it follows from
Lemma~\ref{lem:interpretPartition} that the edge vectors in the path
$\lin_\theta(M'+n)$ are the same as the edge vectors in the path
$\lin_\theta(M)\lin_\theta(n)$, possibly in a different order.
Therefore these two paths have the same endpoints, so
\[
\ceil{(M'+n)\theta} = \ceil{M'\theta} + \ceil{n\theta}.
\]

(b) $\Leftrightarrow$ (c): Observe that (b) is equivalent to:
\begin{description}
\item{(b$'$)}
There are no lattice points above the line $y=\theta x$ and below the line
$y-\ceil{M'\theta} = \theta(x-M')$ with $M'\le x \le M$.
\end{description}

By the interpretation of $\lin_\theta(M)$ as the boundary of a convex
hull, condition (c) is equivalent to the following two conditions: (i)
there are no lattice points below $\lin_\theta(M')\lin_\theta(M-M')$
and above the line $y=\theta x$ with $0\le x \le M$, and (ii) the edges in
the path $\lin_\theta(M')\lin_\theta(M-M')$ have monotonically
increasing slope.  By the definition of $\lin_\theta(M')$ and
$\lin_\theta(M-M')$, condition (b$'$) is equivalent to condition (i).
But condition (i) implies condition (ii).  To see this, note that to
prove (ii), it is enough to show that slope of the last edge in the
path $\lin_\theta(M')$ does not exceed the slope of the first edge in
the path $\lin_\theta(M-M')$.  If this fails, then the fourth vertex of
the parallelogram on these two edges is a lattice point of the type
ruled out by (i).

(b) $\Rightarrow$ (a): Since (b) implies (c), it follows by replacing
$M$ with $M'+n$ that (b) also implies $\lin_\theta(M'+n) =
\lin_\theta(M')\lin_\theta(n)$ for all $n=1,\ldots,M-M'$.  By
Lemma~\ref{lem:interpretPartition}, this implies (a).

(c) $\Leftrightarrow$ (d):
By Lemma~\ref{lem:interpretPartition} and the
convexity of $\lin_\theta(M)$, condition (d) holds if and only if
$\lin_\theta(M')$ is an initial subpath of $\lin_\theta(M)$.  But if
the latter holds, then the rest of the path $\lin_\theta(M)$ is by
definition $\lin_\theta(M-M')$.
\end{proof}

\begin{proof}[Proof of Proposition~\ref{prop:PTC}.]
We will use induction and Lemma~\ref{lem:indkey}.  By symmetry, we can
assume that $M_+ \le M_-$.  (Otherwise we can replace $\theta$ by
$-\theta$ and positive ends by negative ends.)  By
Proposition~\ref{prop:pinpout}, we may further assume that $M_+<M$.
We now proceed in four steps.

{\em Step 1.\/} We begin with some setup and preliminary calculations.
Write
\[
S = (a_1,\ldots,a_k;
a_1',\ldots,a_{k'}' \mid b_1,\ldots,b_l; b_1',\ldots,b_{l'}').
\]
Order the $a_i$'s and $b_j$'s according to the standard convention
\eqref{eqn:OC}, and order the $a_i'$'s and $b_j'$'s
so that
\begin{equation}
\label{eqn:OC'}
\frac{\floor{a_i'\theta}}{a_i'} \ge \frac{\floor{a_{i+1}'\theta}}{a_{i+1}'},
\quad \quad \frac{\ceil{b_j'\theta}}{b_j'} \le
\frac{\ceil{b_{j+1}'\theta}}{b_{j+1}'}.
\end{equation}

For future reference we now compute $\kappa_\theta(S)$.  By
Lemma~\ref{lem:pkappa},
\[
\begin{split}
\kappa_\theta(S) 
&= \ceil{M_+\theta} + \floor{(M-M_+)\theta} + k' - \floor{M_-\theta} -
\ceil{(M-M_-)\theta} + l'.
\end{split}
\]
Since $M_+<M_-$, by the hypothesis \eqref{eqn:star} and Lemma~\ref{lem:tfae}(b)
this becomes
\begin{equation}
\label{eqn:kappatheta}
\kappa_\theta(S) = \left\{\begin{array}{cl} 
k', & M_-=M,\\
k'+l'-1, & M_- < M.
\end{array}\right.
\end{equation}

Let $m$ denote the smallest integer such that $\sum_{j=1}^m b_j \ge
a_1$.  Observe that we must have a strict inequality $\sum_{j=1}^m b_j
> a_1$.  The reason is that since $M_+<M$, the hypothesis \eqref{eqn:star}
implies that $(a_1)$ is a proper subpartition of $\pin_\theta(M)$, while
$(b_1,\ldots,b_m)$ is a subpartition of $\pout_\theta(M)$.  If these
two subpartitions had the same size, then it would follow that
$\kappa_\theta(\pin_\theta(M) \mid
\pout_\theta(M))\ge 2$, contradicting Lemma~\ref{lem:pkappa}.

Next define $\overline{b} \eqdef \sum_{j=1}^m b_j - a_1$ and
\[
\overline{S} \eqdef (a_2,\ldots,a_k; \pout(M-M_+) \mid \overline{b},
b_{m+1},\ldots, b_l ; \pin_\theta(M-M_-)).
\]
By Lemma~\ref{lem:indkey}(c),(d),
\[
\overline{S} = (\pin(M_+-a_1);\pout(M-M_+) \mid
\pout(M_--a_1);\pin(M-M_-)).
\]
Moreover, the hypothesis \eqref{eqn:star} still holds when $(M,M_+,M_-)$ are
replaced by $(M-a_1,M_+-a_1,M_--a_1)$.  The strategy of the induction
will be to deduce the conclusions of the proposition for $S$ from
those for $\overline{S}$.

{\em Step 2.\/} We now show that if $J\subset\{1,\ldots,k\}$ and
$J'\subset\{1,\ldots,l'\}$ satisfy
\begin{equation}
\label{eqn:hypo}
\sum_{j\in J}b_j + \sum_{j\in J'}b_j'\ge a_1,
\end{equation}
then $\{1,\ldots,m\}\subset J$.

To prove this, first note that by 
the hypothesis \eqref{eqn:star} and Lemma~\ref{lem:tfae}(a),
\[
\pout_\theta(M) = (b_1,\ldots,b_l) \cup \pout_\theta(M-M_-).
\]
It follows that the $b_j$'s for $j\in J$, together with
$\pout_\theta(M-M_-)$, comprise a subpartition of $\pout_\theta(M)$.
By \eqref{eqn:hypo}, the sum of the numbers in this subpartition is
\begin{equation}
\label{eqn:badSubpartition}
\sum_{j\in J}b_j + (M-M_-)
\ge \sum_{j\in J}b_j + \sum_{j\in J'}b_j'
\ge a_1.
\end{equation}
By Lemma~\ref{lem:indkey}(a), this subpartition must contain the
minimal initial segment of $\pout_\theta(M)$ whose sum is at least
$a_1$.  By \eqref{eqn:star}, this initial segment is $(b_1,\ldots,b_m)$.

{\em Step 3.\/} We claim now that if $\{S_\nu\}$ is a
$\theta$-decomposition of $S$ (see Definition~\ref{def:TD}), reordered
so that $1\in I_1$ if $k>0$, then it must have the following
properties:
\begin{description}
\item{(i)}
If $k>0$, then $I_1=\{1,\ldots,k\}$; $I_1'=\{i\}$ for some $i$ with
$a'_i=a'_1$; $J_1=\{1,\ldots,q\}$ for some $q$; and $J_1'=\emptyset$.
\item{(ii)}
For all $\nu>1$ (and also for $\nu=1$ if $k=0$), we have
\[
|I_\vu|+|I'_\vu|=|J_\nu|+|J'_\nu|=1.
\]
\item{(iii)}
If $M_-<M$, then there exists $\nu$ such that $I_\nu=J_\nu=\emptyset$.
\end{description}
We prove this claim by induction on $k$.

(Base case.) Suppose that $k=0$ and let $\{S_\nu\}$ be a
$\theta$-decomposition of $S$.  Since $k=0$, the set $I_\nu$ is empty
for each $\nu$.  Since $\nu$ runs from $1$ to $\kappa_\theta(S)$, and
since $I'_\nu$ is nonempty for each $\nu$ by the sum condition
\eqref{eqn:sumCondition}, it follows that $\kappa_\theta(S)\le k'$.
We then deduce from equation
\eqref{eqn:kappatheta} that $l'\le 1$ and $\kappa_\theta(S)=k'$, so
$|I'_\nu|=1$ for each $\nu$.

By the hypothesis \eqref{eqn:star} and Lemma~\ref{lem:tfae}(c), we have
\[
\lout_\theta(M)=\lout_\theta(M_-)\lout_\theta(M-M_-).
\]
Therefore $l\le k'$, and $b_j=a_j'$ for all $j=1,\ldots,l$.  Recall
from \S\ref{sec:pinpout} that the ordering convention \eqref{eqn:OC'}
implies that $a_i'\ge a_{i+1}'$ for all $i$.  Now consider the $\nu$
for which $1\in J_\nu$.  Since $|I'_\nu|=1$, by the sum condition
\eqref{eqn:sumCondition} we must have $I'_\nu=\{i\}$ where
$a_i'=a_1'$, and therefore $J_\nu=\{1\}$ and $J'_\nu=\emptyset$.
Continuing by induction, the $\theta$-decomposition can be reordered
so that $J_\nu=\{\nu\}$ and $J'_\nu=\emptyset$ for $\nu=1,\ldots,l$.

If $l'=0$, then we have described all of $S_1,\ldots,S_\nu$.  If
$l'=1$, then the description of $\{S_\nu\}$ is completed by noting
that under the above reordering, $J_{l+1}=\emptyset$ and
$J'_{l+1}=\{1\}$.  Now points (i)--(iii) follow immediately from the
above description of $\{S_\nu\}$.

(Induction step.)  Suppose $k>0$ and assume that the claim holds for
$k-1$.  To carry out the induction we will relate
$\theta$-decompositions of $S$ to $\theta$-decompositions of
$\overline{S}$.  By equation \eqref{eqn:kappatheta},
$\kappa_\theta(S)=\kappa_\theta(\overline{S})$.
  Thus we can identify
a $\theta$-decomposition of $\overline{S}$ with a decomposition
\[
\begin{split}
\{2,\ldots,k\} &= \overline{I}_1\sqcup \cdots \sqcup
\overline{I}_{\kappa_\theta(S)},\\
\{1,\ldots,k'\} &= \overline{I}_1'\sqcup \cdots \sqcup
\overline{I}_{\kappa_\theta(S)}',\\
\{m,\ldots,l\} &= \overline{J}_1\sqcup \cdots \sqcup
\overline{J}_{\kappa_\theta(S)},\\
\{1,\ldots,l'\} &= \overline{J}_1'\sqcup \cdots \sqcup
\overline{J}_{\kappa_\theta(S)}',
\end{split}
\]
such that for each $\nu=1,\ldots,\kappa_\theta(S)$, the data set
$\overline{S}_\nu$ satisfies the sum condition
\eqref{eqn:sumCondition}.  Here $\overline{S}_\nu$ is defined as in
\eqref{eqn:Snu}, but with $b_m$ replaced by $\overline{b}$.

Given a $\theta$-decomposition $\{\overline{S}_\nu\}$ of
$\overline{S}$, reorder the $\theta$-decomposition so that $m\in
\overline{J_1}$.  We can then define a $\theta$-decomposition
$\{S_\nu\}$ of $S$ by setting
\[
I_1\eqdef \{1\}\cup\overline{I}_1, \quad\quad I'_1\eqdef \overline{I}'_1, \quad
\quad J_1\eqdef \{1,\ldots,m-1\}\cup\overline{J}_1, \quad \quad
J'_1\eqdef \overline{J}'_1
\]
and leaving the components of the $\theta$-decomposition for
$\nu=2,\ldots,\kappa_\theta(S)$ unchanged.  It follows from Step 2
that every $\theta$-decomposition of $S$ is obtained this way from a
$\theta$-decomposition of $\overline{S}$.  Points (i)--(iii) for
$\theta$-decompositions of $S$ then follow from points (i)--(iii) for
$\theta$-decompositions of $\overline{S}$.  Note that
Lemma~\ref{lem:indkey}(d) gurarantees that when $k=1$, the unique element
$i$ of $I_1'$ will satisfy $a_i'=a_1'$.

{\em Step 4.\/} We now complete the proof of the proposition.  
Part (a) is an immediate consequence of points (i)--(iii) from Step 3.
We now prove part (b) by induction on $k$. (Part (c) then follows by
symmetry.)

If $k=0$ then
\[
S=(;\pout_\theta(M)\mid \pout_\theta(M)).
\]
In this case $\kappa_\theta(S)=k'=l$, and a $\theta$-decomposition of
$S$ is equivalent to a permutation of $\pout_\theta(M)$ that preserves
the sizes of the elements.  So it follows immediately from the
definition of $c_\theta$ that $c_\theta(S)=\pout_\theta(M)!$ as
desired.

If $k>0$, then as in the proof of Proposition~\ref{prop:pinpout}, it
follows from assertion (i) of Step 3 and Lemma~\ref{lem:indkey}(c),(d)
that
\[
c_\theta(S) = c_\theta(\overline{S}) \cdot \prod_{n=1}^m
\delta_\theta\left(a_1-\sum_{j=1}^{n-1}b_j\, , \, b_n\right).
\]
By Lemma~\ref{lem:indkey}(b), this becomes
$c_\theta(S)=c_\theta(\overline{S})$.  We are now done by induction.
\end{proof}



\begin{thebibliography}{99}

\bibitem{dbn} D. Bar-Natan, {\em On the Vassiliev knot invariants\/},
  Topology {\bf 34} (1995), 423--472.

\bibitem{behwz} F. Bourgeois, Y. Eliashberg, H. Hofer, K. Wysocki, and
  E. Zehnder, {\em Compactness results in symplectic field theory\/},
  Geom. Topol. {\bf 7} (2003), 799--888.

\bibitem{bm} F. Bourgeois and K. Mohnke, {\em Coherent orientations in
symplectic field theory\/}, Math. Z. {\bf 248} (2004), 123--146.

\bibitem{dragnev} D. Dragnev, {\em Fredholm theory and transversality
  for noncompact pseudoholomorphic maps in symplectizations\/},
  Comm. Pure Appl. Math {\bf 57} (2004), 726--763.

\bibitem{egh} Y. Eliashberg, A. Givental, and H. Hofer, {\em
Introduction to symplectic field theory\/}, Geom. Funct. Anal. (2000),
560--673.

\bibitem{fabert} O. Fabert, {\em Counting trivial curves in rational
  symplectic field theory\/}, in preparation.

\bibitem{fh} A. Floer and H. Hofer, {\em Coherent orientations for
  periodic orbit problems in symplectic geometry\/}, Math. Z. {\bf
  212} (1993), 13--38.

\bibitem{hofer} H. Hofer, {\em Holomorphic curves and dynamics in
  dimension three\/}, Symplectic geometry and topology (Park City, UT,
  1997), 35--101, IAS/Park City Math. Ser. {\bf 7}, AMS, 1999.

\bibitem{hwz2} H. Hofer, K. Wysocki, E. Zehnder, {\em Properties
  of pseudo-holomorphic curves in symplectizations.  II.  Embedding
  controls and algebraic invariants\/}, Geom. Funct. Anal. {\bf 5\/}
  (1995), 270--328.

\bibitem{pfh2} M. Hutchings, {\em An index inequality for embedded
pseudoholomorphic curves in symplectizations\/},
J. Eur. Math. Soc. {\bf 4\/} (2002), 313--361.

\bibitem{pfh3} M. Hutchings and M. Sullivan, {\em The periodic Floer
  homology of a Dehn twist\/}, Algebr. Geom. Topol. {\bf 5\/} (2005),
  301--354.

\bibitem{t3} M. Hutchings and M. Sullivan, {\em Rounding corners of
  polygons and the embedded contact homology of $T^3$\/}, Geometry and
  Topology {\bf 10\/} (2006), 169--266.

\bibitem{obg2} M. Hutchings and C. H. Taubes, {\em Gluing
  pseudoholomorphic curves along branched covered cylinders II\/},
  arxiv:0705.2074.

\bibitem{krmr} P. Kronheimer and T. Mrowka, {\em Monopoles and
three-manifolds\/}, book in preparation.

\bibitem{ozsz} P. Ozsv\'{a}th and Z. Szab\'{o}, {\em Holomorphic disks and
  topological invariants for closed three-manifolds\/}, Ann. of
  Math. {\bf 159} (2004), 1027--1158.

\bibitem{schwarz} M. Schwarz, {\em Cohomology operations from $S^1$
cobordisms in Floer homology\/}, ETH Z\"{u}rich PhD thesis, 1995.

\bibitem{siefring} R. Siefring, {\em The relative asymptotic behavior
  of pseudoholomorphic half-cylinders\/}, math.SG/0702356.

\end{thebibliography}
\end{document}